\magnification=1000
\hsize=11.7cm
\vsize=18.9cm
\lineskip2pt \lineskiplimit2pt
\nopagenumbers

\hoffset=-1truein
\voffset=-1truein

\advance\voffset by 4truecm
\advance\hoffset by 4.5truecm

\newif\ifentete

\headline{\ifentete\ifodd	\count0 
      \rlap{\head}\hfill\tenrm\llap{\the\count0}\relax
    \else
        \tenrm\rlap{\the\count0}\hfill\llap{\head} \relax
    \fi\else
\global\entetetrue\fi}

\def\entete#1{\entetefalse\gdef\head{#1}}
\entete{}

\input amssym.def
\input amssym.tex

\def\-{\hbox{-}}
\def\.{{\cdot}}
\def\O{{\cal O}}
\def\K{{\cal K}}

\def\L{{\cal L}}

\def\P{{\cal P}}

\def\B{{\cal B}}
\def\W{{\cal W}}

\def\int{\frak i\frak n\frak t}

\def\qq{\quad{\rm and}\quad}

\def\too{\longrightarrow}

 3
 2
\font\large=cmr10  scaled \magstep 2
 2
\font\larti=cmti10  scaled \magstep 2
 2
\font\cds=cmr7
\font\cdt=cmti7

\font\it=cmti10

\centerline{\large Parameterization of irreducible characters}
\medskip
\centerline{\large  for {\larti p}-solvable groups}

\vskip 0.5cm

\centerline{\bf Lluis Puig }
\medskip
\centerline{\cds CNRS, Institut de Math\'ematiques de Jussieu}
\smallskip
\centerline{\cds 6 Av Bizet, 94340 Joinville-le-Pont, France}
\smallskip
\centerline{\cds puig@math.jussieu.fr}

\vskip 0.5cm
\noindent
{\bf £1. Introduction}
\bigskip
£1.1. Let $p$ be a prime number, $k$ an algebraically closed field of charac-teristic $p$ and $G$ a $p\-$solvable finite group. In~[16], up to the choice of a {\it polarization $\omega$\/}, we give a {\it natural\/} parameterization of the set of isomorphism classes  of the simple $kG\-$modules, in terms of the set of $G\-$conjugacy classes of the {\it weights\/} of $G$ with respect to
$p\,,$ introduced by Jon Alperin in~[1] in order to formulate his celebrated conjecture.

\medskip
£1.2. In [4] Everett Dade formulates a refinement of Alperin's conjecture involving {\it ordinary\/} irreducible characters ---
with their {\it defect\/} --- instead of simple $kG\-$modules, and in [17] Geoffrey Robinson proves that the new conjecture holds
for $p\-$solvable groups. But this refinement is formulated in terms of a vanishing alternating sum of cardinals of suitable sets of
irreducible characters, {\it without\/} giving any possible refinement for the  {\it weights\/}.

\medskip
£1.3. Our purpose in this note is to show that, in the case of the $p\-$solvable finite groups, the method developed in [16]
can be suitably refined to provide,  up to the choice of a {\it polarization $\omega$\/},   a {\it natural bijection\/} --- namely compatible with the action of the group of {\it outer automorphisms\/}  of $G$ --- between the sets of {\it absolutely irreducible characters\/} of $G$ and of $G\-$conjugacy classes of suitable {\it inductive weights\/}  of $G$ defined in \S8 below. Moreover,
any {\it inductive weight\/} of $G$ admits an associated {\it block\/} of $G$ and a natural definition of a {\it defect\/} 
in such a way that then our bijection preserves the associated {\it blocks\/} and the  {\it defects\/}.

\medskip
£1.4. Let $\O$ be a complete discrete valuation ring of characteristic zero admitting $k$ as the {\it residue field\/} and denote by
$\K$ the {\it quotient field\/} of $\O\,.$ As in [16], the inductive arguments lead to consider a more general situation
including central extensions of the involved finite groups; but, as we will see bellow, in the present case the central 
$k^*\-$extensions considered in [16] have to be replaced by the {\it  central $\O^*\-$extensions\/}, analogously called
{\it $\O^*\-$groups\/}.

\medskip
£1.5. Contrarily  to the  $k^*\-$group case, a {\it block\/} of an {\it $\O^*\-$group\/} need not be isomorphic 
to a block of a finite group; although the arguments are very simi-lar, this forces us to develop to some extend 
the basic results on {\it $\O^*\-$groups\/}; in particular, our results on $p\-$solvable finite $k^*\-$groups in [15] have to
be extended to the  {\it $p\-$solvable finite $\O^*\-$groups.\/}
\eject

\bigskip
\noindent
{\bf £2. Generalities on $\O^*\-$groups}
\bigskip
£2.1. As above, let $p$ be a prime number and $\O$ a complete discrete valuation ring with a {\it quotient\/} field $\K$ of 
characteristic zero and an algebraically closed {\it residue} field $k$ of characteristic $p\,.$ An $\O^*\-${\it group\/} 
is a group $\hat G$ endowed with an injective group homomorphism $\varphi\, \colon {\O}^*\to Z(\hat G)$ 
where $Z(\hat G)$ denotes the center of $\hat G\,.$ As in the $k^*\-$group case [6,~\S5], we call
$\O^*\-${\it quotient\/} of $\hat G$ the quotient $G = \hat G/\hbox{Im}(\varphi)$ and write $\lambda.x$ instead of
$\varphi(\lambda)\,x$  for any $\lambda\in \O^*$ and any~$x\in \hat G\,;$ we say that $\hat G$ is {\it split\/} 
whenever $\hat G\cong \O^*\!\times G\,.$ Note that, the central product
$$\hat G^k = k^*\times_{\O^*} \hat G\cong \hat G/\varphi\bigl(1+J({\O})\bigr)
\eqno £2.1.1\phantom{.}$$ 
is a $k^*\-$group of $k^*\-$quotient $G\,.$ Similarly, if $\K'$ is  a finite field extension of $\K$ 
and~$\O'$ is the {\it ring of integers\/} of $\K'\,,$ the central product 
$$\hat G^{{\O}'} = {\O}'^*\times_{\O^*} \hat G
\eqno £2.1.2\phantom{.}$$ 
has an evident structure of $\cal
O'^*\-$group with the same $\O'^*\-$quotient $G\,;$ moreover, since the residue
field of $\O'$ coincides with $k\,,$ the degree of $\K'$ over~$\K$
coincides with the corresponding {\it ramification index\/} of $\O'$ over $\O\,.$

\medskip
£2.2.  The $\O^*\-$group algebra of an $\O^*\-$group $\hat G$ is the $\O\-$algebra
$${\O}_*\hat G ={\O}\otimes_{\O \O^*}{\O} \hat G
\eqno £2.2.1\phantom{.}$$
where ${\O \O^*}$ and ${\O} \hat G$ denote the respective ordinary group algebras of $\O^*$ 
and~$\hat G$ over $\O\,;$ we also set $\K_*\hat G = \K\otimes_\O \O_*\hat G\,.$
The introduction of $\O^*\-$groups mainly depends on the well-known situation exhibited by the 
next proposition where $\hat G$ need not be {\it split\/} even if $\hat H$ is  so. We use  terminology 
introduced in~£2.6 below.

\bigskip
\noindent
{\bf Proposition £2.3.} {\it Let $G$ be a finite group, $H$ a normal subgroup of $G\,,$ 
$\hat H$ an $\O^*\-$group with $\O^*\-$quotient $H\,.$ Assume that the action of $G$ on $H$
can be lifted to an action of $G$ on $\hat H$ and that $G$ stabilizes an absolutely irreducible 
charac-ter $\zeta$ of $\hat H\,.$ Denote by $\,e_\zeta$ the primitive idempotent of $Z(\K_*H)$
associated with~$\zeta\,.$  If $\K'$ is a finite field extension of $\K$ of degree divisible by~$\big\vert G/[G,G]\big\vert\,,$
denoting by $\O'$ the ring of integers of $\K'$ there are an $\O'^*$-group $\hat G$ with 
$\O'^*$-quotient~$G$ containing and normalizing~$\hat H^{{\O}'}\,,$ and an $\O'^*$-group homomorphism 
$$\hat G\too (\K'_*\hat H^{{\O}'}\!e_\zeta)^*
\eqno £2.3.1\phantom{.}$$
lifting  the action of $G$ on $\K'_*\hat H^{{\O}'}\!e_\zeta$ and extending the structural $\O'^*$-group 
homomorphism from $\hat H^{{\O}'}\,,$ and they are unique up to a unique
$\O'^*$-group isomorphism inducing the identity on $G$ .\/}

\medskip
\noindent
{\bf Proof:} We are assuming that $G$ acts on the $\K\-$algebra $\K_*\hat He_\zeta$ which is a full
matrix algebra over $\K$ and therefore the {\it pull-back\/}
$$\matrix{G& \too &{\rm Aut}(\K_*\hat He_\zeta)\cr
\uparrow&\phantom{\big\uparrow}&\uparrow\cr
{\hat G}^\K&\longrightarrow &(\K_*\hat He_\zeta)^*\cr}
\eqno £2.3.2\phantom{.}$$
determines a central $\K^*\-$extension ${\hat G}^{\K}$ of $G\,;$ note that the inclusion 
$H\subset G$ and~the structural $\O^*\-$group homomorphism $\hat H\to (\K_*\hat He_\xi)^*$ 
determine a $G\-$stable ``inclusion'' $\hat H\i {\hat G}^\K\,.$

\smallskip
On the other hand, the valuation maps $\vartheta\,\colon {\K^*}\to {\Bbb Z}$ and 
$\vartheta'\colon {\K'^*}\to {\Bbb Z}$ determine a commutative diagram of exact sequences
$$\matrix{1\to &{\O^*}&\too &{\K^*}&\buildrel \vartheta\over\too &{\Bbb Z}&\to 1\cr
{}&\cap&{}&\cap&{}&\downarrow&{}\cr
1\to&{\cal O'^*}&\too &{\cal K'^*}&\buildrel \vartheta'\over\too &{\Bbb  Z}&\to 1\cr}
\eqno £2.3.3\phantom{.}$$
where the vertical right-hand arrow is the multiplication by the degree $[\K/\K']$ of the extension; 
hence, since ${\Bbb H}^1(G,{\Bbb Z}) = \{0\}\,,$ they determine the following commutative diagram 
of exact sequences [3,~Chap.~XII,~\S2]
$$\matrix{0\to & {\Bbb H}^2(G,{\O}^*)&\to & {\Bbb H}^2(G,{\K^*})&\to& {\Bbb H}^2(G,{\Bbb Z})\cr
{}&\downarrow &{}&\downarrow &{}&\downarrow &{}\cr
0\to& {\Bbb H}^2(G,{\O}'^*)&\to& {\Bbb H}^2(G,{\K'^*})&\to& {\Bbb H}^2(G,{\Bbb Z})}
\eqno £2.3.4\phantom{.}$$
where the vertical right-hand arrow is the multiplication by $[\K/\K']\,;$ hence, since we have [9,~2.6.2] 
$${\Bbb H}^2(G,{\Bbb Z})\cong {\rm Ext}_{\Bbb Z}^1\big(G/[G,G],{\Bbb Z}\big)
\eqno £2.3.5,$$
 it follows from our hypothesis that this vertical arrow is zero.

 \smallskip
Consequently, the element of ${\Bbb H}^2(G,{\K'^*})$ determined by the central pro-duct  
${\K'^*}\times_{\K^*}{\hat G}^{\K}$ has a trivial image in ${\Bbb H}^2(G,{\Bbb Z})$ or, equivalently,
there is a surjective group homomorphism 
$$\theta :  {\K'^*}\times_{\K^*}{\hat G}^\K\too {\Bbb Z}
\eqno £2.3.6\phantom{.}$$
 fulfillng $\theta\big((\overline{{\O'}^*,1})\big) = \{0\}\,;$ moreover, this homomorphism is unique since 
we still have ${\Bbb H}^1(G,{\Bbb Z}) = \{0\}\,.$ Thus, $\hat G = {\rm Ker}(\theta)$ is an $\O'^*$-group 
with $\O'^*$-quotient $G\,,$ $\hat G$ contains and normalizes $\hat H^{\O'}$ since ${\hat G}^\K$ contains
and normalizes~$\hat H\,,$ and it is easily checked that the bottom arrow in diagram £2.3.2 induces an $\O'^*$-group
homomorphism $\hat G\too (\K'_*\hat H^{{\O}'}\!e_\zeta)^*$ which fulfills all the requirements.

\medskip
£2.4. Let $\hat G$ be an $\O^*$-group with a {\it finite} $\O^*\-$quotient $G\,;$ whereas any {\it finite\/}
 $k^*$-group contains a finite subgroup covering its own $k^*$-quotient [6, Lemma~5.5], $\hat G$ 
need not contain a finite subgroup covering $G\,;$ but,  if $G'$ and~$G''$
are two finite subgroups of $\hat G$ covering $G\,,$ they normalize
each other and therefore $G'.G''$ is also a finite subgroup of $\hat
G$ covering $G$ ; in particular, if $\hat G$ contains finite subgroups covering~$G\,,$ ${\rm Aut}(\hat G)$ stabilizes one
of them.

\bigskip
\noindent
{\bf Proposition £2.5.} {\it For some finite extension $\K'$ of $\K\,,$ denoting by $\O'$ the ring of integers of
$\K'\,,$ the $\cal O'^*$-group $\hat G^{{\cal O}'}$ contains a finite subgroup covering~$G$~.\/}
\medskip
\noindent
{\bf Proof:} Denote by $M$ the group of {\it Schur multipliers\/} of $G$ and by $A$ the maximal Abelian 
quotient of $G\,;$ since we have the exact sequence [9,~2.6.2]
$$0\to {\rm Ext}_{\Bbb Z}^1(A,\O^*)\too {\Bbb H}^2(G,\O^*)\too {\rm Hom}(M,\O^*)\to 0
\eqno £2.5.1,$$
the  element  $h$ of ${\Bbb H}^2(G,{\O}^*)$ corresponding to $\hat G$ induces a
group homomorphism $\sigma\,\colon M\to {\O}^*$ which then determines an evident commutative 
diagram of exact sequences
$$\matrix{ 0\to &{\rm Ext}_{\Bbb Z}^1(A,M)\hskip5pt\to &{\Bbb H}^2(G,M)\hskip5pt\to& {\rm Hom}(M,M)&\to 0\cr
&\downarrow &\downarrow &\downarrow &\cr
0\to &{\rm Ext}_{\Bbb Z}^1(A,{\O}^*)\hskip5pt\to &{\Bbb H}^2(G,{\O}^*)\hskip5pt\to&{\rm Hom}(M,{\O}^*)&\to 0\cr}
\eqno £2.5.2.$$
Then, an element of ${\Bbb H}^2(G,M)$ lifting ${\rm id}_{M}$ determines an $M\-$extension $E$ of~$G$ 
(which is finite!) and the central product ${\O}^*\times_{M} E$ determines $h'\in {\Bbb H}^2(G,{\O}^*)$ 
lifting $\sigma$ , so that the difference  $h'-h$ comes from ${\rm Ext}_{\Bbb Z}^1(A,{\O}^*)\,.$

\smallskip
That is to say, $h'-h$ determines an Abelian $\O^*\-$extension $Z$ of $A$ and, choosing a 
finite subset $X$ of $Z$ such that 
$$A = \prod_{x\in X} <\bar x>
\eqno £2.5.3\phantom{.}$$
 where $\bar x$ is the image of  $x\in X$ in~$A\,,$ we consider the field extension $\K'$ of $\K$ generated 
 by  $\vert <\bar x>\vert\-$th roots $\lambda_x$ of the elements $x^{\vert <\bar x>\vert}$ of $\O^*$ 
 when $x$ runs over $X\,;$ then,  denoting by $\O'$  the ring of integers of $\K'\,,$ the $\O'^*\-$group 
 $Z^{\O'}$ is {\it split\/} since the subgroup 
$$A' = \prod_{x\in X} <\lambda_x^{-1}.x>
\eqno £2.5.4\phantom{.}$$
 of $Z$ provides an $\O'^*\-$section. Consequently, denoting by $\iota\,\colon {\O^*}\to {\O'^*}$ 
the inclusion map, we get the equality ${\Bbb H}^2(G,\iota)(h'-h) = 0$
or, equivalently, an $\cal O'^*$-group isomorphism 
$${\O}'^*\times_{M} E\cong {\hat G}^{{\O}'}
\eqno £2.5.5\phantom{.}$$
inducing the identity on $G\,;$ now, the image of $E$ in $\hat
G^{{\cal O}'}$ is the announced finite subgroup.
\eject

 \medskip
 £2.6. If $\hat G$ and $\hat G'$ are two $\O^*\-$groups, we denote by $\hat G\,\hat\times\, \hat G'$ the quotient
 of the direct product $\hat G\times \hat G'$ by the image in $\hat G\times \hat G'$ of the {\it inverse\/}
 diagonal of~$\O^*\times \O^*\,,$ which has an obvious structure of $\O^*\-$group with 
 $\O^*\-$quotient $G\times G'\,;$ moreover, if $G = G'$ then we denote by $\hat G \star\hat G'$ 
 the $\O^*\-$group obtained from  the inverse image of $\Delta (G)\i G\times G$ in 
 $\hat G\,\hat\times\, \hat G'\,,$ which is nothing but the so-called {\it sum\/} of both central 
 $\O^*\-$extensions of $G\,;$  in particular, we have a {\it canonical\/} 
 $k^*\-$group isomorphism
 $$\hat G \star \hat G^\circ \cong \O^*\times G
 \eqno £2.6.1.$$
 Moreover, an $\O^*\-${\it group homomorphism\/} $\hat f\,\colon
\hat G\to \hat G'$ is a group homomorphism such that $\hat f(\lambda.x) =
\lambda.\hat f(x)$ for any  $\lambda\in \O^*$ and any $x\in \hat G$ ; it is
clear that $\hat f$ determines an ordinary group homomorphism
$f\,\colon G\to G'$ between the corresponding $\O^*\-$quotients.  
 
 \medskip
 £2.7. Note that for any $\O\-$free $\O\-$algebra $A$ of finite $\O\-$rank --- just called {\it $\O\-$algebra\/} in the sequel ---
 the group $A^*$ of invertible elements has a canonical $\O^*\-$group structure; we call {\it point\/} of $A$ 
 any $A^*\-$conjugacy class $\alpha$ of primitive idempotents of $A$ and denote by $A(\alpha)$ the simple
 quotient of $A$ determined by~$\alpha\,,$ and by $\P(A)$ the set of  {\it points\/} of $A\,.$
 If $S$ is a full matrix algebra over $O$ then ${\rm Aut}_\O (S)$ coincides with the $\O^*\-$quotient of $S^*\,;$
 in particular, any finite group $G$ acting on $S$ determines --- by {\it pull-back\/} --- an $\O^*\-$group $\hat G$
 of $k^*\-$quotient $G\,,$ together with an $\O^*\-$group homomorphism
 $$\rho : \hat G\too S^*
 \eqno £2.7.1.$$

\medskip
£2.8.   Let $\hat G$ be an $\O^*\-$group with finite $\O^*\-$quotient $G\,;$ a 
{\it $\hat G\-$interior~algebra\/} is an ${\O}\-$algebra $A$ endowed with an ${\O}^*\-$group
homomorphism $\hat G\to A^*\,;$ we say that $A$ is {\it primitive\/} whenever the unity element is 
primitive in $A^G\,.$ Note that the corresponding $k\-$algebra $k\otimes_{\O} A$ 
is only concerned by the $k^*\-$group $\hat G^k$ (cf.~£2.1.1).  A {\it $\hat G\-$interior algebra 
homomorphism\/} from $A$ to another $\hat G\-$interior algebra $A'$ is a {\it not necessarily unitary\/} 
algebra homomorphism $f\,\colon A\to A'$ fulfilling $f(\hat x\. a) = \hat x\. f(a)$ and $f(a\.\hat x) = f(a)\.\hat x\,;$ 
we say that $f$ is an {\it embedding\/}  whenever ${\rm Ker}(f) = \{0\}$ and ${\rm Im}(f) = f(1)A'f(1)\,.$ 
 Occasionally, it is handy to consider the 
 $(A'^{\hat G})^*\-$conjugacy class of $f$ that we denote by $\tilde f$ and call {\it exterior homomorphism\/}
 from $A$ to $A'\,;$ note that the {\it exterior homomorphisms\/} can be composed [7,~Definition~3.1].
Moreover, $A$ is a $G\-$algebra and all the {\it pointed group terminology\/} applies. If $H_\alpha$ is a
pointed group on the $G\-$algebra $A$ and $\hat H$ is the converse image of $H$ in $\hat G\,,$ 
we call $\hat H_\alpha$ a {\it pointed ${\O}^*\-$group\/} on $A$ and we set $A(\hat H_\alpha) = A(H_\alpha)$
denoting by $s_\alpha\colon A^H\to A(\hat H_\alpha)$ the canonical map. It is clear that $A_\alpha$ becomes an $\hat H\-$interior algebra and that we get a structural
$\hat H\-$interior algebra {\it exterior embedding\/}
$$\tilde f_\alpha : A_\alpha\longrightarrow \hbox{Res}_{\hat H}^{\hat G}(A)
\eqno £2.8.1.$$

\smallskip
£2.9. We similarly proceed with {\it inclusions\/} and {\it localness\/} but, contrarily  to the $k^*\-$group case, 
if $\hat P_\gamma$  is a local pointed $\O^*\-$group on $A$ --- namely, $P_\gamma$ is a local pointed group
on the $G\-$algebra $A$ and $\hat P$ is the converse image of $P$ in~$\hat G$ --- then $\hat P$ need not be split.  
 Recall that all the {\it maximal local pointed $\O^*\-$groups $\hat P_\gamma$
on $A$ contained in $\hat H_\alpha$\/} --- called {\it defect pointed $\O^*\-$groups of~$\hat H_\alpha$\/} --- are
mutually $H\-$conjugate  [7,~Theorem~1.2], and that the $\O\-$algebras $A_\alpha$ and $A_\gamma$ 
are mutually {\it Morita equivalent\/} [7,~Corollary~3.5]. As usual, we consider 
the {\it Brauer quotient\/} and the {\it Brauer algebra homomorphism\/}
$${\rm Br}_P^A : A^P\too A(P) = A^P\big/\sum_Q A_Q^P
\eqno £2.9.1\phantom{.}$$
where $Q$ runs over the set of proper subgroups of $P\,.$

\smallskip
£2.10. It is clear that $G$ acts on the set of pointed $\O^*\-$groups on $A$ and, if $\hat H_\alpha$ and 
$\hat K_\beta$ are two of them, we denote by $E_G(\hat K_\beta,\hat H_\alpha)$ the set of {\it $H\-$conjugacy
classes\/} of $\O^*\-$group homomorphisms from $\hat K$ to $\hat H$ induced by the
elements $x\in G$ such that $(\hat K_\beta)^x\subset \hat H_\alpha$ ; as usual, we set 
$$E_G(\hat H_\alpha) = E_G(\hat H_\alpha,\hat  H_\alpha)\cong N_{\hat G}(\hat H_\alpha)/\hat H.C_{\hat G}(\hat H)
\eqno £2.10.1.$$
If $\hat P_\gamma$ is a local pointed $\O^*\-$group on $A\,,$ as in~[9,~6.2] we denote by $\hat N_G(\hat P_\gamma)$ 
 the {\it $k^*\-$group\/} obtained from the {\it pull-back\/}
$$\matrix{N_G(\hat P_\gamma)&\longrightarrow
&\hbox{Aut}\bigl(A(\hat P_\gamma)\bigr)\cr
\uparrow&\phantom{\big\uparrow}&\uparrow\cr
\hat N_G(\hat P_\gamma)&\longrightarrow &A(\hat
P_\gamma)^*\cr}
\eqno £2.10.2,$$
 and then, since  the structural maps from $\hat P.C_{\hat G}(\hat P)$ to $N_G(\hat P_\gamma)$ and to $A(\hat
P_\gamma)^*$ determine a $N_G(\hat P_\gamma)\-$stable $k^*\-$group homomorphism
$\hat P^{^k}\!.C_{\hat G^{^k}}(\hat P)\to \hat N_G(\hat P_\gamma)\,,$ we still can define the $k^*\-$group
$$\hat E_G(\hat P_\gamma) = \big(\hat N_G(\hat P_\gamma) * N_{{\hat G}^{^k}}
(\hat P_\gamma)^\circ\big)\big/P.C_G(\hat P)
\eqno £2.10.3.$$

\medskip
£2.11. As in [5,~2.5],  we say that an injective $\O^*\-$group homomorphism 
$\hat\varphi\,\colon \hat K\to \hat H$ is an {\it $A\-$fusion from $\hat K_\beta$ to~$\hat H_\alpha$\/} 
whenever there is a  $\hat K\-$interior algebra~{\it embedding\/} 
$$f_{\hat\varphi} : A_\beta\too {\rm Res}_{\hat \varphi} (A_\alpha) 
\eqno £2.11.1\phantom{.}$$
such that the inclusion $A_\beta\i A$ and the composition of $f_{\hat\varphi}$
with the inclusion $A_\alpha\i A$ are $A^*\-$conjugate; we denote by $F_A(\hat K_\beta,\hat H_\alpha)$ 
the set of {\it $H\-$conjugacy classes\/} of  $A\-$fusions from $\hat K_\beta$ to~$\hat H_\alpha$
and, as usual, we write $F_A (\hat H_\alpha)$ instead of $F_A (\hat H_\alpha,\hat H_\alpha)\,.$ 
If $A_\alpha = iAi$ for $i\in \alpha\,,$ as in~[5, Corollary~2.13] we have a group homomorphism
$$F_A (\hat H_\alpha)\too N_{A_\alpha^{^*}}(\hat H\.i)\big/\hat H\.(A_\alpha^H)^*
\eqno £2.11.2\phantom{.}$$
and  it is clear that
$$E_G(\hat K_\beta,\hat H_\alpha)\subset F_A(\hat K_\beta,\hat H_\alpha)
\eqno £2.11.3.$$
\eject
\noindent
Note that if $B$ is another $\hat G\-$interior algebra and $h\,\colon A\to B$ a $\hat G\-$interior algebra 
embedding, we have
$$F_B(\hat K_\beta,\hat H_\alpha) = F_A(\hat K_\beta,\hat H_\alpha)
\eqno £2.11.4.$$

\medskip
£2.12. If $\hat P_\gamma$ is a local pointed $\O^*\-$group on $A$ , we have a canonical isomorphism $\hat P^k\cong k^*\times P\,;$ 
then, choosing $j\in \gamma\,,$ setting $A_\gamma = jAj$ and denoting by $\hat P^{^J}$ the converse image of $1\times P$ in $\hat P\,,$  $\hat F_A(\hat P_\gamma)$
denotes the {\it $k^*\-$group\/} obtained from the {\it pull-back\/} [9,~6.7]
$$\matrix{F_A(\hat P_\gamma)&\too
&N_{(A_\gamma)^{^*}}(\hat P\.j)/(A_\gamma^P)^*\cr
\uparrow&\phantom{\Big\uparrow}&\uparrow\quad\cr
\hat F_A(\hat P_\gamma)&\too &N_{(A_\gamma)^{^*}}(\hat P\.j)
/\hat P^{^J}\!\.\bigl(j+J(A_\gamma^P)\bigr)}
\eqno £2.12.1.$$
As in the ordinary case [9,~Proposition~6.12], the inclusion
$E_G(\hat P_\gamma)\subset F_A(\hat P_\gamma)$ can be lifted to a
{\it canonical $k^*\-$group homomorphism\/} 
$$\hat E_G(\hat P_\gamma)^\circ\too \hat F_A(\hat P_\gamma)
\eqno £2.12.2\phantom{.}$$
as we show in the next proposition.

\bigskip
\noindent
{\bf Proposition £2.13.} {\it Let $\hat P_\gamma$ be a local pointed $\O^*\-$group on $A\,,$
 choose $j\in \gamma$ and set $A_\gamma = jAj\,.$ For any $\hat x\in N_{\hat G}(\hat P_\gamma)$ and any
$\,a\in (A^P)^*$  having the same action  on the simple $k\-$algebra $A(\hat P_\gamma)\,,$ the
element $j(\hat x^{-1}.a)j$ is invertible  in~$A_\gamma$ and
normalizes $\hat P\.j\,.$ Moreover, denoting by $x\,,$ $\hat x^k$ and $\tilde x$ the respective images of $\hat x$ 
in $G\,,$ $\hat G^k$ and $E_G (\hat P_\gamma)\,,$ and setting
$$c_{\hat x, a} = j(\hat x^{-1}.a)j\.\hat P^{^J}\!\.\bigl(j+J(A_\gamma^P)\big)
\eqno £2.13.1,$$
there is a $k^*\-$group homomorphism
$$\hat E_G(\hat P_\gamma)^\circ\longrightarrow \hat F_A(\hat P_\gamma)
\eqno £2.13.2\phantom{.}$$
which maps the image of $\bigl(x,s_\gamma(a)\bigr)\otimes (\hat x^k)^{-1}$ in $\hat E_G(\hat P_\gamma)^\circ$ on
the element $(\tilde x^{-1},c_{\hat x, a})$ of~$\hat  F_A(\hat P_\gamma)\,.$\/}

\medskip
\noindent
{\bf Proof:} Since $a^{-1}\.\hat x$ acts trivially on $A(\hat P_\gamma)\,,$ the idempotents $j$ and $j^{\,a^{-1}\.\hat x}$ have the same image in
$A^P/J(A^P)$ and thus the image in this quotient of
$$b = jj^{\,a^{-1}\.\hat x}+(1-j)(1-j)^{\,a^{-1}\.\hat x}
\eqno £2.13.3\phantom{.}$$
is the unity element, so that $b$ is invertible  in $A^P\,;$ hence, since we have
$$b = \big(j(\hat x^{-1}\.a)j + (1- j)(\hat x^{-1}\.a)(1 -j)\big)(a^{-1}\.\hat x)
\eqno £2.13.4,$$
\eject
\noindent
 $j(\hat x^{-1}\.a)j$ is invertible in $A_\gamma\,;$ moreover, it is clear that, for any $\hat u\in \hat P\,,$ we have 
$$j(\hat x^{-1}\.a)j\.\hat u = j(\hat x^{-1}\hat u\.a)j =  \hat u^x\.j(\hat x^{-1}\.a)j
\eqno £2.13.5$$
so that $j(\hat x^{-1}\.a)j$ normalizes $\hat P\.j\,.$

\smallskip 
If $a'\in (A^P)^*$ has the same image as  $a$ in~$A(\hat P_\gamma)\,,$ and moreover we have 
$j^{\,a'^{-1}\.\hat x} = j\,,$ we get
$$\bigl(j(\hat x^{-1}\.a)j\bigr)\bigl(j(a'^{-1}\.\hat x)j\bigr) =
j(\hat x^{-1}\.aa'^{-1}\.\hat x)j = j(aa'^{-1})^xj
\eqno £2.13.6\phantom{.}$$
and this element  belongs to $j + J(A_\gamma^P)$ since $\,s_\gamma(aa'^{-1}) =  1\,;$ hence, the class 
$c_{\hat x, a}$ does not depend on our choice of $a$ and we may assume
that $j^{\,a^{-1}.\hat x} = j\,;$ in particular, we get the
announced map from $\hat E_G(\hat P_\gamma)^\circ$  to~$\hat F_A(\hat P_\gamma)\,.$ 
Finally, if $\hat x'\in N_{\hat G}(\hat P_\gamma)$ and $a'\in (A^P)^*$ have the same action on
$A(\hat P_\gamma)$ , we obtain
$$\eqalign{\big(j(\hat x^{-1}\.a)j\bigr)\bigl(j(\hat x'^{-1}.a')j\bigr) &= 
j(\hat x^{-1}\.a\.\hat x'^{-1}\.a')j\cr
& = j\bigl((\hat x'\hat x)^{-1}\.a'a^{\,\hat x'^{-1}\.a'}\big)j\cr}
\eqno £2.13.7\phantom{.}$$
and therefore, since $s_\gamma(a'a^{\,\hat x'^{-1}.a'}) = s_\gamma(aa')\,,$ it is easily checked that 
this map is a $k^*\-$group
homomorphism.

\medskip
£2.14. If $A'$ is another $\hat G\-$interior algebra and $f\,\colon A\to A'$ a $\hat G\-$interior algebra embedding, 
it follows from Proposition~£3.5 below that, as in [16,~2.9], denoting by $\gamma'$ the point of 
$P$ on $A'$ containing $f(\gamma)\,,$ we have a canonical $\O^*\-$group isomorphism
$$\hat F_{\tilde f} (P_\gamma) : \hat F_A (P_\gamma)\cong \hat F_{A'} (P_{\gamma'})
\eqno £2.14.1.$$
More precisely, let $Q_\delta$ be another local pointed group on $A$ and denote by $\delta'$
the point of $Q$ on $A'$ containing $f(\delta)$ and by
$$f_\gamma^{\gamma'} : A_\gamma\too A'_{\gamma'}\qq f_\delta^{\delta'} : A_\delta\too A'_{\delta'}
\eqno £2.14.2\phantom{.}$$
the induced embeddings; if there is  an {\it $A\-$fusion\/} $\varphi\,\colon Q\cong P$ 
from~$Q_\delta$ to $P_\gamma$ then, according to equality~£2.11.4 above, $\varphi$ is also  an 
{\it $A'\-$fusion\/} from $Q_{\delta'}$ to~$P_{\gamma'}\,,$
so that we have two $Q\-$interior algebra isomorphisms
$$f_\varphi : A_\delta\cong {\rm Res}_\varphi (A_\gamma)\qq f'_\varphi : A'_{\delta'}\cong {\rm Res}_\varphi (A'_{\gamma'})
\eqno £2.14.3\phantom{.}$$
and the uniqueness of the {\it exterior isomorphisms\/} $\tilde f_\varphi$ and $\tilde f'_\varphi$ forces 
the equality
$$\tilde f'_\varphi\circ \tilde f_\delta^{\delta'} = {\rm Res}_\varphi( \tilde f_\gamma^{\gamma'})\circ \tilde f_\varphi
\eqno £2.14.4\,.$$
In particular, since by our very definition we have 
$$\hat F_{{\rm Res}_\varphi (A_\gamma)}(Q_\delta) = \hat F_A (P_\gamma)\qq \hat F_{{\rm Res}_\varphi (A'_{\gamma'})}(Q_{\delta'}) = \hat F_{A'} (P_{\gamma'})
\eqno £2.14.5,$$
we get the following commutative diagram of $\O^*\-$group isomorphisms
$$\matrix{\hat F_A (Q_\delta) &\buildrel \hat F_{\tilde f_\varphi} (Q_\delta)\over \cong& \hat F_A (P_\gamma)\cr
\hskip-40pt {\scriptstyle \hat F_{\tilde f} (Q_\delta)}\hskip4pt\wr\!\Vert&\phantom{\Big\uparrow}
&\wr\Vert\hskip4pt{\scriptstyle \hat F_{\tilde f} (P_\gamma)}\hskip-40pt\cr
\hat F_{A'} (Q_{\delta'}) &\buildrel \hat F_{\tilde f'_\varphi} (Q_{\delta'})\over \cong& \hat F_{A' }(P_{\gamma'})\cr}
\eqno £2.14.6.$$

\medskip
£2.15. If $\hat H$ is an $\O^*\-$subgroup of $\hat G$ and $B$ an $\hat H\-$interior algebra, we call
{\it induced $\hat G\-$interior algebra of $B$\/} the
${\O}_*\hat G\-$bimodule
$$\hbox{Ind}_{\hat H}^{\hat G}(B) = {\O}_*\hat G\otimes_{{\O}_*\hat H} B\otimes_{{\O}_*\hat H} {\O}_*\hat G
\eqno £2.15.1\phantom{.}$$
endowed with the distributive multiplication defined by
$$(x\otimes b\otimes y)(x'\otimes b'\otimes y') =
\Bigg\{\matrix{x\otimes b\.yx'\.b'\otimes y'&\hbox{if}\;yx'\in \hat H\cr {}&{}\cr
0&\hbox{otherwise}\cr}
\eqno £2.15.2\phantom{.}$$
and with the $\O^*\-$group homomorphism $\hat G\longrightarrow
\hbox{Ind}_{\hat H}^{\hat G}(B)^*$ mapping $\hat x\in \hat G$ on the element
$\hat x.\hbox{Tr}_H^G(1\otimes 1\otimes 1)\,.$

\medskip
£2.16.   If $G^{\skew6\hat {}}\,$ is another $\O^*\-$group with  $\O^*\-$quotient $G$ and $A'$ is 
a $G^{\skew6\hat {}}\,\-$in-terior algebra, setting  $\,^{\skew4\hat {}}G =\hat G * G^{\skew6\hat {}}\,,$
the tensor product $A'' = A\otimes_{\O} A'$ has an evident structure of $\,^{\skew4\hat {}}G\-$interior
algebra. Then, for any pointed $\O^*\-$group $\,^{\skew4\hat {}}H_{\alpha''}$ on $A''$ , there are points $\alpha$
 of $\hat H$ on $A$ and $\alpha'$ of $H^{\skew6\hat {}}\,$ on~$A'$ , and an $\,^{\skew4\hat {}}H\-$interior algebra embedding
$A''_{\alpha''}\to A_\alpha\otimes_{\O} A'_{\alpha'}$ such that we
have a commutative diagram
$$\matrix{&{\rm Res}_{\,^{\skew4\hat {}}H}^{\,^{\skew4\hat {}}G}(A'')\hskip-10pt\cr
\nearrow\hskip-20pt&&\hskip-30pt\nwarrow\cr
A''_{\alpha''}&\too& A_\alpha\otimes_{\O} A'_{\alpha'}\cr}
\eqno £2.16.1.$$ 
 Moreover, if $\,^{\skew4\hat {}}P_{\gamma''}$ is a local
pointed $\O^*\-$group on $A''$ , the corresponding points $\gamma$
of $\hat P$ on $A$ and $\gamma'$ of $P^{\skew6\hat {}}\,$ on $A'$
have to be local.

\medskip
£2.17. As a matter of fact, we are mainly interested in the case where a Sylow $p\-$subgroup of $G$
 stabilizes an  $\O\-$basis of $A\,;$ then, $\gamma''$ is the unique
local point of $P$ on $A_\gamma\otimes_{\O}A'_{\gamma'}$ [7,~Proposition~5.6] and therefore, for any 
$\tilde\varphi\in F_A(\hat P_\gamma)$ and any  $\tilde\varphi'\in F_{A'}(P^{\skew6\hat {}}
\,_{\gamma'})$ having the same image in $\widetilde{\hbox{Aut}}(P)\,,$
the corresponding $\O^*\-$group {\it outer automorphism\/} $\tilde\varphi'' = \tilde\varphi * \tilde\varphi'$ of
$\,^{\skew4\hat {}}P$ belongs to $F_{A''}(\,^{\skew4\hat {}}P_{\gamma''})\,.$ In particular, 
if $F$ is  a finite group  and we have group homomorphisms
$$\theta : F\too F_A(\hat P_\gamma)\qq \theta' : F\too F_{A'}(P^{\skew6\hat {}}\,_{\gamma'})
\eqno £2.17.1\phantom{.}$$
 inducing the same group homomorphism $F\to \widetilde{\rm Aut}(P)\,,$ 
we get a group homomorphism $\theta''\,\colon F\to F_{A''}(\,^{ \skew4\hat {}}P_{\gamma''})$ and, setting
$$\eqalign{\widehat F^{^\gamma} &= {\rm Res}_{\,\theta}\big(\hat F_A(\hat P_\gamma)\big)\cr
\widehat F^{^{\gamma'}} = {\rm Res}_{\,\theta'}\big(\hat F_A(\hat P_{\gamma'})\big)
&\qq \widehat F^{^{\gamma''}} = {\rm Res}_{\,\theta''}\big(\hat F_A(\hat P_{\gamma''})\big)\cr}
\eqno £2.17.2,$$ 
the corresponding embedding in the bottom of diagram~2.16.1 induces a $k^*\-$group isomorphism
[7,~Proposition~5.11]
$$\widehat F^{^\gamma} * \widehat F^{^{\gamma'}}\cong\widehat F^{^{\gamma''}}
\eqno £2.17.3.$$

\bigskip
\noindent
{\bf £3. The $\O^*\-$group algebra}
\bigskip
£3.1. Let $\hat G$ be an $\O^*\-$group with finite $\O^*\-$quotient $G\,,$ and $\O'$ 
the ring of integers of a finite field extension $\K'$ of $\K$ such that the $\O'^*\-$group $\hat G^{{\O}'}$ 
contains a finite subgroup $G'$ covering~$G$ (cf.~Proposition~£2.5). It is clear that the product 
${\O'^*}\otimes \hat G$ in $\,{\O}'\otimes_{\O} {\O}_*\hat G\,$ can be identified to $\hat G^ {{\O}'}\!\,,$ so that we have 
${\O}'_*\hat G^{{\O}'} =  {\O}'\otimes_{\O} {\O}_*\hat G\,,$ and we will analyse the surjective 
$G\-$algebra homomorphism  ${\O}'G'\too {\O}'_*\hat G^{{\O}'}$ induced by the inclusion $G'\subset \hat G^{{\O}'}$ . 
We respectively denote by $Z'$ the Sylow $p\-$subgroup of $({\O}'^*.1)\cap G'$ and by $\Lambda'$ the converse image in
$\O'^*$ of the Hall $p'\-$subgroup of this intersection, and in the group algebra ${\O}'G'$ we set
$$e_{\Lambda'} = {1\over \vert \Lambda'\vert}\sum_{\lambda'\in \Lambda'}\lambda'^{-1}(\lambda'.1)
\eqno £3.1.1;$$
then, we still get a surjective $G\-$algebra homomorphism over $\O'$
$$s_{Z'} : {\O}'G'e_{\Lambda'}\longrightarrow {\O}'_*\hat G^{{\O}'}
\eqno £3.1.2\phantom{.}$$
which can be actually  considered as a homomorphism of  $G'\-$interior algebras; more precisely,
considering the obvious group homomorphism $Z'\to \O'^*\,,$ $s_{Z'}$~induces  a $G'\-$interior algebra isomorphism
$${\O}'\otimes_{{\O}'Z'}{\O}'G'e_{\Lambda'}\cong  {\O}'_*\hat G^{{\O}'}
\eqno £3.1.3.$$
Note that, if we were working over $k\,,$ only the $k^*\-$group $\hat G^k$ would be concerned
since we have
$$k\otimes_{\O} {\O}_*\hat G\cong k_*\hat G^k\cong
k\bar G'\bar e_{\Lambda'}
\eqno £3.1.4\phantom{.}$$
where we set  $\bar G' = G'/Z'$ and $\bar e_{\Lambda'}$ is the image of $e_{\Lambda'}$ in $k\bar G'\,.$
\eject

\bigskip
\noindent
{\bf Proposition £3.2.} {\it The radical $J({\O}'G'e')$ contains ${\rm Ker}(s_{Z'})\,,$  and for any subgroup $H$ of $G$ we have
$$s_{Z'}\big(({\O}'G'e_{\Lambda'})^H\big) = ({\O}'_*\hat G^{{\O}'})^H
\eqno £3.2.1.$$
In particular, $s_{Z'}$ is a strict covering homomorphism of $G\-$algebras. Moreover, for any $p\-$subgroup $P$ of $G\,,$ 
${\rm Br}_P^{{\O}_*\hat G}$ induces an $N_G(P)\-$algebra isomorphism
$$k_*C_{\hat G^{^k}}(\hat P)\cong ({\O}_*\hat G)(P)
\eqno £3.2.2\phantom{.}$$
where $\hat P$ denotes the converse image of $P$ in $\hat G\,.$\/}

 \medskip
 \noindent
 {\bf Proof:} Since $Z'$ is a finite $p\-$group, from isomorphism~£3.1.3 we get
$${\rm Ker}(s_{Z'})\subset J({\O}'Z').({\O}'G'e_{\Lambda'})
\subset J({\O}'G'e_{\Lambda'})
\eqno £3.2.3;$$
 moreover, since any subgroup $H$ of $G$ stabilizes the following direct sum
decomposition 
$${\O}'_*\hat G^{{\O}'} = \bigoplus_{x\in G} \,{\O}'x'
\eqno £3.2.4\phantom{.}$$ 
where $x'\in G'$ is a lifting of $x\in G\,,$ $({\O}'_*\hat G^{{\O}'})^H$ is linearly generated 
by the elements $Tr_{C_H(x')}^H(x')$  where $x'$ runs over the set of elements of $G'$ such that the
stabilizer of ${\O}'x'$ in $H$ centralizes $x'\,,$ and this  implies the equality~£3.2.1. 
Consequently, $s_{Z'}$ is a {\it strict covering homomorphism\/} of $G\-$algebras [10,~\S4]. In particular, if $H = P$ is a
$p\-$group and $\hat P$ is the converse image of $P$ in $\hat G$ , it
is easily checked that
$$\eqalign{({\O}_*\hat G)(P)&\cong ({\O}'_*\hat G^{{\O}'})(P)\cr
& = \bigoplus_{x\in C_G(\hat P)} k\hbox{Br}_P(x')\cong k_*C_{\hat G^{^k}}(\hat P)\cr}
\eqno £3.2.5\phantom{.}$$
where again $x'\i C_{G'}(\hat P)$ lifts $x\in C_G(\hat P)\,.$

\bigskip
\noindent
{\bf Corollary £3.3.} {\it The inclusion ${\O}_*\hat G\i {\O}'_*\hat G^{{\O}'}$ and the homomorhism 
$s_{Z'}$ induce a bijection, preserving inclusions and localness, between the set of
pointed $\O^*\-$groups $\hat H_\beta$ on ${\O}_*\hat G$ and
the set of pointed groups $H'_{\beta'}$ on ${\O}'G'e_{\Lambda'}$ such that
$H'$ contains $\Lambda'\.Z'\,.$ Moreover, if $\hat H_\beta$ and $H'_{\beta'}$
correspond to each other, then we have $\hat H^{{\O}'}\cap G' = H'$ and this equality together with
 $s_{Z'}$ determine the following group isomorphism and  $H'\-$interior algebra isomorphism
 $$\eqalign{E_{G'}(H'_{\beta'}) &\cong E_G(\hat H_\beta)\cr
  {\O}'\otimes_{{\O}'Z'}({\O}'G'e_{\Lambda'} )_{\beta'}&\cong {\O}'\otimes_{\O} ({\O}_*\hat G)_\beta\cr}
  \eqno £3.3.1.$$\/}
  \eject

\par
\noindent
{\bf Proof:} Since ${\O}'_*\hat G^{{\O}'}\cong \O'\otimes_\O \O_*\hat G\,,$ for any subgroup $H$ of $G$
we have 
$$({\O}'_*\hat G^{{\O}'})^H\cong \O'\otimes_\O (\O_*\hat G)^H
\eqno £3.3.2\phantom{.}$$
and therefore, since $k$ is algebraically closed, the inclusion ${\O}_*\hat G\subset {\O}'_*\hat G^{{\O}'}$ 
induces a bijection, preserving inclusions and localness, between the sets of corresponding pointed 
$\O^*\-$ and $\O'^*\-$groups, which maps any pointed $\O^*\-$group $\hat H_\beta$ on ${\O}_*\hat G$
on a pointed $\O'^*\-$group $(\hat H^{{\O}'})_{\beta^{{\O}'}}$ on ${\O}'_*\hat G^{{\O}'}$ fulfilling
$\beta\subset \beta^{{\O}'}\,;$ moreover, since $s_{Z'}$ is a {\it strict covering homomorphism\/} 
of $G\-$algebras, it does the same between the sets of pointed groups on ${\O}'G'e_{\Lambda'}$ 
and ${\O}'_*\hat G^{{\O}'}$  [10, Propositions~4.4 and~4.18], namely $H_{\beta^{{\O}'}}$ or,
equivalently, $(\hat H^{{\O}'})_{\beta^{{\O}'}}$ comes from the pointed group $H'_{\beta'}$ on
${\O}'G'e_{\Lambda'}$ fulfilling $s_{Z'}(\beta')\subset\beta^{{\O}'}\,,$ where $H' = \hat H^{{\O}'}\cap G'$ is the
converse image of $H$ in $G\,;$ now, it is clear that
$$N_G(H'_{\beta'}) = N_G(\hat H_\beta)\qq C_G(H') = C_G(\hat H)
\eqno £3.3.3,$$ 
and the bottom isomorphism in~£3.3.1 follows from isomorphism~£3.1.3.

\medskip
£3.4. More generally, for any $\hat G\-$interior algebra $A\,,$ it is clear that the extension 
$A' = \O'\otimes_\O A$ becomes a $G'\-$interior algebra and, as in the corollary above, the inclusion
$A\i A'$ induces a bijection, preserving inclusions and localness, between the sets of pointed 
$\O^*\-$groups $\hat H_\beta$ on $A$ and of pointed groups $H'_{\beta'}$ on $A'$
fulfilling $\Lambda'\.Z'\i H'\,,$
which maps $\hat H_\beta$ on $H'_{\beta'}$ if and only if we have 
$$H' = \hat H^{{\O}'}\cap G' \qq \beta\i \beta'
\eqno £3.4.1.$$
Indeed, from the left-hand equality and from any $\O\-$basis of $\O'$ we easily get 
$A'^{H'} \cong \O'\otimes_\O A^H$ and therefore, since $k$ is also the {\it residue\/} 
field of $\O'\,,$ we obtain an obvious bijection  between the sets of points of $H'$ on $A'$ 
and of~$H$ on $A\,.$ Actually, the bijection above also preserves the {\it fusions\/} as it shows the next result.

\bigskip
\noindent
{\bf Proposition~£3.5.} {\it With the notation above, for any pair of pointed $\O^*\-$groups
$\hat H_\beta$ and $\hat K_\gamma$ on $A\,,$ we have
$$F_A (\hat K_\gamma,\hat H_\beta) = F_{A'} (K'_{\gamma'},H'_{\beta'})
\eqno £3.5.1\phantom{.}$$ 
where we are setting  $H' = \hat H^{{\O}'}\cap G' $ and $K' = \hat K^{{\O}'}\cap G'$
and the points $\beta'$ of~$H'$ and $\gamma'$ of $K'$ respectively contain $\beta$ and $\gamma\,.$\/}

\medskip
\noindent
{\bf Proof:} Since $A'_{\beta'}\cong\O'\otimes_\O A_\beta$ and $A'_{\gamma'}\cong\O'\otimes_\O A_\gamma\,,$
if $\hat\varphi\,\colon \hat K\to \hat H$ is an {\it $A\-$fusion\/} from $\hat K_\gamma$ to~$\hat H_\beta\,,$
the corresponding  $\hat K\-$interior algebra~{\it embedding\/} 
$$f_{\hat\varphi} : A_\gamma\too {\rm Res}_{\hat\varphi} (A_\beta) 
\eqno £3.5.2\phantom{.}$$
determines a $K'\-$interior algebra~{\it embedding\/} 
$$f'_{\varphi'} : A'_{\gamma'}\too {\rm Res}_{\varphi'} (A'_{\beta'}) 
\eqno £3.5.3\phantom{.}$$
\eject
\noindent
where $\varphi'\,\colon K'\to H'$ is the group homomorphism determined by $\hat\varphi\,,$
and it is clear that the inclusion $A'_{\gamma'}\i A'$ and the composition of $f'_{\varphi'}$
with the inclusion $A'_\beta\i A'$ are also $A'^*\-$conjugate.

\smallskip
Conversely, since any  {\it $A'\-$fusion\/} from $K'_{\gamma'}$ to~$H'_{\beta'}$ can be decomposed 
as a composition of an isomorphism with an inclusion [8,~2.11], we may assume that $\varphi'$
is an isomorphism; in this case, choosing $j\in \beta\i \beta'$ and $\ell\in \gamma\i \gamma'\,,$
it follows from [8,~2.11] that the right-hand multiplication by a suitable invertible element
$a'$ of $A'^*$ determines an isomorphism from $A'\ell$ considered as an 
$A'\otimes_{\O'} \O'K'\-$module by left- and right-hand multiplication, onto ${\rm Res}_{{\rm id}_{A'}\otimes\varphi'}(A'j)$ where $A'j$ is similarly considered as an  $A'\otimes_{\O'} \O'H'\-$module by left- and right-hand multiplication.

\smallskip
But, considering an $\O\-$basis $\Delta'$ of $\O'\,,$ it is quite clear that $A'\ell$ and $A'j$ are respectively
isomorphic to $(A\ell)^{\Delta'}$ and $(Aj)^{\Delta'}$ considered as $A\otimes_{\O} \O_*\hat K\-$ 
and $A\otimes_{\O} \O_*\hat H\-$modules. Consequently, $A\ell$ and ${\rm Res}_{{\rm id}_{A}\otimes\hat\varphi}(Aj)$
are isomorphic considered as $A\otimes_{\O} \O_*\hat K\-$modules; then, the quotients
$$A/A\ell\cong A(1 -\ell)\qq {\rm Res}_{{\rm id}_{A}\otimes\hat\varphi}(A/Aj)\cong {\rm Res}_{{\rm id}_{A}\otimes\hat\varphi}\big(A(1 -j)\big)
\eqno £3.5.4\phantom{.}$$
are also isomorphic considered as $A\otimes_{\O} \O_*\hat K\-$modules and therefore  the right-hand 
multiplication by a suitable invertible element $a$ of $A^*$ still determines an isomorphism
$$A\ell\cong {\rm Res}_{{\rm id}_{A}\otimes\hat\varphi}(Aj)
\eqno £3.5.5,$$
so that it follows again from [8,~2.11] that $\hat\varphi$ an {\it $A\-$fusion\/} from 
$\hat K_\gamma$ to~$\hat H_\beta\,.$ We are done.

\medskip
£3.6. As in the ordinary case, we call {\it block of\/} $\hat G$ any primitive idempotent $b$ 
of $Z({\O}_*\hat G)$~; thus, $\alpha = \{b\}$ is a point of $G$ (or of $\hat G\,$) on ${\O}_*\hat G$
and, if $P_\gamma$ is a defect pointed group  of~$G_\alpha$~, we call {\it defect of\/} $b$ 
the integer ${\rm d} = {\rm d}(b)$ such that $p^{\rm d} = \vert P\vert$ . 
Then, according to the corollary above, $G_\alpha$ and $P_\gamma$ respectively  determine 
pointed groups $G'_{\alpha'}$ and $P'_{\gamma'}$ on ${\O}'G'e_{\Lambda'}$ where we denote by $P'$ 
the Sylow $p\-$subgroup of $\hat P^{{\O}'}\cap G'\,,$ so that $P'_{\gamma'}$ is a defect
pointed group of~$G'_{\alpha'}\,;$ once again, we have $\alpha' = \{b'\}$
for a block $b'$ of $G'$ such that $b'e_{\Lambda'} = b'\,,$ but note that this
block has a defect ${\rm d}(b')$ such that $p^{{\rm d}(b')} = p^{\rm d}\vert Z'\vert\,;$
 in particular, since  isomorphism~£3.1.3 implies that
$${\rm rank}_{\O'}\bigl(({\O}'G'e')_{\gamma'}\bigr) = 
{\rm rank}_{\O}\bigl(({\O}_*\hat G)_{\gamma}\bigr)\vert Z'\vert
\eqno £3.6.1\phantom{.}$$
and we have $E_{G'}(P'_{\gamma'})\cong E_G(\hat P_\gamma)$ (cf.~isomorphisms~£3.3.1), 
as in the ordinary case $E_G(\hat P_\gamma)$ is a $p'\-$group, $p^{\rm d}$ divides $\,\hbox{rank}_{\O}\bigl(({\cal
O}_*\hat G)_{\gamma}\bigr)$ and we get
$$p^{-d}\,\hbox{rank}_{\O}\bigl(({\O}_*\hat
G)_{\gamma}\bigr)\equiv \vert E_G(\hat P_\gamma)\vert\quad
(\hbox{mod}\;p)
\eqno £3.6.2.$$
Moreover, note that the image $\bar b'$ of $b'$ in $k\bar G'$ is a block of $\bar G'$ and corresponds to the image 
$b'$ of $b$ in $k_*\hat G^k$ by isomorphism~£3.1.4; thus, $\bar b$ is a $k\-$block of~$\hat G^k$ and 
this correspondence clearly determines a bijection between the sets of blocks of $\hat G$ and $\hat G^k$
which preserves the {\it defects\/}.
\eject

\medskip
£3.7. It is clear that the $\K\-$algebra $\K_*\hat G\,b$ is semisimple and let us  assume
that $\K$ is a {\it splitting field\/} for $\hat G$ or, equivalently, that $\K_*\hat G\,b$ is 
a {\it direct product of full matrix algebras over $\K\,.$\/}
Since the $\O\-$algebra $\O_*\hat G\,b$ is symmetric, the $\O\-$module  of 
the {\it symmetric $\O\-$linear forms $Z^\circ ({\O}_*\hat G\,b)$\/} is a free 
$Z({\O}_*\hat G\,b)\-$module of rank one [12,~Proposition~2.2] and therefore, denoting by
 ${\rm Irr}_{\K}(\hat G,b)$ the set of characters of the simple ${\K}_*\hat G\,b\-$modules --- simply called
 {\it irreducible characters of $\hat G$\/} in the sequel ---
 the $Z({\O}_*\hat G\,b)\-$submodule
$$Z_{\rm ch}^\circ (\O_*\hat G\,b) = \bigoplus_{\chi  \in {\rm Irr}_\K(\hat G,b)}\,\O\.\chi
\eqno £3.7.1,$$    
 determines an ideal $Z_{\rm ch} ({\O}_*\hat G\,b)$ of $Z({\O}_*\hat G\,b)$ [12,~Proposition~2.9]. 
 Note that, if $\hat x\in \hat G$ and $C_G (\hat x)\not= C_G(x)$ where $x$ 
denotes the image of $\hat x$ in $G\,,$  the elements $\hat x$ 
and $\lambda\.\hat x$  are $G\-$conjugate  for some $\lambda\in \O^*-\{1\}$ and therefore any $\mu\in Z^\circ (\O_*\hat G)$ 
vanishes on $\hat x\,.$

\medskip
£3.8. Moreover, for any $\chi\in\hbox{Irr}_{\K}(\hat G,b)\,,$ recall that the image of the restriction of $\chi$ to 
$Z(\O_*\hat G\,b)$ is contained in $\chi(1)\O$ and that $\chi(1)^{-1}\chi$ defines an $\O\-$algebra homomorphism
from $Z(\O_*\hat G\,b)$ to $\O\,;$ we call {\it defect of\/} $\chi$ the integer
${\rm d}(\chi)$ such that
$${\chi\over\chi(1)}\bigl(Z_{\rm ch}({\O}_*\hat
G\,b)\bigr) = p^{{\rm d}(\chi)}{\O}
\eqno £3.8.1.$$
More explicitly, according to Corollary~£2.4, there is a finite field extension $\K'$ of~$\K$ such that, denoting by $\O'$
the ring of integers of~$\K'\,,$  the $\O'^*\-$group $\hat G^{{\O}'}$ contains a finite subgroup $G'$ covering~$G\,;$
then, $\chi$ can be extended to an irreducible character 
of $\hat G^{{\O}'}$ which, by restriction, determines an irreducible character $\chi'$
of $G'$ fulfilling $\chi'\bigl((\lambda'.1)\,b'\bigr) =\lambda'\,
\chi(1)$ for any $\lambda'\in\Lambda'\,,$  where $b'$ is the
corresponding block of $G'\,,$ and therefore we get
$${\chi'\over\chi'(1)}\bigl(Z_{\rm ch}({\O}'G'b')\bigr) = {\vert G'\vert\over\chi'(1)}\,{\O}'
\eqno £3.8.2.$$
But, it follows from isomorphisms~£3.1.3 and~£3.1.4 that
$$s_{Z'}\big(Z_{\rm ch}({\O}'G'b')\big) = \vert Z'\vert \,{\O}'\otimes_{\O}Z_{\rm ch}({\O}_*\hat G\,b)
\eqno £3.8.3.$$
Consequently, we actually get
$$p^{{\rm d}(\chi')} = p^{{\rm d}(\chi)}\,\vert Z'\vert\qq
{ \rm d}(\chi) = \vartheta_p\big({\vert G\vert\over \chi(1)}\big)
\eqno £3.8.4\phantom{.}$$
where $\vartheta_p\,\colon \Bbb Z\to \Bbb Z$ denotes  the {\it $p\-$adic valuation\/}; in particular, 
since ${\rm d}(b')$ coincides with the maximal value of $\,{\rm d}(\chi')$ when $\,\chi'$ 
runs over the set of irreducible characters of ${\K}G'b'\,,$ it is not difficult to prove that ${\rm d}(b)$ 
coincides with the maximal value of ${\rm d}(\chi)$ when $\chi$ runs over ${\rm Irr}_{\K}(\hat G,b)\,.$
\eject

\bigskip
\noindent
{\bf Proposition £3.9.} {\it Let $b$ be a block of $\,\hat G$ and $\hat P_\gamma$ a defect pointed $\cal
O^*\-$group of~$\hat G_{\{b\}}\,.$ Then, the following statements
are equivalent:\/}
\smallskip
\noindent
{\bf £3.9.1}\quad {\it The block $b$ has defect zero.\/}
\smallskip
\noindent
{\bf £3.9.2}\quad {\it The source algebra $({\O}_*\hat
G)_\gamma$ has $\O\-$rank one.\/}
\smallskip
\noindent
{\bf £3.9.3}\quad {\it The block algebra ${\O}_*\hat G\,b$ is
a full matrix algebra over $\O\,.$\/}
\smallskip
\noindent
{\bf £3.9.4}\quad {\it The block $\K\-$algebra ${\K}_*\hat
G\,b$ is a direct product of full matrix algebras over $\K$ and
there is $\chi\in\hbox{Irr}_{\K}(\hat G,b)$ with defect zero.
\smallskip
\noindent
In this case, we have $\Bbb O_p (G) = \{1\}\,.$\/}

\medskip
\noindent
{\bf Proof:} With the notation above, if $\hbox{d}(b) = 0$ then $Z'$ is the unique defect group of 
the corresponding block  $b'$ of $G'$ and therefore the $Z'\-$interior source algebra of $b'$ is isomorphic to
${\O}'Z'$ [9, Proposition~14.6], so that statement~£3.9.2 follows from isomorphism~£3.1.3. 
Since the block algebra is always {\it Morita equivalent\/}  to the source algebra [7,~Corollary~3.5],  
statement~£3.9.3 follows from statement~£3.9.2.

 \smallskip
 If statement~£3.9.3 holds then $\hbox{Irr}_{\K}(\hat G,b)$ has a unique element $\chi$ 
 and, since an indecomposable direct summand of ${\O}_*\hat G\,b$ as ${\O}_*\hat G\,b\-$module 
has $\O\-$rang equal to  $\chi(1)\,,$ $\vert G\vert_p$ divides $\chi(1)$ and therefore we have 
${\rm d}(\chi) = 0\,.$ Finally, statement~£3.9.4 implies that the sum
$$e_\chi = {\chi(1)\over \vert G\vert}\sum_{x\in G}\chi(\hat x)\.\hat x^{-1}
\eqno £3.9.5\phantom{.}$$
where $\hat x\in \hat G$ is a lifting of $x\in G\,,$ belongs to $Z({\O}_*\hat G)\,;$ but, $e_\chi$ is an idempotent
  and therefore we get $e_\chi = b\,,$ so that $\chi$ is the unique element in~$\hbox{Irr}_{\K}(\hat G,b)\,,$ which forces 
$\hbox{d}(b) = \hbox{d}(\chi) = 0\,;$ moreover, in this case the image $\bar b$ of $b$ in $k_*\hat G^k$
is a block of {\it defect zero\/} of $\hat G^k\,,$ so that there is a projective simple $k_*\hat G^k\-$module 
which forces  $\Bbb O_p (G) = \{1\}\,.$

\bigskip
\noindent
{\bf Proposition £3.10.} {\it For any pair of local pointed $\O^*\-$groups $\hat P_\gamma$ 
and $\hat Q_\delta$ on ${\O}_*\hat G\,,$ we have
$$E_G(\hat Q_\delta,\hat P_\gamma) = F_{{\O}_*\hat G}(\hat Q_\delta,\hat P_\gamma)
\eqno £3.10.1.$$\/}
\par
\noindent
{\bf Proof:} We already know that the left-hand member is contained in the right-hand one (cf.~£2.11.3). 
For any $\tilde\varphi\in  F_{{\O}_*\hat G}(\hat Q_\delta,\hat P_\gamma)$ 
and any representative $\varphi$ of~$\tilde\varphi\,,$  it is quite clear that there is a local point 
$\varepsilon$ of $\hat R = \varphi(\hat Q)$ on ${\O}_*\hat G$ such that $\hat R_\varepsilon\subset \hat
P_\gamma$ and that $\varphi$ determines an element of 
$F_{{\O}_*\hat G}(\hat Q_\delta,\hat R_\varepsilon)$ ; thus, we may assume that the 
$\O^*\-$quotients $P$ of $\hat P$ and $Q$ of $\hat Q$ have the same
order; in this case, according to £2.11 and choosing $i\in \gamma$
and $j\in \delta$~, there is $a\in  ({\O}_*\hat G)^*$ fulfilling
$(u.j)^a = \varphi(u).i\,$ for any $u\in\hat Q\,.$

\smallskip
With the notation above, let $P'_{\gamma'}$ and $Q'_{\delta'}$ be the corresponding local pointed groups 
on ${\O}'G'e'$ (cf. Corollary~£3.3), and choose $i'\in\gamma'$ and $j'\in\delta'$ respectively
 lifting $i$ and $j\,;$ then, ${\O}'Q'j'$ is an indecomposable direct summand of
$j'({\O}'G')j'$ as ${\O}'(Q'\times Q')\-$modules and therefore the quotient
$$({\O}'_*\hat Q^{{\O}'})j\cong {\O}'\otimes_{{\O}'Z'} {\O}'Q'j'
\eqno £3.10.2\phantom{.}$$ 
is a direct summand of the quotient
$$j({\O}'_*\hat G^{{\O}'})j\cong  {\O}'\otimes_{{\O}'Z'} j'({\O}'G')j'
\eqno £3.10.3\phantom{.}$$
always as ${\O}'(Q'\times Q')\-$modules. The ${\O}'(Q'\times Q')\-$module $({\O}'_*\hat Q^{{\O}'})j$
is still indecomposable since,
denoting by $\Delta(Q')$ the diagonal subgroup of $Q'\times Q'$ and
considering $\O'$ as an ${\O}'\bigl((Z'\times Z').\Delta(Q')
\bigr)\-$module where $\Delta(Q')$ acts trivially and $Z'\times Z'$
throughout  the {\it canonical\/} homomorphism $Z'\to \O'^*\,,$ we have
$$({\O}'_*\hat Q^{{\O}'})j\cong \hbox{Ind}_{(Z'\times Z').\Delta(Q')}^{Q'\times Q'}({\O}')
\eqno £3.10.4.$$

\smallskip
Consequently, since we clearly have
$$\,j({\O}_*\hat G)ja = j({\O}_*\hat G)i
\eqno £3.10.5,$$
$\,({\O}'_*\hat Q^{{\O}'})ja\,$ becomes an
indecomposable direct summand of $j({\O}'_*\hat G^{{\O}'})i$
as ${\O}'(Q'\times P')\-$modules and therefore, since we have the
 ${\O}'(Q'\times P')\-$module isomorphism 
$$j({\O}'_*\hat G^{{\O}'})i\cong {\O}'\otimes_
{{\O}'Z'} j'({\O}'G')i'
\eqno £3.10.6,$$
the ${\O}'(Q'\times P')\-$module $j'({\O}'G')i'$ admits a direct summand of $\O'\-$rank equal to 
$\vert P\vert\vert Z'\vert = \vert P'\vert$ ; thus, the ${\O}'(Q'\times P')\-$module ${\O}'G'e'$
admits such a direct summand, which forces the existence of $x'\in G'$
fulfilling ${\O}'Q'x' = x'{\O}'P'\,,$ together with an ${\O}'
(Q'\times P')\-$module isomorphism
$$({\O}'_*\hat Q^{{\O}'})ja\cong {\O}'\otimes_{{\O}'Z'} {\O}'Q'x'
\eqno £3.10.7.$$

\smallskip
But, the left-hand member of this isomorphism admits a vertex $V'$ equal to
$\{\big(z'u',\varphi(u')\big)\}_{z'\in Z',u'\in Q'}$ and a $V'\-$source
of $\O'\-$rank one where the subgroup   $\{\big(u', \varphi(u')\big)\}_{u'\in Q'}\i V'$ is the kernel, 
whereas the right-hand member admits a vertex equal
to $\{(z'u',u'^{x'})\}_{z'\in Z',u'\in Q'}$ and the corresponding
kernel is equal to $\{(u',u'^{x'})\}_{u'\in Q'}\,;$ hence, up to modifying our
choice of $x'\,,$ we may assume that $\varphi(u') = u'^{x'}$ for any $u'\in Q'\,.$
In conclusion, $x'$ induces $\varphi\,\colon \hat Q\to \hat P$
and therefore $\tilde\varphi$ belongs to $E_G(\hat Q_\delta, \hat P_{\delta^{x'}})\,;$ 
in particular, since $E_G(\hat Q_\delta, \hat
P_{\delta^{x'}})\subset F_{{\O}_*\hat G}(\hat Q_\delta, \hat
P_{\delta^{x'}})$ and $\varphi$ is an isomorphism,
$\widetilde{\hbox{id}}_{\hat P}$ belongs to $F_{{\O}_*\hat
G}(\hat P_{\delta^{x'}},\hat P_\gamma)\,,$ which forces $\,\delta^{x'}
= \gamma\,$ (cf.~£2.11) and $\tilde\varphi\in E_G(\hat Q_\delta,\hat P_\gamma)\,.$
We are done.
\eject

\bigskip
£3.11. Let $\hat Q_\delta$ be a local pointed $\O^*\-$group on ${\O}_*\hat G$~; if~$\hat H$ is an
$\O^*\-$sub-group of $\hat G$ containing $\hat Q.C_{\hat G}(\hat Q)$ then from isomorphism~£3.2.2 we get (cf.~£2.8)
$$s_\delta\big(({\O}_*\hat G)^{H}\big)\i s_\delta\big(({\O}_*\hat G)^{Q.C_{G}(\hat Q)}\big)
 = ({\O}_*\hat G)(Q)^{C_G (\hat Q)}\cong k
\eqno £3.11.1\phantom{.}$$
and therefore there is a unique point $\beta$ of $\hat H$ on ${\cal
O}_*\hat G$  such that $\hat Q_\delta\subset \hat H_\beta$ ; in
particular, there are unique points $\zeta$ of $\hat Q.C_{\hat
G}(\hat Q)$ and $\nu$ of $N_{\hat G}(\hat Q_\delta)$ on ${\cal
O}_*\hat G$ fulfilling
$$\hat Q_\delta\i \hat Q.C_{\hat G}(\hat Q)_\zeta\i N_{\hat G}(\hat Q_\delta)_\nu
\eqno £3.11.2\phantom{.}$$
and, according to Corollary~£3.3, $\hat Q_\delta$ is a {\it defect pointed $\O^*\-$group\/}
 of $N_{\hat G}(\hat Q_\delta)_\nu$ if and
only if it is a {\it maximal local pointed $\O^*\-$group\/} on ${\O}_*\hat G$ [7,~Corollary~1.4]. 
We denote by $b(\delta)$ the block of $\hat Q.C_{\hat G}(\hat Q)$ fulfilling $b(\delta)\.\zeta = \zeta\,,$ by $\bar C_{\hat G^k}(\hat Q)$ the
$k^*\-$group $C_{\hat G^k}(\hat Q)/C_Q(\hat Q)$ , and by $\bar
b(\delta)$ the image of $b(\delta)$ in $k_*\bar C_{\hat G^k}(\hat Q)$
, which is a block of $\bar C_{\hat G^k}(\hat Q)$ over $k$ .   We say that $\hat Q_\delta$ is {\it selfcentralizing\/},
or that it is a {\it selfcentralizing pointed $\O^*\-$group\/} on $\O_*\hat G\,,$  if $\hat Q$
is a {\it defect pointed $\O^*\-$group\/} of~$\hat Q.C_{\hat G}(\hat Q)_\zeta\,;$ moreover, as in [16,~3.1], 
we say that $\hat Q_\delta$ is {\it radical\/} if it is {\it selfcentralizing\/} and we have
$$\Bbb O_p\big(E_G (\hat Q_\delta)\big) = \{1\}
\eqno £3.11.3.$$

\bigskip
\noindent
{\bf Proposition £3.12.} {\it For any local pointed $\O^*\-$group $\hat Q_\delta$ on ${\O}_*\hat G$ , the following
statements are equivalent:\/}
\smallskip
\noindent
{\bf £3.12.1}\quad {\it $\hat Q_\delta$ is selfcentralizing.\/}
\smallskip
\noindent
{\bf £3.12.2}\quad {\it The structural homomorphism induces $kC_Q(\hat Q)\cong ({\O}_*\hat G)_\delta (Q) \,.$\/} 
\smallskip
\noindent
{\bf £3.12.3}\quad {\it The block $\bar b(\delta)$ of $\bar C_{\hat G^k}(\hat Q)$ has defect zero.\/}
\smallskip
\noindent
{\bf £3.12.4}\quad {\it For any local pointed $\O^*\-$group $\hat P_\gamma$ on ${\O}_*\hat G$ 
containing $\hat Q_\delta\,,$ $\hat Q$ contains $C_{\hat P}(\hat Q)\,.$\/}
\medskip
\noindent
{\bf Proof:} With the notation above, it follows from Proposition~£3.2 that $s_{Z'}$ determines a surjective
$k\-$algebra homomorphism
$$k\otimes_{\O'}(\O'G'\,e_{\Lambda'})^Q\too (\O'_*\hat G^{\O'})(Q)
\eqno £3.12.5;$$
but, setting $Q' = \hat Q^{\O'}\cap G'$ and denoting by $Q'_{\delta'}$ and by $Q'.C_{G'}(Q')_{\zeta'}$ 
the respective pointed groups on~$\O'G'\,e_{\Lambda'}$ corresponding to $\hat Q_\delta$ 
and to $\hat Q.C_{\hat G}(\hat Q)_\zeta$ (cf.~Corollary~£3.3), statement~£3.12.1 implies that
$Q'_{\delta'}$ is a defect pointed group of~$Q'.C_{G'}(Q')_{\zeta'}\,;$ in this case, it follows from
[9,~statements~14.5.1 and~2.9.2] that the structural homomorphism induces an isomorphism
$$kZ(Q')\cong (\O'G'\,e_{\Lambda'})_{\delta'}(Q')
\eqno £3.12.6;$$
consequently, the structural homomorphism (cf.~isomorphism~£3.2.2)
$$kC_Q (\hat Q)\too (\O'_*\hat G^{\O'})(Q)\cong ({\O}_*\hat G)_\delta (Q)
\eqno £3.12.7\phantom{.}$$
is surjective, which proves statement~£3.12.2.
\eject

\smallskip
Now, if we have the $k\-$algebra isomorphism 
$$ kC_Q(\hat Q)\cong ({\O}_*\hat G)_\delta (Q)
\eqno £3.12.8,$$ 
then we have a $k\-$algebra embedding $k\to k_*\bar C_{\hat G^k}(\hat Q)\bar b(\delta)$ 
 (cf.~isomorphism £3.2.2) which determines an indecomposable projective $k_*\bar C_{\hat G^k}
 (\hat Q)\bar b(\delta)\-$module isomorphic to its own socle; thus, the block algebra 
 $k_*\bar C_{\hat G^k} (\hat Q)\bar b(\delta)$ admits a simple projective module
and therefore the block $\bar b(\delta)$ has defect zero.

\smallskip
Assume that statement~£3.12.3 holds and let $\hat P_\gamma$ be a local pointed $\O^*\-$ group on 
${\O}_*\hat G$ containing $\hat Q_\delta\,;$ then, there are successively local points $\varepsilon$ 
of $\hat Q.C_{\hat P}(\hat Q)$ and $\varphi$ of $\hat Q$ on ${\O}_*\hat G$ fulfilling 
$$\hat Q_\varphi\subset \hat Q.C_{\hat P}(\hat Q)_\varepsilon\subset \hat P_\gamma
\eqno £3.12.9\phantom{.}$$ 
which forces $b(\varphi) = b(\delta)$ [2,~Theorem~1.8]
and therefore, since we have  (cf.~isomorphism~£3.2.2)
$$k_*\bar C_{\hat G^k}(\hat Q)\bar b(\delta)\cong ({\O}_*\hat G)(\hat Q_\delta)
\eqno £3.12.10,$$
we get $\varphi = \delta$ ; in particular, since we have [2,~Proposition~1.5]
$$\big((\O_*\hat G)(Q)\big)\big(C_P (\hat Q)\big)\cong (\O_*\hat G)\big(Q.C_P(\hat Q)\big)
\eqno £3.12.11,$$
$\hbox{Br}_Q(\varepsilon)$ is a local point of $C_P(\hat Q)$ on $(\O_*\hat G)(Q)\,;$ but,
it follows from [13,~Theorem~2.9] that the canonical map
$$k_* C_{\hat G^k}(\hat Q)\too k_*\bar C_{\hat G^k}(\hat Q)
\eqno £3.12.12\phantom{.}$$
is a {\it semicovering $C_{\hat G^k}(\hat Q)\-$interior algebra homomorphism\/}; hence,
according to [6,~Proposition~3.15] and to isomorphism~£3.2.2 above, ${\rm Br}_Q(\varepsilon)$
determines a local point of $C_P(\hat Q)/C_Q(\hat Q)$ on $k_*\bar C_{\hat G^k}(\hat Q)\bar b(\delta)\,,$ 
which forces this $p\-$group to be trivial, so that $C_{\hat P}(\hat Q) = Z(\hat Q)$ .

\smallskip
Finally, it follows from Corollary~£3.3 that statement~£3.12.4 implies that, for any local pointed groups $P'_{\gamma'}$
containing $Q'_{\delta'}\,,$ we have $C_{P'}(Q') =Z(Q')\,;$ then, according to [14,~4.8 and~Corollary~7.3],  $Q'_{\delta'}$ is
selfcentralizing and therefore it is a defect pointed group of $Q'.C_{G'}(Q')_{\zeta'}\,;$ hence, once again by Corollary~£3.3,
$\hat Q_\delta$ is selfcentralizing too.

\bigskip
\noindent
{\bf Corollary £3.13.} {\it If $\hat P_\gamma$ and $\hat Q_\delta$ are
local pointed $\O^*\-$groups on ${\O}_*\hat G$ such that $\hat Q_\delta\i \hat P_\gamma$ and $\hat Q_\delta$ is
selfcentralizing then $\hat P_\gamma$ is selfcentralizing too.\/}
\medskip
\noindent
{\bf Proof:} In this situation, it is clear that  $\hat P_\gamma$ fulfills also statement~£3.12.4.

\bigskip
\noindent
{\bf Proposition £3.14.} {\it For any local pointed $\O^*\-$group $\hat P_\gamma$ on ${\O}_*\hat G$
 there is a unique local point $\gamma^k$ of $P$  on $k_*\hat G^k$ containing the image of $\gamma\,,$
 the canonical group homomorphism
 $$E_G(\hat P_\gamma)\too E_G(P_{\gamma^k})
 \eqno £3.14.1\phantom{.}$$
 is surjective and its kernel is an abelian $p\-$group. In particular, $P_{\gamma^k}$ is selfcentralizing 
 if and only if $E_G(\hat P_\gamma)$ is isomorphic to $E_G(P_{\gamma^k})$  and $\hat P_\gamma$ is selfcentralizing too. Moreover, for any local pointed group $Q_{\delta^k}$ on $k_*\hat G^k\,,$
 there is a  local point $\delta$ on $\O_*\hat G$ of the converse image $\hat Q$ of $Q$ in $\hat G$ such that 
 $\delta^k$ contains the image of $\delta$ and that the canonical map
$$E_G(\hat Q_\delta,\hat P_\gamma)\too E_G(Q_{\delta^k},P_{\gamma^k})
\eqno £3.14.2$$
is surjective.\/}
\medskip
\noindent
{\bf Proof:} With the notation above, it follows from Corollary~£3.3 that $\hat P_\gamma$ determines
a local pointed group $P'_{\gamma'}$ on ${\O}'G'e_{\Lambda'}\,;$ but, it is clear that the image $\gamma'^k$
of $\gamma'$ in $kG'e_{\Lambda'}$ is a local point of $P$ on this $k\-$algebra; moreover, it follows from
[13,~Theorem~2.9] that the canonical map
$$kG'e_{\Lambda'}\too k_*\hat G^k
\eqno £3.14.3\phantom{.}$$
is a {\it semicovering $G\-$algebra homomorphism\/}; hence, according to [6,~Proposition~3.15], the image of
$\gamma'^k$ is contained in a local point $\gamma^k$ of $P$  on  $k_*\hat G^k\,,$ 
which then contains the image of $\gamma\,;$ moreover, if follows from  [13,~Corollary~2.13] that
 $P_{\gamma^k}$ is selfcentralizing 
 if and only if $E_G(\hat P_\gamma)$ is isomorphic to $E_G(P_{\gamma^k})$  and $\hat P_\gamma$ is selfcentralizing too. 

\smallskip
The existence of the local point $\delta$ of $\hat Q$ in the converse image of $\delta^k$ follows 
from  [6,~Proposition~3.15] and from the fact  that homomorphism~£3.14.3 is a {\it semicovering\/};
moreover, if  $E_G(Q_{\delta^k},P_{\gamma^k})$ is {\it not empty\/} then it is quite clear that we
can choose $\delta$ in such a way that $E_G(\hat Q_\delta,\hat P_\gamma)$ is also {\it not empty\/};
in this case,  the surjectivity of the canonical map~£3.14.2 follows from [13,~Theorem~2.9]{\footnote{\dag}{\cds 
In [13,~Theorem~2.9] the necessary assumption that there is at least an {\cdt A}-fusion from {\cdt Q}$_{_\delta}$
to {\cdt P}$_{_\gamma}$ has been forgotten.}}
 and from Proposition~£3.10 above.

\bigskip
\noindent
{\bf £4. Fong reduction for interior algebras over an $\O^*\-$group}
\bigskip
£4.1.  As a matter of fact,  all our arguments in~[15,~\S3] on interior algebras over finite $k^*\-$groups
can be translated to interior algebras over finite $\O^*\-$groups; we will explicit the translation of the statements since
they demand some modifications, and some indications on the proofs if necessary; but, it seems useless
to repeat the identical part of the proofs.

\medskip
£4.2.  The first modification concerns, for any finite $p\-$group $P\,,$ the involved {\it Dade $P\-$algebras\/}; in our
present situation, we call  {\it Dade $P\-$algebra\/} a full matrix algebra $S$ over $\O$ endowed with an action of
$P$ which stabilizes an $\O\-$basis of $S$ containing $1_S\,;$ recall that two Dade $P\-$algebras $S$ and $S'$ are {\it similar\/} if $S$ can be {\it embedded\/}  (cf.~£2.8) in the {\it tensor product\/} ${\rm End}(N)\otimes_k S'$ for a suitable $kP\-$module 
$N$ with a $P\-$stable $\O\-$basis~[11,~1.5 and~2.5.1]. Moreover, recall that $S(P)$ is a simple $k\-$algebra [11,~1.8];
in particular, if $S$ is {\it primitive\/} then we have 
$${\rm rank}_\O (S) \equiv 1 \bmod p
\eqno £4.2.1;$$
hence, in all the cases, it follows from [11,~3.13] that the action of $P$ on $S$ can be lifted to a group homomorphism
$P\to S^*$ and we will consider $S$ as a $P\-$interior algebra. Since  $S(P)$ is simple, $P$ has a unique local point 
$\pi$ on~$S$ [11,~1.8] that very often we omit, respectively writing  $F_S (P)$ and $\hat F_S (P)$ instead of 
$F_S (P_\pi)$ and $\hat F_S (P_\pi)$ 
(cf.~£2.11 and~£2.12). Moreover, since we have an evident $k^*\-$group homomorphism
$$\hat F_S (P)\too \hat F_{k\otimes_\O S} (P)
\eqno £4.2.2,$$
any {\it polarization\/} $\omega$ considered in [16,~2.15] still supplies  a $k^*\-$group homomorphism
$$\omega_{(\hat P,S)} : \hat F_S (P)\too k^*
\eqno £4.2.3.$$

\medskip
£4.3. Let $\hat P$ be an $\O^*\-$group with $\O^*\-$quotient $P\,;$ we call {\it monomial\/} any $\O_*\hat P\-$module $M$  such that 
the sources of any indecomposable  direct summand of~$M$ have $\O\-$rang {\it one\/} or, equivalently, 
it admits an $\O\-$module decomposition $M = \bigoplus_{i\in I} L_i$ such that ${\rm rank}_\O(L_i) = 1$ 
and that $\hat P$ stabilizes the family~$\{L_i\}_{i\in I}\,;$ then, $k\otimes_\O M$ is a {\it permutation\/}
$kP\-$module. More precisely, we call {\it twisted diagonal\/} $\O_*(\hat P\,\hat\times\, \hat P^\circ)\-$module 
any {\it monomial\/} $\O_*(\hat P\,\hat\times\, \hat P^\circ)\-$mo-dule  $M$ such that the intersections 
of any vertex of any indecomposable  direct summand of~$M$ with the images in 
$\hat P\,\hat\times\, \hat P^\circ$ of $\hat P$ and $\hat P^\circ$  coincide with the image of $\O^*\,.$

\medskip
£4.4. Let $\hat G$ be an $\O^*\-$group with finite $\O^*\-$quotient $G$ and $A$ a $\hat G\-$interior algebra;
as in [13,~2.8], we denote by $\L_A$ the {\it local category of $A$\/} where the objects are the {\it local
pointed groups\/} on $A$ and the morphisms are the {\it $A\-$fusions\/} (cf.~£2.11) between them with
the usual composition [8,~Definition~2.15]. Let~$S$ be a $G\-$stable unitary subalgebra of~$A$ isomorphic to a {\it direct product of full matrix algebras 
over  $\O$\/} and assume that $G$ acts transitively on the set~$I$ of primitive idempotents of~the center
$Z(S)$ of $S\,;$ let $i$ be an element of~$I$ and denote by~$\hat H$ the stabilizer~of~$i$ in~$\hat G\,.$ 
Thus, the $\O^*\-$quotient $H$ of $\hat H$ acts on the full $\O\-$matrix algebra $Si$ determining 
an $\O^*\-$group $^{\skew4\hat {}}H\,,$ together with an $\O^*\-$group homomorphism 
$\rho\,\colon {}^{\skew4\hat {}}H\to (Si)^*$ (cf.~£2.7) and we set (cf.~£2.6)
$$ H^{\skew4\hat {}}\, = \hat H * (\,{}^{\skew4\hat {}}H)^\circ
\eqno £4.4.1;$$
moreover, if $P$ is a $p\-$subgroup of $H$ and ${\rm Res}_P^H (Si)$ is a Dade $P\-$algebra,
according to~£4.2 the converse image  ${}^{\skew4\hat {}}P$ of $P$ in $\,^{\skew6\hat {}}H$ is {\it split} 
and therefore, up to the choice of a splitting ${}^{\skew4\hat {}}P\cong \O^*\times P\,,$ we can identify to 
each other  the converse images in  $\hat H$ and~$H^{\skew4\hat {}}$ of any subgroup $Q$ of $P\,;$ 
if $p$ does not divide ${\rm rank}_\O (Si)\,,$ we fix our choice assuming that 
$\rho (1\times P)\i {\rm Ker}({\rm det}_{Si})\,.$

\bigskip
\noindent
{\bf Proposition £4.5.} {\it With the notation and the hypothesis above, there
exists an $H^{\skew4\hat {}}\,$-interior algebra $B\,,$ unique up to isomorphisms,
such that we have a $\hat G$-interior algebra isomorphism
$$A\cong {\rm Ind}_{\hat H}^{\hat G}(Si\otimes_k B)
\eqno £4.5.1$$ 
mapping $s\in Si$ on $1\otimes (s\otimes 1_B)\otimes 1\,.$ In
particular,
$A$ and $B$ are Morita equivalent.\/}
\medskip
\noindent
{\bf Proof:} Since [7,~Proposition~2.1] holds over $\O\,,$ the proof of [15,~Proposition~3.2] applies.

\bigskip
\noindent
{\bf Corollary £4.6.} {\it With the notation and the hypothesis above, assume that
$B$ has a unique $H\-$conjugacy class of maximal local pointed $\O^*\-$groups
$\hat P_\gamma\,,$ that ${\rm Res}_P^H (Si)$ is a Dade $P\-$algebra and that~$A$ and $B$ are
twisted diagonal $\O_*(\hat P\,\hat\times\,\hat P^\circ)\-$modules by left- and right-hand multiplication. 
Choosing a splitting $\,{}^{\skew4\hat {}}P\cong \O^*\times P\,,$   for any local pointed $\O^*\-$group 
$\hat Q_\delta$ on $B$ we have a local point $\iota (\delta)$ of $\hat Q$ on $A$ such that isomorphism~{\rm £4.5.1} induces 
a $\hat Q\-$interior algebra embedding
$$A_{\iota (\delta)}\longrightarrow {\rm Res}_Q^H(Si)\otimes_\O B_\delta
\eqno £4.6.1\phantom{.}$$
and this correspondence determines an equivalence of categories $\iota\colon
\L_B\to \L_A$ between the local categories of $B$ and $A\,.$ In particular, $A$
 has a unique $G\-$conjugacy class of maximal local pointed $\O^*\-$groups.\/}

\medskip
\noindent
{\bf Proof:} The proof of [15,~Corollary~3.3] shows the existence of a map $\iota$ between
the sets of local pointed $\O^*\-$groups on $B$ and on $A\,.$ Then, in order to prove the equality
$$F_B (\hat R_\varepsilon,\hat Q_\delta) = F_A (\hat R_{\iota(\varepsilon)},\hat Q_{\iota (\delta)})
\eqno £4.6.2\phantom{.}$$
for any pair  of local  pointed $\O^*\-$groups $\hat Q_\delta$ and $\hat R_\varepsilon$ on $B\,,$ note that
it follows from Proposition~£3.5 above that [8,~Proposition~2.14] and [10,~Theorem~5.3]
apply to our present situation; thus, we still have
$$F_{iAi} (\hat R_{\iota(\varepsilon)},\hat Q_{\iota (\delta)}) =F_A (\hat R_{\iota(\varepsilon)},\hat Q_{\iota (\delta)})
\eqno £4.6.3\phantom{.}$$
and, since   we have $iAi\cong Si\otimes_\O B\,,$ it suffices to prove that $F_{Si} (R,Q)$
contains both $F_{iAi} (\hat R_{\iota(\varepsilon)},\hat Q_{\iota (\delta)})$ and 
$F_B (\hat R_\varepsilon,\hat Q_\delta)\,.$

\smallskip
Let $\hat \varphi\,\colon \hat R\to \hat Q$ be an $\O^*\-$group homomorphism which  belongs either to
$F_A (\hat R_{\iota(\varepsilon)},\hat Q_{\iota (\delta)})$ or to 
$F_B (\hat R_\varepsilon,\hat Q_\delta)\,;$
once again, according  to  Proposition~£3.5 above, [8,~Proposition~2.18] apply to our present situation and
therefore, since ${\rm Res}_P^H (Si)$ is a Dade $P\-$algebra, if follows from [11,~statement~2.5.1] that
 it suffices to prove that the {\it Dade $R\-$algebras\/} ${\rm Res}_R^H (Si)$ and 
 ${\rm Res}_{\varphi}\big({\rm Res}_Q^H (Si)\big)$ are {\it similar\/} (cf.~£4.2). But, it follows from 
 [11,~statement~1.5.2] that these {\it Dade $R\-$algebras\/} are {\it similar\/} if and only if the corresponding
 {\it Dade $R\-$algebras\/} over~$k$ are so; moreover, according to our hypothesis, the actions of $P\times P$ 
 on~$k\otimes_\O A$ and $k\otimes_\O B$ by left and right multiplication stabilize bases where $P\times \{1\}$ 
 and~$\{1\}\times P$ act freely, so that the hypothesis of [15,~Corollary~3.3] are fulfilled.

\smallskip
Consequently,   the proof of  [15,~Corollary~3.3] applies, proving that the corresponding 
{\it Dade $P\-$algebras over $k$\/}
$${\rm Res}_R^H (k\otimes_\O Si)\qq {\rm Res}_{\varphi}\big({\rm Res}_Q^H (k\otimes_\O Si)\big)
\eqno £4.6.4\phantom{.}$$
are indeed {\it similar\/} to each other  which shows that the functor $\iota\,\colon \L_B\to \L_A$ is {\it fully faithful\/};
once again,  the proof of  [15,~Corollary~3.3] applies showing that this functor is {\it essentially surjective\/}, 
so that it is an equivalence of categories. We are done.

\medskip
£4.7. The main point in {\it Fong reduction\/} is that if $A$ is a
{\it block algebra\/} $\O_*\hat Gb$ for a block $b$ of $\hat G$ then $i$ is a block
of $\hat H$ and if moreover $p$ does note divide  ${\rm rank}_\O (Si)$ then $B$ is also a 
{\it block algebra\/}. Denote by $V$ a simple $k\otimes_\O Si$-module, which becomes an
$k_*\,{}^{\skew4\hat {}} H^k$-module through~$\rho$ (cf.~£4.4).

\bigskip
\noindent
{\bf Proposition £4.8.} {\it  With the notation and the hypothesis above, 
if $A\cong \O_*\hat Gb$ for a block $b$ of $\hat G$ then $i$ is a block
of $\hat H$ which belongs to a point $\beta$ of $\hat H$ on~$A$ and we have
$i(\O_*\hat G)i = \O_*\hat Hi\,.$ In particular, we have an equivalence of categories
$\L_{\O_*\hat Hi}\cong \L_{\O_*\hat Gb}\,.$\/}
\medskip
\noindent
{\bf Proof:} Replacing [15,~Corollary~3.3] by Corollary~£4.6 above, the proof of [15,~Proposition~3.5] applies
and then the last statement follows from Corollary~£4.6 above.

\bigskip
\noindent
{\bf Theorem £4.9.} {\it With the notation and the hypothesis above, assume that we have
$A\cong \O_*\hat Gb$ for a block $b$ of $\hat G$ and that $p$ does not divide  ${\rm rank}_\O (Si)\,.$
Then, we have $B\cong \O_*H^{\skew4\hat {}}c$ for a block $c$ of $H^{\skew4\hat {}}$ and 
$V$ is a simple $k_*\,{}^{\skew4\hat {}} H^k\-$mo-dule.  Moreover, if  $\,{\rm Res}_P^H (Si)$ is a Dade $P\-$algebra, 
we have an equivalence of categories $\L_{\O_*H^{\skew4\hat {}}\,c}\cong \L_{\O_*\hat Gb}\,.$\/} 
\medskip
\noindent
{\bf Proof:} It follows from Proposition~£3.5 above that [10,~\S4]
apply to our present situation; thus, replacing [15,~Corollary~3.3] by Corollary~£4.6 above, 
the proof of [15,~Theorem~3.6] applies.

\bigskip
\noindent
{\bf Theorem £4.10.} {\it With the notation and the hypothesis above, assume that
$A\cong \O_*\hat Gb$ for a block $b$ of $\hat G$ and that  $S = \O_*\hat Kb$
  for a normal $\O^*\-$subgroup~$\hat K$ of~$\hat G$~having a block $d$ of defect
zero such that $db \not= 0\,.$ Then, $K$ is a normal subgroup of $H^{\skew4\hat
{}}$ and we have $B\cong \O_*(H^{\skew4\hat {}}/K)\bar c$ for a block $\bar c$ of
$H^{\skew4\hat {}}/K\,.$\/}

\medskip
\noindent
{\bf Proof:} The corresponding part of the proof of [15,~Theorem~3.7] applies.

\medskip
£4.11. It is well-known that, up to replacing $\O$ by $\K\,,$ the so-called {\it Clifford's reduction\/}
 can be viewed as a particular case of Fong's reduction. Explicitly, let $\hat N$ be a normal 
 $\O^*$-subgroup of $\hat G$  and $\nu$ an {\it absolutely\/} irreducible character  of $\hat N\,;$
denote by ${\rm Irr}_\K(\hat G,\nu)$ the set of $\chi\in {\rm Irr}_\K(\hat G)$ such that $\nu$ is involved in
${\rm Res}_{\hat N}^{\hat G}(\chi)\,,$ by $e_\nu$ the primitive  idempotent of $Z({\K}_*\hat N)$ associated with $\nu$ and by
$H$ and $\hat H$ the respective stabilizers of $\nu$ in $G$ and $\hat G\,.$  If~$\hat\K$ is an extension
of $\K$ of degree divisible by $\vert G\vert\,,$ it follows from Proposition~£2.3 that, 
denoting by $\hat\O$ the ring of integers of $\hat\K$ and setting
$$\skew2\hat{\hat G} = \hat G^{\hat\O}\quad ,\quad \skew2\hat{\hat H} = \hat H^{\hat\O}\qq
\skew3\hat{\hat N} = \hat N^{\hat\O}
\eqno £4.11.1,$$
 there are an $\hat\O^*$-group $\,\skew1\hat{\hat{}} H$ with  $\hat\O^*$-quotient~$H$ containing and 
 normalizing~$\skew3\hat{\hat N}\,,$ and an $\hat\O^*$-group homomorphism 
$$\skew1\hat{\hat{}} H\too (\hat\K_*\skew3\hat{\hat N}e_\zeta)^*
\eqno £4.11.2\phantom{.}$$
lifting  the action of $H$ on $\hat\K_*\skew3\hat{\hat N}e_\zeta$ and extending the structural 
$\hat\O^*$-group  homomorphism from $\skew3\hat{\hat N}\,.$ It is clear that this homomorphism determines an
irreducible character ${\rm Ext}(\nu)$ of  $\,\skew1\hat{\hat{}} H\,;$ moreover, we set
$H\,\skew1\hat{\hat{}} = \skew2\hat{\hat H} * (\,\skew1\hat{\hat{}} H)^\circ $
and identify $N\i \skew3\hat{\hat N} * \skew3\hat{\hat N}^\circ$ with its canonical image in 
$H\,\skew1\hat{\hat{}}\,.$

\bigskip
\noindent
{\bf Proposition~£4.12.} {\it With the notation above, we have a $\skew2\hat{\hat G}\-$interior algebra isomorphism
$${\hat\K}_*\skew2\hat{\hat G}\,{\rm Tr}_H^G(e_\nu)\cong
{\rm Ind}_{\skew2\hat{\hat H}}^{\skew2\hat{\hat G}} \big(\hat\K_*\skew3\hat{\hat N}e_\nu\otimes_{\hat\K} 
\hat\K_*(H\,\skew1\hat{\hat{}}\,/N)\big)
\eqno 4.12.1\phantom{.}$$
mapping $ye_\nu$ on $1\otimes (ye_\nu\otimes 1)\otimes 1$ for any $y\in \skew3\hat{\hat N}\,.$
In particular, if $\K$ is a splitting field for $\hat G$ then $\hat K$ is a splitting field for 
$H\,\skew1\hat{\hat{}}\,/N$ and the map sending any  $\zeta\in {\rm Irr}_{\hat\K}(H\,\skew1\hat{\hat{}}\,/N)$ to 
$$\chi = {\rm Ind}_{\skew2\hat{\hat H}}^{\skew2\hat{\hat G}}\big({\rm Ext} (\nu)\otimes {\rm Res}(\zeta)\big)
\eqno £4.12.2,$$
where ${\rm Res}(\zeta)$ denotes the corresponding irreducible character of $H\,\skew1\hat{\hat{}}\,,$ 
determines a bijection 
$${\rm Irr}_{\hat\K}(H\,\skew1\hat{\hat{}}\,/N)\cong {\rm Irr}_\K(\hat G,\nu)\eqno £4.12.3\phantom{.}$$ 
fulfilling ${\rm d}(\chi) = {\rm d}(\nu)+{\rm d}(\zeta)\,.$\/}
\eject

\medskip
\noindent
{\bf Proof:} Considering a set $X\i G$ of representatives for $G/H$
and the $G\-$stable pairwise orthogonal  set of idempotents $\{(e_\nu)^{x^{-1}}\}_{x\in X}\,,$ from [9,~2.14.2] we get a 
 $\skew2\hat{\hat G}\-$interior algebra isomorphism
 $${\hat\K}_*\skew2\hat{\hat G}\,\hbox{Tr}_H^G(e_\nu)\cong 
 {\rm Ind}_{\skew2\hat{\hat H}}^{\skew2\hat{\hat G}} \big(e_\nu ({\hat\K}_*\skew2\hat{\hat G})e_\nu\big)
 \eqno £4.12.4\phantom{.}$$
and, as in the proof of [15,~Proposition~3.5], we actually get
$$e_\nu ({\hat\K}_*\skew2\hat{\hat G})e_\nu = {\hat\K}_*\skew2\hat{\hat H} e_\nu
\eqno £4.12.5;$$
moreover, since $\hat\K_*\skew3\hat{\hat N}e_\nu$ is a full matrix algebra over $\hat\K$ and a unitary subalgebra of 
${\hat\K}_*\skew2\hat{\hat H} e_\nu\,,$ the multiplication in this $\hat\K\-$algebra induces
a $\hat\K\-$algebra isomorphism
$$\hat\K_*\skew3\hat{\hat N}e_\nu\otimes_{\hat\K} C\cong {\hat\K}_*\skew2\hat{\hat H} e_\nu
\eqno £4.12.6\phantom{.}$$
where $C$ is the centralizer of $\hat\K_*\skew3\hat{\hat N}e_\nu$ in ${\hat\K}_*\skew2\hat{\hat H}e_\nu$
[7,~Proposition~2.1]; then,  as in the proof of [15,~Proposition~3.2], $C$ becomes an 
$H\,\skew1\hat{\hat{}}\,/N\-$interior algebra and it is easily checked that the structural $\hat\K\-$algebra
homomorphism 
$$\hat\K_*(H\,\skew1\hat{\hat{}}\,/N)\too C
\eqno £4.12.7\phantom{.}$$
is an isomorphism.

\smallskip
In particular, if $\K$ is a splitting field for $\hat G$ then $\K_*\hat G\,\hbox{Tr}_H^G(e_\nu)$ is a direct product
of full matrix algebras over $\K$ and therefore ${\hat\K}_*\skew2\hat{\hat G}\,\hbox{Tr}_H^G(e_\nu)$
 is a direct product of full matrix algebras over $\hat\K\,;$  since $\hat\K_*\skew3\hat{\hat N}e_\nu$ 
 is a full matrix algebra over $\hat\K\,,$ this forces $\hat\K_*(H\,\skew1\hat{\hat{}}\,/N)$ to be also 
 a full matrix algebra over $\hat\K$ and induces bijection~£4.12.3. Finally, setting 
 $\xi = {\rm Ext} (\nu)\otimes {\rm Res}(\zeta)\,,$  we have
$${{\rm d}(\chi)} = \vartheta_p\big({\vert G\vert\over \chi(1)}\big) = \vartheta_p\big({\vert H\vert\over \xi(1)}\big)
=  \vartheta_p\big({\vert N\vert\vert H/N\vert\over \nu(1)\zeta(1)}\big) = {\rm d}(\nu)+{\rm d}(\zeta)
\eqno £4.12.8.$$
We are done.

\bigskip
\bigskip
\noindent
{\bf  £5. The $p$-solvable $\O^*\-$group case}
\bigskip 
£5.1. As above, $\hat G$ is an $\O^*\-$group with finite $\O^*\-$quotient $G$ and in
this section {\it we assume that $G$ is $p$-solvable\/}. Let $b$ be a block of 
$\hat G$ and  $S$  a $G\-$stable unitary subalgebra of~$A$ isomorphic to a {\it direct product 
of full matrix algebras  over  $\O\,,$\/} maximal such that $p$ does not divide the $\O\-$rank 
of its indecomposable factors. Since $b$ is primitive in~$Z(\O_*\hat Gb)\,,$ $G$~acts transitively on
the set~$I$ of primitive idempotents of~$Z(S)$ and we borrow the notations
$i\,,$ $\hat H\,,$ $^{\skew4\hat {}}H$ and~$H^{\skew4\hat {}}\,$ from~£4.4. According to 
Propositions~£4.5 and~£4.8, and to Theorem~£4.9,\break
\eject
\noindent
$i$ is a block of $\hat H$ which belongs to a point $\beta$ of~$\hat H$ on $\O_*\hat G$ and, 
for a suitable block $c$ of~$H^{\skew4\hat{}}\,,$ we have $\hat G$- and $\hat H$-interior
algebra isomorphisms
$$\O_*\hat G b\cong {\rm Ind}_{\hat H}^{\hat G} (\O_*\hat H i)\qq 
(\O_*\hat G)_\beta \cong \O_*\hat Hi\cong Si\otimes_\O \O_*H^{\skew4\hat {}}\,c
\eqno £5.1.1\phantom{.}$$
and an equivalence of categories $\iota\,\colon \L_{\O_*H^{\skew4\hat {}}\,c}\to 
\L_{\O_*\hat G b}\,;$ in particular, there is a defect pointed $\O^*\-$group $\hat P_\gamma$
of $b$ contained in $\hat H_\beta\,.$ Moreover, we denote by ${\Bbb O}_{p'}(\hat H)\,,$
${\Bbb O}_{p'} (\,^{\skew4\hat {}}H)$ and ${\Bbb O}_{p'}(H^{\skew4\hat {}}\,)$ the 
respective inverse images in $\hat H\,,$ $\,^{\skew6\hat {}}H$ and~$H^{\skew4\hat {}}$ of~${\Bbb O}_{p'}(H)\,.$

\bigskip
\noindent
{\bf Proposition £5.2.} {\it With the notation and the hypothesis above, $P$ is a Sylow
$p\-$subgroup of $H\,,$ $Si$ coincides with $\O_* {\Bbb O}_{p'}(\hat H) i$ and is a Dade $P\-$algebra,
 and the inclusion of  ${\Bbb O}_{p'}(\hat H)$ in $(Si)^*$ induces an $H$-stable $\O^*$-group isomorphism 
$\O^*\!\times {\Bbb O}_{p'}(H) \cong {\Bbb O}_{p'}(H^{\skew4\hat {}}\,)$ such that, identifying ${\Bbb O}_{p'}(H)$
with its image, we have
$$c = {1\over\vert {\Bbb O}_{p'}(H)\vert}\sum_{y\in{\bf O}_{p'}(H)}\! y\qq 
\O_*H^{\skew4\hat {}}\,c \cong  \O_*\bigl(H^{\skew4\hat {}}\,/{\Bbb O}_{p'}(H)\big)
\eqno £5.2.1.$$
In particular, we have an equivalence of categories $\L_{\O_*\hat G\,b}\cong 
\L_{\O_*(H^{\skew4\hat {}}\,/{\Bbb O}_{p'}(H))}$ and the $\O\-$algebras $\O_*\hat G\,b$ and 
$\O_*\bigl(H^{\skew4\hat {}}\,/{\Bbb O}_{p'}(H)\big)$ are Morita equivalent.
Moreover, setting $Q = P\cap {\Bbb O}_{p',p}(H)\,,$ $c$ is primitive in $(\O_*H^{\skew4\hat {}}\,c)^Q\,.$\/}
\medskip
\noindent
{\bf Proof:} Note that the Brauer First Main Theorem still holds for blocks of $\O^*\-$groups (cf.~£2.7) 
and that it follows from Proposition~£3.5 above that [10,~Proposition~5.6 and~Corollary~5.8]
apply to our present situation; thus, the proof of [15,~Proposition~4.2] applies. Moreover,  since
$Si$ is generated by a finite $p'\-$group [9,~Lemma~5.5], $P$ stabilizes an $\O\-$basis of $Si$ and, since $p$ 
does not divide ${\rm rank}_\O(Si)\,,$ $P$ fixes some element in any  $P\-$stable  $\O\-$basis; then,
the last statement follows from Proposition~£4.5 and Theorem~£4.9.

\bigskip
\noindent
{\bf Corollary £5.3.} {\it  With the notation and the hypothesis above, denote by  $Q$ the intersection
$P\cap  {\Bbb O}_{p',p}(H)\,,$  by $\hat Q$ the converse image of $Q$ in $\hat P$ and by $\delta$ a
local point of $\hat Q$ on $\O_*\hat G$ such that $\hat Q_\delta\i \hat P_\gamma\i 
\hat H_\beta$ . Then, $\hat Q_\delta$ is the unique local pointed $\O^*\-$group on~$\O_*\hat G$
fulfilling the following conditions
\smallskip
\noindent
{\rm £5.3.1}\quad We have  $\,\hat Q_\delta \,\triangleleft \,\hat P_\gamma\,,$ $C_P (Q) =
Z(Q)$ and ${\Bbb   O}_p\big(E_G(\hat Q_\delta)\big) = \{1\}$
\smallskip
\noindent
{\rm £5.3.2}\quad We have $E_G (\hat R_\varepsilon,\hat P_\gamma) = E_{N_G (\hat Q_\delta)} 
(\hat R_\varepsilon,\hat P_\gamma)$ for any local pointed $\O^*\-$ group $\hat R_\varepsilon$ on
$\O_*\hat G$ contained in~$\hat P_\gamma\,.$
\smallskip
\noindent
Moreover, denoting by $b$ and $f$ the respective blocks of $\hat G$ and $N_{\hat G} (\hat Q_\delta)$ 
determined by $\hat P_\gamma\,,$ the $\O\-$algebras $\O_*\hat G\,b$ and 
$\O_*N_{\hat G} (\hat Q_\delta)f$ are Morita equivalent.\/}

\medskip
\noindent
{\bf Proof:} It follows from Proposition~£5.2 that $Si$ is a {\it Dade $P\-$algebra\/};
hence, Corollary~£4.6 and Theorem~£4.9 above apply, and therefore
the proof of [15, Corollary~4.3] applies too. Moreover, it follows from isomorphisms~£5.1.1 
and Proposition~£5.2 above that $\O_*\hat G\,b$ is  Morita equivalent to 
$\O_*\bigl(H^{\skew4\hat {}}\,/{\Bbb O}_{p'}(H)\big)\,;$ similarly, setting  $\hat N = N_{\hat H} (\hat Q_\delta)$ 
and $N = N_H (\hat Q_\delta)\,,$  $\O_*N_{\hat G} (\hat Q_\delta)f$ is Morita equivalent to 
$\O_*\bigl(N^{\skew4\hat {}}\,/{\Bbb O}_{p'}(N)\big)\,;$ but, the uniqueness of $\delta$ forces 
$N = N_H(\hat Q)$ and thus, by the {\it Frattini argument\/}, we get $H = {\Bbb O}_{p'}(H).N\,,$ so that we still get
$$H/{\Bbb O}_{p'}(H)\cong N/{\Bbb O}_{p'}(N)
\eqno £5.3.3.$$

\smallskip
Finally, we claim that we have an $\O^*\-$group isomorphism
$$H^{\skew4\hat {}}/{\Bbb O}_{p'}(H)\cong N^{\skew4\hat {}}/{\Bbb O}_{p'}(N)
\eqno £5.3.4;$$
indeed, it suffices to prove that $^{\skew4\hat {}}N$ is an $\O^*\-$subgroup of $^{\skew4\hat {}}H\,;$
but, these $\O^*\-$groups come from the respective actions of $N$ and $H$ on $\Bbb O_{p'}(\hat N)$
and $\Bbb O_{p'}(\hat H)$ and therefore  it suffices to prove that $^{\skew4\hat {}}N^k$ is an $k^*\-$subgroup of $^{\skew4\hat {}}H^k\,;$ moreover, it is clear that 
$$\Bbb O_{p'}(\hat N) = C_{\Bbb O_{p'}(\hat H)}(\hat Q)
\eqno £5.3.5.$$
Consequently, the corresponding block $j$ of $\hat N$ fulfills ${\rm Br}_Q (j) = {\rm Br}_Q (i)$
and therefore we have to consider  the respective actions of $N$ and $H$ on $(Si)(Q)$ and $k\otimes_\O Si\,;$
thus, our claim follows from the existence of a splitting (cf.~£4.2)
$$\omega_{(Q,k\otimes_\O Si)} : \hat F_{k\otimes_\O Si}(Q)\too k^*
\eqno £5.3.6.$$

\medskip 
£5.4.   Similarly as in [16], we call {\it Fitting pointed $\O^*\-$group\/} of $\hat G$ any local pointed 
$\O^*\-$group $\hat O_\eta$ on $\O_*\hat G$ fulfilling  conditions~£5.3.1 and~£5.3.2 above 
with respect to some maximal local pointed $\O^*\-$group  $\hat P_\gamma$ on $\O_*\hat G\,.$

\bigskip
\noindent
{\bf Theorem £5.5.} {\it   With the notation and the hypothesis above,   there are a $p\-$solvable finite $\O^*\-$group
$\hat L$ containing $\hat P$ as an $\O^*\-$subgroup and a primitive Dade $P\-$algebra~$T\,,$ both unique up to
isomorphisms, fulfilling the following conditions
\smallskip
\noindent
{\rm £5.5.1}\quad $C_L\big({\Bbb O}_p(L)\big) = Z\big({\Bbb O}_p(L)\big)$ where $L$ is the $\O^*\-$quotient of $\hat L\,.$
\smallskip
\noindent
{\rm £5.5.2}\quad There is a $\hat P\-$interior algebra embedding $e_\gamma\,\colon (\O_*\hat G)_\gamma\too 
T\otimes_\O \O_*\hat L\,.$
\smallskip
\noindent
In particular, $e_\gamma$ induces an equivalence of categories $\L_{\O_*\hat G\,b}\cong \L_{\O_*\hat L}$
and a Morita equivalence between the $\O\-$algebras $\O_*\hat G\,b$ and $\O_*\hat L\,.$
Moreover, for any local pointed  $\O^*\-$group $\hat Q_\delta$ on $\O_*\hat G$ contained in 
$\hat P_\gamma\,,$ we have
$$F_{\O_*\hat G} (\hat Q_\delta,\hat P_\gamma) \i F_T (Q,P)
\eqno £5.5.3\phantom{.}$$
and, denoting by $\hat Q_{\dot\delta}$ the  corresponding local pointed  $\O^*\-$group on $\O_*\hat L\,,$  
any {\it polarization\/} $\omega$ determines a $k^*\-$group isomorphism
$$\hat F_{\O_*\hat G} (\hat Q_\delta) \cong \hat F_{\O_*\hat L} (\hat Q_{\dot\delta}) 
\eqno £5.5.4.$$\/}
\eject

\medskip
\noindent
{\bf Proof:} As in [15,~4.4], the existence of $\hat L\,,$ $T$ and $e_\gamma$
follows from the results above; moreover, it follows from [11,~statement~1.5.2] that, in order to prove the
uniqueness of the {\it Dade $P\-$algebra\/} $T\,,$ it suffices to prove the uniqueness of the {\it Dade $P\-$algebra\/}
$k\otimes_\O T$ over $k$ and this uniqueness follows from [15,~Lemma~4.5 and~Theorem~4.6].

\smallskip

Moreover, from Corollary~£4.6, any local pointed $\O^*\-$group $\hat Q_\delta$ 
on $\O_*\hat G$ contained in $\hat P_\gamma$ determines a local pointed $\O^*\-$group 
$\hat Q_{\dot\delta}$ on~$\O_*\hat L\,,$ and this correspondence is bijective; then, since 
$(\O_*\hat G)_\gamma$ and $\O_*\hat L$  are twisted diagonal 
$\O_*(\hat P\,\hat\times\,\hat P^\circ)\-$modules by left- and right-hand multiplication, 
it follows again from Corollary~£4.6 above and from the uniqueness of the $\hat P_\gamma\-$source pair
$(T,\hat L)$   that we have
$$F_{\O_*\hat G} (\hat Q_\delta,\hat P_\gamma) = F_{\O_*\hat L} (\hat Q_{\dot\delta},\hat P_{\dot\gamma}) 
\i F_T (Q,P)
\eqno £5.5.5.$$

\smallskip
 At this pont, according to [10,~Proposition~5.11], any {\it polarization\/}
$\omega$ determines a $k^*\-$group isomorphism
$$\hat F_{\O_*\hat G} (\hat Q_\delta) \cong \hat F_{\O_*\hat L} (\hat Q_{\dot\delta}) 
\eqno £5.5.6.$$
Finally, isomorphism~£5.6.4 below together with Proposition~£4.5 implies the Morita equivalence between
$\O_*\hat G\,b$ and~$\O_*\hat L\,.$ We are done.

\medskip
£5.6.  Let $b$ be a block of $\hat G$ and $\hat P_\gamma$  a maximal local pointed $\O^*\-$group 
on $\O_*\hat G\,b\,;$ we call {\it $\hat P_\gamma\-$source\/} of $b$  any pair $(T,\hat L)$ 
formed by a {\it primitive Dade $P\-$algebra\/} $T$ --- considered as a $P\-$interior algebra fulfilling 
$${\rm det}_T (P\.{\rm id}_T) = \{1\}
\eqno £5.6.1\phantom{.}$$
 --- and by a $p\-$solvable finite $\O^*\-$group $\hat L$ containing $\hat P$ as an $\O^*\-$subgroup, which fulfills  conditions~£5.5.1 and~£5.5.2 above.
Then, since we have  $\hat P\-$interior algebra embeddings
 $$\O_*\hat L\too T^\circ \otimes_\O T\otimes_\O \O_*\hat L \longleftarrow T^\circ
 \otimes_\O (\O_*\hat G)_\gamma
 \eqno £5.6.2\phantom{.}$$
 and $\hat P$ has a unique local point $(T^\circ\otimes_\O T)\times \dot\gamma$  on 
 $T^\circ\otimes_\O T\otimes_\O \O_*\hat L$ [10, Theorem~5.3], we still have a $\hat P\-$interior algebra embedding
 $$e_\gamma^\circ : \O_*\hat L\too T^\circ \otimes_\O (\O_*\hat G)_\gamma
 \eqno £5.6.3.$$
Note that, according to Corollary~£4.6 and Proposition~£5.2 above,  the uniqueness of the 
$\hat P_\gamma\-$source pair $(T,\hat L)$ forces the existence of an $\O^*\-$group isomorphism
$$\hat L\cong H^{\skew4\hat {}}\,/{\Bbb  O}_{p'}(H)
\eqno £5.6.4\phantom{.}$$
inducing the canonical $\O^*\-$group homomorphism $\hat P\to H^{\skew4\hat {}}\,/{\Bbb O}_{p'}(H)\,,$ 
and of a $P\-$algebra embedding $\,T\to {\rm Res}_P^H (Si)\,.$

 \medskip
 £5.7. Let $\hat Q_\delta$ be a  local pointed $\O^*\-$subgroup on $\O_*\hat G$ contained in $\hat P_\gamma$
 and containing the Fitting pointed $\O^*\-$group $\hat O_\eta$ of $\hat G$ contained in $\hat P_\gamma\,;$\break
 \eject
 \noindent
we know that  there is a unique point $\beta$ of $N_{\hat G}(\hat Q_\delta)$ on $\O_*\hat G$ fulfilling (cf.~£3.11)
 $$\hat Q_\delta\i N_{\hat G}(\hat Q_\delta)_\beta
 \eqno £5.7.1\phantom{.}$$
 and that, up to replacing $\hat P_\gamma$ by a $G\-$conjugate, we may assume that 
 $N_{\hat P} (\hat Q_\delta)_\rho$ is a {\it defect pointed $\O^*\-$group\/} of $N_{\hat G}(\hat Q_\delta)_\beta$
 for a suitable local point $\rho$ of $N_{\hat P} (\hat Q_\delta)$ on $\O_*\hat G\,;$ then, denoting by $f$ the block
 of $C_{\hat G}(\hat Q)$ determined by $\delta\,,$ note that $f$ is still a block of $N_{\hat G}(\hat Q_\delta)$
 and that $N_{\hat P} (\hat Q_\delta)_\rho$ determines a {\it maximal local pointed $\O^*\-$group\/}
 $N_{\hat P} (\hat Q_\delta)_{\bar\rho}$ on $\O_* N_{\hat G}(\hat Q_\delta)f$

 \bigskip
 \noindent
 {\bf Proposition~£5.8.} {\it  With the notation and the hypothesis above, if $(T^{^\delta},\hat L^{^\delta})$
 is a $N_{\hat P} (\hat Q_\delta)_{\bar\rho}\-$source of  the block $f$ of $N_{\hat G}(\hat Q_\delta)\,,$
we have an $N_{\hat P} (\hat Q_\delta)\-$algebra embedding and an $\O^*\-$group isomorphism
inducing the identity on $N_{\hat P} (\hat Q_\delta)$
 $$k\otimes_\O T^{^\delta}\too T (Q)\quad and\quad \hat L^{^\delta}\cong N_{\hat L}(\hat Q)
 \eqno £5.8.1.$$\/}

\par
\noindent
{\bf Proof:} It follows from Proposition~£5.2 that $Si = \O_*\Bbb O_{p'}(\hat H)i$ is a {\it Dade $P\-$algebra\/} 
and therefore the $k\-$algebra (cf.~isomorphism~£3.2.2)
$$(Si)(Q)\cong k_*\Bbb O_{p'}\big(C_{\hat H^k}(\hat Q)\big){\rm Br}_Q (i)
\eqno £5.8.2\phantom{.}$$
is simple [11,~1.8]; thus, denoting by $f^k$ the image of $f$ in $k_*N_{\hat G^k}(\hat Q_\delta)\,,$ it is clear that 
$N_G(\hat Q_\delta)$ stabilizes the semisimple $k\-$subalgebra
$$\prod_x k_*\Bbb O_{p'}\big(C_{\hat H^k}(\hat Q)\big){\rm Br}_Q (i^x)
\eqno £5.8.3\phantom{.}$$
 of $k_*N_{\hat G^k}(\hat Q_\delta)f^k\,,$ where $x\in N_G(\hat Q_\delta)$ runs over a set of representatives
  for~$N_G(\hat Q_\delta)/N_H(\hat Q_\delta)\,.$

\smallskip
Hence, it follows from [15,~Proposition~3.5] that ${\rm Br}_Q (i)$ is a block of $N_{\hat H^k}(\hat Q_\delta)$
over $k$ and therefore it comes from a block $j$ of $N_{\hat H}(\hat Q_\delta)$ (cf.~£3.6); actually, it is easily 
checked that $j$ belongs to $\O_*\Bbb O_{p'}\big(C_{\hat H}(\hat Q)\big)i\i Si\,;$  in particular,
$N_G (\hat Q_\delta)$ stabilizes in $\O_*N_{\hat G}(\hat Q_\delta)f$ the following  direct product
of full matrix algebras over $\O$
$$S^{^\delta} = \prod_x \O_*\Bbb O_{p'}\big(C_{\hat H}(\hat Q)\big)j^x
\eqno £5.8.4,$$
where $x\in N_G(\hat Q_\delta)$ runs over a set of representatives for $N_G(\hat Q_\delta)/N_H(\hat Q_\delta)\,,$
it acts transitively over the set $J$ of primitive idempotents of $Z(S^{^\delta})\,,$ and  $N_H (\hat Q_\delta)$ 
is the stabilizer of $j\in J\,.$

\smallskip
Consequently, setting $\hat K = \Bbb O_{p'}\big(C_{\hat H}(\hat Q)\big)$ and  considering the 
corres-ponding $\O^*\-$groups $^{\skew4\hat {}}N_H (\hat Q_\delta)$ 
and $N_H (\hat Q_\delta)^{\skew4\hat {}}\,\,,$ it follows from Proposition~£5.2\break
 \eject
 \noindent
 that we have 
a $N_{\hat G}(\hat Q_\delta)\-$interior algebra isomorphism
$$\O_*N_{\hat G}(\hat Q_\delta)f\cong {\rm Ind}_{N_{\hat H}(\hat Q_\delta)}^{N_{\hat G}(\hat Q_\delta)}
\big(\O_*\hat Kj\otimes_\O  \O_* (N_H (\hat Q_\delta)^{\skew4\hat {}}\,/K)\big)
\eqno £5.8.5;$$
thus, it follows from~£5.6 above that we have an $N_{\hat P} (\hat Q_\delta)\-$algebra embedding
$T^{^\delta}\to \O_*\hat Kj$ and therefore, since we have a $P\-$algebra embedding $T\to Si\,,$
from isomorphism~£5.8.2 we get $N_{\hat P} (\hat Q_\delta)\-$algebra embeddings
$$k\otimes_\O T^{^\delta}\too (Si)(Q)\longleftarrow T(Q)
\eqno £5.8.6\phantom{.}$$
which induces the $N_{\hat P} (\hat Q_\delta)\-$algebra embedding in~£5.8.1 since the 
$N_{\hat P} (\hat Q_\delta)\-$al-gebra $T^{^\delta}$ is primitive.

\smallskip
Finally, since the $k^*\-$group $^{\skew4\hat {}}N_H (\hat Q_\delta)^k$ comes from the action 
of $N_H (\hat Q_\delta)$ on the simple $k\-$algebra $(Si)(Q)\,,$  the existence of a splitting (cf.~£4.2)
$$\omega_{(Q,k\otimes_\O Si)} : \hat F_{k\otimes_\O Si}(Q)\too k^*
\eqno £5.8.7\phantom{.}$$
shows that $^{\skew4\hat {}}N_H (\hat Q_\delta)^k$ can be identified to the $k^*\-$group 
$N_{\,^{\skew4\hat {}}H^k}(\hat Q_\delta)$ and therefore we get a $k^*\-$group isomorphism
$$\big(N_H (\hat Q_\delta)^{\skew4\hat {}}\;\big)^k\cong N_{(H^{\skew4\hat {}}\,)^k}(\hat Q_\delta)
\eqno £5.8.8;$$
moreover, since both $\O^*\-$groups $N_H (\hat Q_\delta)^{\skew4\hat {}}\,$ and 
$N_{H^{\skew4\hat {}}\,}(\hat Q_\delta)$ contain $N_{\hat P} (\hat Q_\delta)\,,$ and 
$N_P (\hat Q_\delta)$ is a Sylow $p\-$subgroup of $N_H (\hat Q_\delta)\,,$ it follows from
[3,~Ch.~XII,~Theorem~10.1] that these  $\O^*\-$groups are isomorphic, so that we get $\O^*\-$group
isomorphisms
$$\hat L^{^\delta}\cong N_H (\hat Q_\delta)^{\skew4\hat {}}\,/K\cong 
N_{H^{\skew4\hat {}}\,}(\hat Q_\delta)/K \cong N_{\hat L}(\hat Q)
\eqno £5.8.9.$$
We are done.

\medskip
£5.9. Let us denote by ${\rm Out}_{\hat P}\big((\O_*\hat G)_\gamma\big)$
the group of {\it exterior automorphisms\/} (cf. £2.8) of the $\hat P\-$interior algebra 
$(\O_*\hat G)_\gamma\,,$ which is Abelian according to [9,~Proposition~14.9].

\bigskip
\noindent
{\bf Proposition~£5.10.} {\it With the notation above, there are group isomorphisms
$${\rm Out}_{\hat P}\big((\O_*\hat G)_\gamma\big)\cong {\rm Out}_{\hat P} (\O_*\hat L)\cong {\rm Hom} (L,k^*)
\eqno £5.10 .1\phantom{.}$$
mapping $\tilde\sigma\in {\rm Out}_{\hat P}\big((\O_*\hat G)_\gamma\big)$ on an element 
$\dot{\tilde\sigma}\in {\rm Out}_{\hat P} (\O_*\hat L)$ such that, for any $\hat P\-$interior algebra embedding 
$e_\gamma\,\colon (\O_*\hat G)_\gamma\to T\otimes_\O \O_*\hat L$ we have
$$\tilde e_\gamma\circ \tilde\sigma = (\widetilde{\rm id}_T \otimes \dot{\tilde\sigma}) \circ \tilde e_\gamma
\eqno £5.10 .2,$$
and mapping $\zeta\in {\rm Hom} (L,k^*)$ on the exterior class of the $\hat P\-$interior algebra automorphism
of $\O_*\hat L$ sending $\hat y\in \hat L$ to $\zeta (y)\.\hat y$ where $y$ is the image of $\hat y$ in $L\,.$
Moreover, ${\rm Out}_{\hat P}\big((\O_*\hat G)_\gamma\big)$ acts regularly over the set of exterior embeddings from
$(\O_*\hat G)_\gamma$ to $T\otimes_\O \O_*\hat L\,.$\/}

\medskip
\noindent
{\bf Proof:} The proof of [16,~Proposition~5.2] applies.
\eject

\bigskip
\bigskip
\noindent
{\bf  £6.  Charactered  pointed $\cal O^*$-groups}
\bigskip 
£6.1. Let $\hat G$ be an $\O^*$-{\it group\/} with a finite $p\-$solvable $\cal O^*$-{\it quotient\/} $G$ and
 assume that $\K$ is a splitting field for all the $\O^*$-subgroups of $\hat G$. 
The main difference from [16] is that, in our present situation, any $\O^*\-$subgroup $\hat Q$ of~$\hat G$
with a finite $p\-$group $\O^*\-$quotient $Q$ has to be always ``accompanied''  with an {\it irreducible character\/}.

\medskip
£6.2. Thus, let us call {\it charactered $\O^*\-$subgroup\/} of $\hat G$ any pair $(\hat Q,\mu)$ formed by
an $\O^*\-$subgroup $\hat Q$ of~$\hat G$ with a finite $p\-$group $\O^*\-$quotient $Q$ and an 
{\it irreducible $\K\-$character\/} $\mu$ of $\hat Q\,,$ and call {\it defect\/} of $(\hat Q,\mu)$ the
defect of~$\mu$ (cf.~£3.8);  denote by $e_\mu$ the central idempotent 
of ${\K}_*\hat Q$ determined by $\mu$ and by~$N_{\hat G} (\hat Q,\mu)$ the stabilizer of $\mu$ in 
$N_{\hat G} (\hat Q)\,.$

\medskip
£6.3. Now,  as in~£4.11 above, if~$\hat\K$ is an extension of $\K$ of degree divisible 
by~$\vert G\vert$ and $\hat\O$ the ring of integers of $\hat\K\,,$ it follows from Proposition~£2.3 that
 there are an $\hat\O^*$-group $\,\skew3\hat{\hat N}_{G} (\hat Q,\mu)$ with  
 $\hat\O^*\-$quotient~$N_{G} (\hat Q,\mu)$ containing and  normalizing~$\skew3\hat{\hat Q} = \hat Q^{\hat\O}\,,$ 
 and an $\hat\O^*$-group homomorphism 
$$\skew3\hat{\hat N}_{G} (\hat Q,\mu)\too (\hat\K_*\skew3\hat{\hat Q}e_\mu)^*
\eqno £6.3.1\phantom{.}$$
lifting  the action of $N_{G} (\hat Q,\mu)$ on $\hat\K_*\skew3\hat{\hat Q}e_\mu$ and 
extending the structural $\hat\O^*\-$ group  homomorphism from $\skew3\hat{\hat Q}\,;$ 
then,  setting $\skew2\hat{\hat G} = \hat G^{\hat\O}\,,$ the group $Q\i \skew3\hat{\hat Q} * \skew3\hat{\hat Q}^{^\circ}$ becomes a {\it normal\/} subgroup 
of $N_{\skew2\hat{\hat G}} (\hat Q,\mu) * \skew3\hat{\hat N}_{G} (\hat Q,\mu)^\circ$ and we set
$$\skew3\hat{\bar N}_G (\hat Q,\mu) = \big(N_{\skew2\hat{\hat G}} (\hat Q,\mu) * 
\skew3\hat{\hat N}_{G} (\hat Q,\mu)^\circ\big)\big/Q
\eqno £6.3.2\phantom{.}$$
which is an $\hat\O^*\-$group with $\hat\O^*\-$quotient $\bar N_G (\hat Q,\mu) = 
N_{\hat G} (\hat Q,\mu)/\hat Q\,;$ moreover, according to Proposition~£4.12, we have a bijection
$${\rm Irr}_\K (N_{\hat G} (\hat Q),\mu)\cong {\rm Irr}_{\hat\K}\big(\skew3\hat{\bar N}_G (\hat Q,\mu)\big)
\eqno £6.3.3$$
fulfilling ${\rm d}(\zeta) = {\rm d}(\mu) + {\rm d}(\bar\zeta)$ if $\zeta\in {\rm Irr}_\K (N_{\hat G} (\hat Q),\mu)$ 
maps on $\bar\zeta\,.$

\medskip
£6.4. We say that a {\it charactered $\O^*\-$subgroup\/} $(\hat R,\nu)$ of $\hat G$ is {\it normal\/} 
in~$(\hat Q,\mu)$ if $\hat R$ is normal in $\hat Q\,,$ $Q$ stabilizes $\nu\,,$ and  $\nu$ is involved 
in ${\rm Res}_{\hat R}^{\hat Q}(\mu)$ or, equivalently, we have ${\displaystyle{\mu(1)\over \nu(1)}}\.\nu 
= {\rm Res}_{\hat R}^{\hat Q}(\mu)\,;$ then, setting $\skew3\hat{\hat R} = \hat R^{\hat\O}\,,$  
it follows from Propositions~£3.2 and~£4.12 that the action of $Q$ on ${\hat\K}_*\skew3\hat{\hat R}\,e_\nu$ 
determines an $\hat\O^*\-$group $\skew3\hat{\hat Q}^{^\nu}$ of $\hat\O^*\-$quotient $Q\,,$ that $R$ can
be identified to a normal subgroup of $\skew3\hat{\hat Q} *(\skew3\hat{\hat Q}^{^\nu})^\circ$  and that, setting
$$\bar Q = Q/R\qq        \skew3\hat{\bar Q}
= \big(\skew3\hat{\hat Q} *(\skew3\hat{\hat Q}^{^\nu})^\circ\big)/R
\eqno £6.4.1\phantom{.}$$
\eject
\noindent
and denoting by $\mu^{\hat\O}$ the irreducible character of $\skew3\hat{\hat Q}$ determined by $\mu\,,$
there is a unique irreducible character $\bar \mu$ of $ \skew3\hat{\bar Q}$ fulfilling
$$\mu^{\hat\O} = {\rm Ext} (\nu)\otimes {\rm Res}(\bar\mu)\qq {\rm d}(\mu) = {\rm d}(\nu) + {\rm
d}(\bar\mu)
\eqno £6.4.2$$
where ${\rm Ext}(\nu)$ and ${\rm Res}(\bar\mu)$ are the corresponding irreducible characters of 
$\skew3\hat{\hat Q}^{^\nu}$ and $\skew3\hat{\hat Q} *(\skew3\hat{\hat Q}^{^\nu})^\circ\,.$

\bigskip
\noindent
{\bf Proposition~£6.5.} {\it With the notation and the hypothesis above, let $(\hat Q,\mu)$ and $(\hat R,\nu)$
be charactered $\O^*\-$subgroups of $\hat G$ such that $(\hat R,\nu)$ is normal in~$(\hat Q,\mu)\,.$
Set $N  = N_G (\hat R,\nu)$ and $\bar N = N/R\,.$ Denote by~$\skew3\hat{\skew3\hat\K}$ an extension 
of $\hat\K$ of degree divisible by $\bar N$ and  by $\skew3\hat{\skew3\hat\O}$ the ring of integers 
of~$\skew3\hat{\skew3\hat\K}\,.$ Then,  $(\skew3\hat {\bar Q},\bar\mu)$ is a  charactered 
$\hat\O^*\-$subgroup of $\skew3\hat{\bar N}_G (\hat R,\nu)$ and the natural group isomorphism $\bar N_N (\hat Q,\mu)\cong 
\bar N_{\bar N} (\skew3\hat {\bar Q},\bar \mu)$ can be canonically lifted to an
  $\skew3\hat{\skew3\hat\O}\-$group isomorphism
$$\skew3\hat{\bar N}_N (\hat Q,\mu)^{\skew3\hat{\skew3\hat\O}}\cong 
\skew3\hat{\bar N}_{\bar N} (\skew3\hat {\bar Q},\bar \mu)
\eqno £6.5.1.$$\/}

\par
\noindent
{\bf Proof:} Clearly $\skew3\hat {\bar Q}$ is contained in $\skew3\hat{\bar N} =\skew3\hat{\bar N}_G (\hat R,\nu)$
(cf.~£6.3.2); thus,  $(\skew3\hat {\bar Q},\bar\mu)$ is a {\it charactered $\hat\O^*\-$subgroup\/} of $\skew3\hat{\bar N}$ 
 and, according to our definition, we have
$$\skew3\hat{\bar N}_{\bar N} (\skew3\hat {\bar Q},\bar\mu) =  \big(N_{\skew3\hat{\skew3\hat{\bar N}}}
(\skew3\hat {\bar Q},\bar\mu) * \skew3\hat{\hat N}_{\bar N} (\skew3\hat {\bar Q},\bar\mu)^\circ\big)\big/\bar Q
\eqno £6.5.2\phantom{.}$$
where we set $\skew3\hat{\skew3\hat{\bar N}} = \skew3\hat{\bar N}^{\skew3\hat {\skew3\hat\O}}\,.$
But, it follows from Proposition~£4.12 that we have
$$\hat\K_*\skew3\hat{\hat Q} e_\mu\cong 
\hat\K_*\skew3\hat{\hat R} e_\nu\otimes_{\hat\K} \hat\K_*\skew3\hat {\bar Q} e_{\bar\mu}
\eqno £6.5.3\phantom{.}$$
and that $N_N (\hat Q,\mu)$ stabilizes each factor of this tensor product. Consequently, we get a canonical
$\skew3\hat {\skew3\hat\O}^*\-$group isomorphism
$$\skew3\hat {\skew3\hat {N}}_N (\hat Q,\mu)^{\skew3\hat {\skew3\hat\O}}\cong 
N_{\skew3\hat {\skew3\hat {N}}_G (\hat R,\nu)} (\hat Q,\mu)^{\skew3\hat {\skew3\hat\O}}
* \skew3\hat{\hat N}_{N} (\skew3\hat {\bar Q},\bar \mu) 
\eqno £6.5.4,$$
where $\skew3\hat{\hat N}_{N} (\skew3\hat {\bar Q},\bar \mu) $ denotes the corresponding {\it pull-back\/} from
$\skew3\hat{\hat N}_{\bar N}  (\skew3\hat {\bar Q},\bar\mu)\,,$ and therefore we still get
$$\skew3\hat{\hat N}_{\!\bar N} (\skew3\hat {\bar Q},\bar \mu) \cong \big(\skew3\hat {\skew3\hat {N}}_{\!N} (\hat Q,\mu)
* N_{\skew3\hat {\skew3\hat {N}}_{\!G} (\hat R,\nu)} (\hat Q,\mu)^\circ\big)^{\skew3\hat {\skew3\hat\O}}\big/R
\eqno £6.5.5.$$

\smallskip
Moreover, always according to our definition, we have
$$\skew3\hat{\bar N} \cong \big(N_{\skew2\hat{\hat G}} (\hat R,\nu) * 
\skew3\hat{\hat N}_{G} (\hat R,\nu)^\circ\big)\big/R
\eqno £6.5.6\phantom{.}$$
and therefore we still have
$$N_{\skew3\hat{\skew3\hat{\bar N}}}(\skew3\hat {\bar Q},\bar\mu)\cong \big(N_{N_{\skew2\hat{\hat G}}(\hat Q,\mu)} (\hat R,\nu) *  \skew3\hat{\hat N}_{\!N_G (\hat Q,\mu)} (\hat R,\nu)^\circ\big)^{\skew3\hat {\skew3\hat\O}}/R
\eqno £6.5.7.$$
\eject
\noindent
Hence, setting $\skew3\hat{\hat N} = \hat N^{\hat\O}\,,$ we finally obtain
$$\eqalign{\skew3\hat{\bar N}_{\bar N} (\skew3\hat {\bar Q},\bar\mu)&\cong 
\big(N_{N_{\skew2\hat{\hat G}}(\hat Q,\mu)}  (\hat R,\nu) * 
\skew3\hat {\skew3\hat {N}}_{\!N} (\hat Q,\mu)^\circ\big)^{\skew3\hat {\skew3\hat\O}}\big/Q\cr
&= \big(N_{\skew3\hat{\hat N}} (\hat Q,\mu) * 
\skew3\hat {\skew3\hat {N}}_{\!N} (\hat Q,\mu)^\circ\big)^{\skew3\hat {\skew3\hat\O}}\big/Q
= \skew3\hat{\bar N}_N (\hat Q,\mu)^{\skew3\hat{\skew3\hat\O}}\cr}
\eqno £6.5.8.$$
We are done.

\medskip
£6.6. Coherently, let us call {\it charactered weight\/} of $\hat G$ any triple $(\hat Q,\mu,\zeta)$
formed by a {\it charactered $\O^*\-$subgroup\/} $(\hat Q,\mu)$ of $\hat G$ and, considering the ring of integers
$\hat\O$ of an  extension
$\hat\K$ of $\K$ of degree divisible by~$\vert G\vert\,,$  by an irreducible character $\zeta$ of {\it defect zero\/}  
of the $\hat\O^*\-$group $\skew3\hat{\bar N}_G (\hat Q,\mu)\,,$ and we call {\it defect\/} of~$(\hat Q,\mu,\zeta)$ the
defect of~$\mu$ (cf.~£3.8). Some of them will form our set of parameters and, for a  {\it charactered weight\/}
 $(\hat Q,\mu,\zeta)$ of $\hat G\,,$ this depends on the following {\it local points\/} of~$\hat Q$ on $\O_*\hat G$
 determined by $\zeta\,.$

 \medskip
 £6.7. First of all, denoting by $e_\zeta$ the idempotent of  $Z\big(\hat\K_*\skew3\hat{\bar N}_G (\hat Q,\mu)\big)$ determined by $\zeta\,,$ it follows from Proposition~£3.9 that $e_\zeta$ is actually a block of  
 $\hat\O_*\skew3\hat{\bar N}_G (\hat Q,\mu)$ and that  $\hat\O_*\skew3\hat{\bar N}_G (\hat Q,\mu)e_\zeta$
 is a full matrix algebra over $\hat\O\,.$ Moreover, $C_G (\hat Q)$ is contained in $N_G (\hat Q,\mu)$ and acts 
 trivially on $\hat\K_*\skew3\hat{\hat Q}e_\mu\,,$ so that the converse image of $C_G (\hat Q)$ in 
 $\skew3\hat{\hat N}_{G} (\hat Q,\mu)$ is {\it split\/}, and we denote by~$\bar C_{\skew3\hat{\skew3\hat G}} (\hat Q)$ 
the converse image of $\bar C_G (\hat Q)\cong C_G (\hat Q)/C_Q(\hat Q)$ in  $\skew3\hat{\bar N}_G (\hat Q,\mu)\,.$ 
Then, it follows from 
 Proposition~£4.12 that all the irreducible characters of $\bar C_{\skew3\hat{\skew3\hat G}} (\hat Q)$ 
 involved in the restriction of~$\zeta$  to this normal $\hat\O\-$subgroup still have {\it defect zero\/}; 
 hence, it follows again  from Proposition~£3.9 that the 
 $\hat\O\-$subalgebra~$\hat\O_*\bar C_{\skew3\hat{\skew3\hat G}} (\hat Q)e_\zeta$ 
of~$\hat\O_*\skew3\hat{\bar N}_G (\hat Q,\mu)e_\zeta$ is isomorphic to direct product of full 
matrix algebras over~$\hat\O$ and, in particular, there is a block of {\it defect zero\/}
$\bar f$ of $\bar C_{\skew3\hat{\skew3\hat G}} (\hat Q)$ such that $\bar f e_\zeta\not= 0\,.$

\medskip
£6.8. Consequently, since the image $\bar f^k$ of $\bar f$ in $k_*\bar C_{\skew3\hat G^k} (\hat Q)$ 
is a block  of {\it defect zero\/}  of $\bar C_{\skew3\hat G^k} (\hat Q)\,,$ according to isomorphism~£3.2.2,
we get a local point $\delta$ of $\hat Q$ on $\O_*\hat G$ fulfilling
$$(\O_*\hat G)(\hat Q_\delta)\cong k_*\bar C_{\skew3\hat G^k} (\hat Q)\bar f^k
\eqno £6.8.1,$$
and, according to Proposition~£3.9, the local pointed $\O^*\-$group $\hat Q_\delta$ is {\it selfcentralizing\/}.
Moreover, it follows from Proposition~£4.10 that there are an $\hat\O^*\-$group 
$\skew3\hat F_{\O_*\hat G}(\hat Q_{\delta},\mu)$ with an $\hat\O^*\-$quotient equal to the stabilizer of $\mu$
in~$F_{\O_*\hat G}(\hat Q_{\delta})\,,$ and an irreducible character of {\it defect zero\/} $\bar\zeta$
of this $\hat\O^*\-$group such that we have an $\skew3\hat{\bar N}_G (\hat Q,\mu)\-$interior 
algebra isomorphism
$$\hat\O_*\skew3\hat{\bar N}_G (\hat Q,\mu)e_\zeta\cong {\rm Ind}_{\skew3\hat{\bar N}_G (\hat Q_\delta,\mu)}^{\skew3\hat{\bar N}_G (\hat Q,\mu)}\big(\bar C_{\skew3\hat{\skew3\hat G}} (\hat Q)\bar f
\otimes_{\hat\O} \hat\O_*\skew3\hat F_{\O_*\hat G}(\hat Q_{\delta},\mu)e_{\bar\zeta}\big)
\eqno £6.8.2\phantom{.}$$
where $\skew3\hat{\bar N}_G (\hat Q_\delta,\mu)$ denotes the stabilizer of $\bar f$ or, equivalently, of $\delta$ in 
$\skew3\hat{\bar N}_G (\hat Q,\mu)\,.$

\medskip
£6.9. More generally, we call {\it charactered pointed $\O^*\-$group\/} on $\O_*\hat G$ any pair 
$(\hat Q_\delta,\mu)$ formed by  a {\it selfcentralizing\/} pointed $\O^*\-$group $\hat Q_\delta$  (cf.~£3.11)
on ${\cal O}_*\hat G$  and an irreducible $\K\-$character~$\mu$ of~$\hat Q\,,$ and call {\it defect\/} of 
$(\hat Q_\delta,\mu)$ the defect of $\mu$ (cf.~£3.8); denote by 
$N_{\hat G} (\hat Q_\delta,\mu)$ and $F_{\O_*\hat G}(\hat Q_\delta,\mu)$ the respective stabilizers
of $\mu$ in $N_{\hat G}(\hat Q_\delta)$ and $F_{\O_*\hat G} (\hat Q_\delta)\,,$ and by $f$ the block of
$C_{\hat G} (\hat Q)$ determined by $\delta\,.$ Since $\hat Q_\delta$ is selfcentralizing, it follows from
Proposition~£3.9 that the image $\bar f$ of $f$ in $\hat\O_*\bar C_{\skew3\hat{\skew3\hat G}}(\hat Q)$
is a block of {\it defect zero\/} and therefore, considering the blocks of $\skew3\hat{\bar N}_G (\hat Q,\mu)$
involved in ${\rm Tr}_{\bar C_G (\hat Q)}^{\bar N_G (\hat Q,\mu)}(\bar f)\,,$ we can claim as above Proposition~£4.10 
in order to exhibit an $\hat\O\-$group $\skew3\hat F_{\O_*\hat G}(\hat Q_{\delta},\mu)\,;$ let us give
an alternative definition of this $\hat\O\-$group.

\medskip
£6.10. Let $\hat P_\gamma$ be a {\it maximal local pointed $\O^*\-$group\/} on  $\O_*\hat G$ 
containing~$\hat Q_\delta\,;$ as in~£5.7 above, we may assume that $N_{\hat P}(\hat Q_\delta)_{\bar\rho}$ 
is a {\it maximal local pointed $\O^*\-$group\/} on  $\O_*N_{\hat G} (\hat Q_\delta)f$ and then
consider a {\it $N_{\hat P}(\hat Q_\delta)_{\bar\rho}\-$source\/}  $(S^{^\delta},\hat L^{^\delta})$ of 
the block $f$ of~$N_{\hat G} (\hat Q_\delta)\,;$ thus, $\hat L^{^\delta}$ contains $N_{\hat P}(\hat Q_\delta)$
 which contains $\hat Q\,,$ and it actually follows 
from isomorphism~£5.6.4 that $\hat L^{^\delta}$ normalizes   $\hat Q$ and that  
$\K$ is also a splitting field for all the $\O^*$-subgroups  of~$\hat L^\delta\,.$ 
On the one hand, Proposition~£2.3 applies to $N_{\hat L^\delta}(\hat Q,\mu)\,,$ 
$\hat Q\,,$ $\mu$ and $\hat\K\,,$ and as in~£6.3.3 above  we get a $\hat\O^*\-$group 
$\skew3\hat{\bar N}_{\!L^{^\delta}} (\hat Q,\mu)$ with $\hat\O^*\-$quotient 
$\bar N_{L^\delta}(\hat Q,\mu)$ where $L^{^\delta}$ denotes the $\O^*\-$quotient 
of~$\hat L^{^\delta}\,.$ On the other hand, since we have
$$Q\i {\Bbb O}_p (L^{^\delta})\i N_P(\hat Q_\delta)\qq C_P(Q)\i Q
\eqno £6.10.1,$$
any $p'\-$subgroup of $L^{^\delta}$ which centralizes $Q$ still centralizes ${\Bbb O}_p (L^{^\delta})$
[5,~Ch.~5, Theorem~3.4] and therefore, according to condition~£5.5.1, it is trivial; hence, since 
$N_P(\hat Q_\delta)$ is a  Sylow $p\-$subgroup of $L^{^\delta}\,,$ we still have
$$C_L (Q) = Z(Q)
\eqno £6.10.2.$$
In particular, $\hat Q$ has a unique local point $\dot\delta$ on $\O_*\hat L^{^\delta}$ (cf.~isomorphism~£3.2.2)
and actually we have $\dot\delta = \{1\}$ [14,~1.19].

\medskip
£6.11. Consequently, from Proposition~£3.10 we get
$$\bar N_{\!L^{^\delta}} (\hat Q,\mu) = N_{L^\delta}(\hat Q,\mu)/Q = E_{L^{^\delta}} (\hat Q_{\dot\delta},\mu)
= F_{\O_*\hat L^{^\delta}}(\hat Q_{\dot\delta},\mu) 
\eqno £6.11.1;$$
but, it follows from Theorem~£5.5 that the choice of a $\hat P\-$interior algebra embedding
$$e_{\bar\gamma} : \big(\O_*N_{\hat G} (\hat Q_\delta)\big)_{\bar\gamma}\too S^{^\delta}\!\otimes_\O 
\O_*\hat L^{^\delta}
\eqno £6.11.2\phantom{.}$$
induces an equivalence of categories $\L_{\O_*N_{\hat G} (\hat Q_\delta)f}\cong \L_{\O_*\hat L^{^\delta}}\,;$ 
consequently, denoting by $\hat Q_{\bar\delta}$ the local pointed  $\O^*\-$group on  
$\O_*N_{\hat G} (\hat Q_\delta)f$  determined by~$\hat Q_\delta\,,$ we easily get the equalities
(cf.~Proposition~£3.10)
$$F_{\O_*\hat L^{^\delta}}(\hat Q_{\dot\delta},\mu) = F_{\O_*N_{\hat G} (\hat Q_\delta)}
(\hat Q_{\bar\delta},\mu) = E_G (\hat Q_{\delta},\mu)  = F_{\O_*\hat G}(\hat Q_{\delta},\mu)
\eqno £6.11.3.$$
Thus,   the $\hat\O^*\-$group $\skew3\hat{\bar N}_{\!L^{^\delta}} (\hat Q,\mu)$ 
has $\hat\O^*\-$quotient $F_{\O_*\hat G}(\hat Q_{\delta},\mu)$ and, since this $\hat\O^*\-$group 
contains $N_{\hat P}(\hat Q_\delta)\,,$ up to the choice of a {\it polarization\/} $\omega$ 
it easily follows from~isomorphism~£5.5.4, from Proposition~£5.10 
and from the uniqueness of the {\it $N_{\hat P}(\hat Q_\delta)_{\bar\rho}\-$source\/}  
$(S^{^\delta},\hat L^{^\delta})$ that this $\hat\O^*\-$group is independent of our choices 
up to {\it unique $\hat\O^*\-$group exterior isomorphisms\/}; hence, it makes sense to define
$$\skew3\hat F_{\O_*\hat G}(\hat Q_{\delta},\mu) = \skew3\hat{\bar N}_{\!L^{^\delta}} (\hat Q,\mu)
\eqno £6.11.4;$$
note that the $\hat\O^*\-$group isomorphism~£5.6.4 shows that this definition agrees with the definition above.

\medskip
£6.12. Moreover, it still follows from Theorem~£5.5 that the choice of $e_{\bar\gamma}$ induces 
a Morita equivalence between the $\O\-$algebras $\O_*N_{\hat G} (\hat Q_\delta)f$ and 
$\O_*\hat L^{^\delta}\,,$ and therefore we clearly get a bijection between the sets 
${\rm Irr}_\K\big(N_{\hat G} (\hat Q_\delta),\mu\big)$ and ${\rm Irr}_\K (\hat L^{^\delta},\mu)$
preserving the {\it defect\/} (cf.~£3.8.1); but, in the  present situation  bijection~£6.3.3 yields a new bijection
$${\rm Irr}_\K (\hat L^{^\delta},\mu)\cong {\rm Irr}_{\hat\K}\big(\skew3\hat{\bar N}_{\!L^{^\delta}} (\hat Q,\mu)\big)
\eqno £6.12.1;$$
finally, with the definition above, we get a bijection
$$\Gamma^\omega_{(\hat Q_\delta,\mu)} : {\rm Irr}_\K\big(N_{\hat G} (\hat Q_\delta),\mu\big) \cong 
{\rm Irr}_{\hat\K} \big(\skew3\hat F_{\O_*\hat G}(\hat Q_{\delta},\mu)\big)
\eqno £6.12.2\phantom{.}$$
fulfilling $${\rm d}(\xi) = {\rm d}(\mu) + {\rm d}\big(\Gamma^\omega_{(\hat Q_\delta,\mu)}(\xi)\big)
\eqno £6.12.3\phantom{.}$$
for any $\xi\in {\rm Irr}_\K\big(N_{\hat G} (\hat Q_\delta),\mu\big)\,;$ once again, up to the 
choice of a {\it polarization\/} $\omega$ it easily follows from~isomorphism~£5.5.4, from Proposition~£5.10  
and from the uniqueness of the {\it $N_{\hat P}(\hat Q_\delta)_{\bar\rho}\-$source\/}  
$(S^{^\delta},\hat L^{^\delta})$ that this bijection is independent of our choices.

\medskip
£6.13. More precisely, denote by $b$ the block of $\hat G$ determined by $\hat Q_\delta\,;$ 
in the case that $\hat Q_\delta$ is a {\it Fitting pointed $\O^*\-$group\/} (cf.~£5.4), 
it follows from Corollary~£5.3 that the $\O\-$algebras $\O_*\hat G\,b$
and $\O_*N_{\hat G} (\hat Q_\delta)f$ are Morita equi-valent and therefore, choosing a set of representatives
${\rm Fit}(\hat G)$ for the set of\break
\eject
\noindent
 $G\-$conjugacy classes of {\it charactered Fitting pointed $\O^*\-$groups\/} of
$\hat G\,,$ it is easily checked that we get bijection
$$\Gamma^\omega_{\hat G} : {\rm Irr}_\K (\hat G)\cong \bigsqcup_{(\hat O_{\eta},\nu)\in {\rm Fit}(\hat G)} 
{\rm Irr}_{\hat\K}  \big(\skew3\hat F_{\O_*\hat G}(\hat O_{\eta},\nu)\big)
\eqno £6.13.1\phantom{.}$$
fulfilling 
$${\rm d}(\chi) = {\rm d}(\nu) + {\rm d}\big(\Gamma^\omega_{\hat G}(\chi)\big)
\eqno £6.13.2\phantom{.}$$
for any  $\chi\in {\rm Irr}_\K (\hat G)$ such that $\Gamma^\omega_{\hat G}(\chi)$ belongs to 
${\rm Irr}_{\hat\K}  \big(\skew3\hat F_{\O_*\hat G}(\hat O_{\eta},\nu)\big)\,.$ Note that if $(\hat O_{\eta},\nu)$
and $(\hat O'_{\eta'},\nu')$ are $G\-$conjugate {\it charactered Fitting pointed\break
 $\O^*\-$groups\/}
then any element $x\in G$ fulfilling $(\hat O_{\eta},\nu) = (\hat O'_{\eta'},\nu')^x$ induces a  {\it unique $\hat\O^*\-$group exterior isomorphism\/}
$$\skew3\hat F_{\O_*\hat G}(\hat O_{\eta},\nu)\cong \skew3\hat F_{\O_*\hat G}(\hat O'_{\eta'},\nu')
\eqno £6.13.3\phantom{.}$$
and therefore it induces the same bijection
$${\rm Irr}_{\hat\K}  \big(\skew3\hat F_{\O_*\hat G}(\hat O'_{\eta'},\nu')\big)\cong 
{\rm Irr}_{\hat\K}  \big(\skew3\hat F_{\O_*\hat G}(\hat O_{\eta},\nu)\big)
\eqno £6.13.4\phantom{.}$$
which allows us to identify to each other both sets.

\medskip
£6.14. On the other hand, we have proved that any {\it charactered weight\/} $(\hat Q,\mu,\zeta)$ of 
$\hat G$ determines an $\bar N_G(\hat Q,\mu)\-$orbit of {\it local points\/} $\delta$ of $\hat Q$ 
on $\O_*\hat G\,b$ for a suitable block $b$ of~$\hat G$ --- in such a way  that $\hat Q_\delta$ 
is {\it selfcentralizing\/} --- and  an {\it irreducible  $\hat\K\-$character\/} of {\it defect zero\/} 
$\bar\zeta$ of $\skew3\hat F_{\O_*\hat G}(\hat Q_{\delta},\mu)\,.$
Conversely, for any block $b$ of $\hat G\,,$ let us call {\it charactered $b\-$weight\/} of $\hat G$ 
any triple $(\hat Q_\delta,\mu,\bar\zeta)$ formed by a {\it selfcentralizing\/} pointed $\O^*\-$group 
on $\O_*\hat G\,b\,,$ by an irreducible $\K\-$character $\mu$ of $\hat Q$ and by an irreducible 
$\hat\K\-$character of {\it defect zero\/} $\bar\zeta$ of~$\skew3\hat F_{\O_*\hat G}(\hat Q_{\delta},\mu)\,;$ 
then, it follows from Propositions~£4.10 and~£4.12 that isomorphism~£6.8.2 determines an irreducible 
$\hat\K\-$character of {\it defect zero\/} $\zeta$ of the $\hat\O^*\-$group  
$\skew3\hat{\bar N}_G (\hat Q,\mu)\,,$ so that the triple $(\hat Q,\mu,\zeta)$ becomes a 
{\it charactered weight\/}  of $\hat G\,.$

\medskip
£6.15. Actually, these correspondences define a bijection between the set of $G\-$conjugacy classes 
${\rm Chw}_\K (\hat G)$ of {\it charactered weights\/}  of $\hat G$ and the union of sets of 
 $G\-$conjugacy classes ${\rm Chw}_\K (\hat G,b)$ of {\it charactered $b\-$weights\/}   of~$\hat G$ when $b$ runs 
over the set of blocks of $\hat G\,.$ Let $b$ be a block of~$\hat G\,,$ $\hat O_\eta$ 
a Fitting pointed $\O^*\-$group on $\O_*\hat G\,b$ and $\nu$ an irreducible character of~$\hat O\,;$ 
let us denote by ${\rm Chw}_\K (\hat G,\hat O_\eta,\nu)$ the subset of~${\rm Chw}_\K (\hat G,b)$ 
determined by the {\it charactered $b\-$weights\/} $(\hat Q_\delta,\mu,\bar\zeta)$ of~$\hat G$ fulfilling
$$\hat O_\eta\i \hat Q_\delta \qq
{\rm Res}_{\hat O}^{\hat Q}(\mu) = {\mu(1)\over \nu(1)}\.\nu
\eqno£6.15.1;$$
then, $(\hat O,\nu)$ is {\it normal\/} in $(\hat Q,\mu)$ and in~£6.4 above we have defined a 
{\it charactered $\hat\O^*\-$group\/} $(\skew3\hat{\bar Q},\bar \mu)\,.$
\eject

\bigskip
\noindent
{\bf Proposition~£6.16.} {\it With the notation and the hypothesis above, $(\skew3\hat{\bar Q},\bar \mu)$ 
is a charactered $\hat\O^*\-$subgroup of  $\hat F_{\O_*\hat G}(\hat O_\eta,\nu)\,,$ and any polarization 
$\omega$ determines an $\skew3\hat{\hat\O}^*\-$group isomorphism
$$\hat F_{\O_*\hat G}(\hat Q_\delta,\mu)^{\skew3\hat{\hat\O}}\cong 
\skew3\hat{\bar N}_{F_{\O_*\hat G}(\hat O_\eta,\nu)} (\skew3\hat{\bar Q},\bar \mu)
\eqno £6.16.1.$$
In particular, denoting by $\skew3\bar{\bar\zeta}$ the character of the right-hand member determined 
by $\bar\zeta\,,$ the triple $(\skew3\hat{\bar Q},\bar \mu,\skew3\bar{\bar\zeta})$ is a charactered weight of 
$\hat F_{\O_*\hat G}(\hat O_\eta,\nu)\,.$ Moreover, the correspondence mapping $(\hat Q_\delta,\mu,\bar\zeta)$
on $(\skew3\hat{\bar Q},\bar \mu,\skew3\bar{\bar\zeta})$ induces a bijection
$$\Delta^\omega_{(\hat O_\eta,\nu)} : {\rm Chw}_\K (\hat G,\hat O_\eta,\nu)\cong {\rm Chw}_{\hat\K}\big (\hat F_{\O_*\hat G}(\hat O_\eta,\nu)\big)
\eqno £6.16.2.$$\/}

\par
\noindent
{\bf Proof:} Let $\hat P_\gamma$ be a {\it maximal local pointed $\O^*\-$group\/} on $\O_*\hat G\,b$
containing $\hat Q_\delta$ and $(T,\hat L)$ a {\it $P_\gamma\-$source\/} of $b\,;$ according to Proposition~£6.5,
in particular $(\skew3\hat{\bar Q},\bar \mu)$ is a charactered $\hat\O^*\-$subgroup of $\skew3\hat{\bar N}_L (\hat O,\nu)\,.$
But, denoting by $g$ the block of $N_{\hat G} (\hat O_\eta,\nu)\,,$ by~$\hat P_{\bar\gamma}$
the maximal local pointed $\O^*\-$group on $\O_*N_{\hat G} (\hat O_\eta,\nu)g$ determined by~$\hat P_\gamma$
and by $(T^{^\eta},\hat L^{^\eta})$ a $\hat P_{\bar\gamma}\-$source of $g\,,$ it follows from Proposition~£5.8
that we may assume that $\hat L^{^\eta} = \hat L$ and then, by its very definition (cf.~£6.11.4), we get
$$\hat F_{\O_*\hat G}(\hat O_\eta,\nu) = \skew3\hat{\bar N}_L (\hat O,\nu)
\eqno £6.16.3.$$

\smallskip
Consequently, $(\skew3\hat{\bar Q},\bar \mu)$ is a charactered $\hat\O^*\-$subgroup of 
$\hat F_{\O_*\hat G}(\hat O_\eta,\nu)$ and from Proposition~£6.5 we obtain an 
$\skew3\hat{\skew3\hat\O}^*\-$group isomorphism
$$\skew3\hat{\bar N}_{N_L (\hat O,\nu)} (\hat Q,\mu)^{\skew3\hat{\skew3\hat\O}}\cong 
\skew3\hat{\bar N}_{\bar N_L (\hat O,\nu)} (\skew3\hat {\bar Q},\bar \mu)
\eqno £6.16.4;$$
moreover, with the notation and the hypothesis in~£5.7 above,  it follows from Proposition~£5.8
that we may assume that $\hat L^{^\delta} = N_{\hat L}(\hat Q)$ and then, by its very definition (cf.~£6.11.4), we get
$$\hat F_{\O_*\hat G}(\hat Q_\delta,\mu) = \skew3\hat{\bar N}_L (\hat Q,\mu)
\eqno £6.16.5;$$
finally, since $\mu$ determines $\nu\,,$ we still have $\skew3\hat{\bar N}_L (\hat Q,\mu) 
=\skew3\hat{\bar N}_{N_L (\hat O,\nu)} (\hat Q,\mu)$ and isomorphism~£6.16.1 follows from
isomorphism~£6.16.4 above.

\smallskip
Moreover, it is easily checked that this correspondence is compatible with $G\-$ and 
$F_{\O_*\hat G}(\hat O_\eta,\nu)\-$conjugation and therefore it induces a map
$$\Delta^\omega_{(\hat O_\eta,\nu)} : {\rm Chw}_\K (\hat G,\hat O_\eta,\nu)\too
 {\rm Chw}_{\hat\K}\big (\hat F_{\O_*\hat G}(\hat O_\eta,\nu)\big)
\eqno £6.16.6;$$
in order to prove that it is bijective, let us define the inverse map. If 
$(\skew3\hat{\bar Q}',\bar \mu',\skew3\bar{\bar\zeta}')$ is a {\it charactered weight\/} of 
$\hat F_{\O_*\hat G}(\hat O_\eta,\nu) = \skew3\hat{\bar N}_L (\hat O,\nu)$ then the converse image\break
\eject 
\noindent
$Q'$ of $\bar Q'$ in $L$ fulfills $O\i Q'\i P$ and determines an $\O^*\-$group $\hat Q'$ such that
$\hat O\i \hat Q'\i \hat P\,;$ since $\hat O_\eta$ is selfcentralizing and contained in $\hat P_\gamma\,,$
it is easily checked that we still have $\hat O_\eta\i \hat Q'_{\delta'}\i \hat P_\gamma$ for a {\it unique\/}
local point $\delta'$ of $\hat Q'$ on $\O_*\hat G\,.$

\smallskip
On the other hand, since $\K$ is a splitting field for $\hat Q'\,,$ it follows from Proposition~£4.12 
that there is an irreducible $\K\-$character $\mu'$ of $\hat Q'$ fulfilling (cf.~£6.4.2)
$$\mu'^{\hat\O} = {\rm Ext} (\nu)\otimes {\rm Res}(\bar\mu')
\eqno £6.16.7$$
and then, since $\hat\K$ is a splitting field for $\hat F_{\O_*\hat G}(\hat Q'_{\delta'},\mu')\,,$
there is an irreducible $\hat\K\-$character $\bar\zeta'$ determining $\skew3\bar{\bar\zeta}'$
{\it via\/} the $\skew3\hat{\hat\O}^*\-$group isomorphism~£6.16.1. Thus, 
$(\skew3\hat{\bar Q}',\bar \mu',\skew3\bar{\bar\zeta}')$ determines the {\it charactered $b\-$weight\/}
$(\hat Q'_{\delta'},\mu',\bar\zeta')$ of $\hat G$ and it is quite clear that this correspondence induces the
announced inverse map. We are done.

\bigskip
\bigskip
\noindent
{\bf £7. The charactered Fitting sequences}

\bigskip
£7.1  Let $\hat G$ be an $\O^*$-{\it group\/} with a finite $p\-$solvable $\cal O^*\-${\it quotient\/} $G$ and
 assume that $\K$ is a splitting field for all the $\O^*\-$subgroups of $\hat G$.  
 In order to exhibit bijections between the sets of  irreducible $\K\-$characters of $\hat G$ and of  
$G\-$conjugacy classes of the  {\it inductive weights\/} of~$\hat G$ defined below, we need a third set ---
namely the set of $G\-$conjugacy classes of {\it charactered Fitting sequences\/} of~$\hat G$ --- which
depends on the choice of a {\it polarization\/} $\omega\,.$
We call  {\it charactered Fitting $\omega\-$sequence\/} of~$\hat G$ any sequence 
$$\B = \{(\K_n,\hat G_n,\hat O^n_{\,\eta_n},\nu_n)\}_{n\in \Bbb N}
\eqno £7.1.1\phantom{.}$$
 of quadruples formed by a field extension $\K_n$ of $\K\,,$ by an $\O^*_n\-$group $\hat G_n$ 
 where $\O_n$ denotes the ring of integers of $\K_n\,,$ by a {\it Fitting pointed $\O_n^*\-$group\/} of $\hat G_n$
(cf.~£5.4), and by an irreducible $\K_n\-$character $\nu_n$ of $\hat O^n\,,$ in such a way that 
$\K_0 = \K\,,$ that $\hat G_0 = \hat G$  and that, for any $n\in \Bbb N\,,$ $\K_{n+1}$ is the field extension of degree 
$\vert G_n\vert$ of $\K_n$ and we have (cf.~definition~£6.11.4)
$$\hat G_{n+1} = \hat F_{{\O_n}_*\hat G_{n}}(O^n_{\,\eta_n},\nu_n)
\eqno £7.1.2.$$ 
Note that $\vert G_{n+1} \vert\le \vert G_{n} \vert$ and actually we have $\vert G_{n+1} \vert= \vert G_{n} \vert$
 if and only if~$O^n = \{1\}$ (cf.~condition~£5.3.1). Moreover, for any $h\in \Bbb N\,,$ the sequence 
$$\B_h = \{(\K_{h+n},\hat G_{h +n},\hat O^{h+n}_{\,\eta_{h+N}},\nu_{h+ n})\}_{n\in \Bbb N}
\eqno £7.1.3\phantom{.}$$
 is clearly a   {\it charactered Fitting $\omega\-$sequence\/} of~$\hat G_h\,.$

\medskip
£7.2. If  $\hat G'$ is an $\O^*\-$group isomorphic to $\hat G$ and $\theta\,\colon \hat G\cong \hat G'$
an $\O^*\-$group isomorphism, it is quite clear that, from any  {\it charactered Fitting $\omega\-$sequence\/}  
$$\B = \{(\K_n,\hat G_n,\hat O^n_{\,\eta_n},\nu_n)\}_{n\in \Bbb N}
\eqno £7.2.1\phantom{.}$$ 
of~$\hat G\,,$ we
are able to construct a {\it charactered Fitting $\omega\-$sequence\/} 
$$\B' = \{(\K_n,\hat G'_n,\hat O'^n_{\,\eta'_n},\nu'_n)\}_{n\in \Bbb N}
\eqno £7.2.2\phantom{.}$$
 of~$\hat G'$ inductively defining a sequence of $\O_n^*\-$group isomorphisms~$\theta_n\,\colon 
 \hat G_n\cong \hat G'_n$ by $\theta_0 = \theta$ and, for any  $n\in \Bbb N\,,$ by (cf.~£2.14)
$$\matrix{\hat F_{\tilde\theta_n} (\hat O^n_{\,\eta_n},\nu_n)\hskip-7pt &:  \hat F_{{\O_n}_*\hat G_{n}}(\hat O^n_{\,\eta_n},\nu_n)&\cong &\hat F_{{\O_n}_*\hat G'_{n}}(\hat O'^n_{\,\eta'_n},\nu'_n)\cr
\Vert &\Vert&\phantom{\big\uparrow}&\Vert\cr
\theta_{n+1}&\hat G_{n+1}&\cong &\hat G'_{n+1}\cr}
\eqno £7.2.3\phantom{.}$$
where we sill denote by $\theta_n : \O_*\hat G_n\cong \O_*\hat G'_n$ the corresponding $\O\-$algebra 
isomorphism and  set 
$$\hat O'^n = \theta_n (\hat O^n)\quad,\quad \eta'_n = \theta_n(\eta_n)\qq \nu'_n = {\rm res}_{\theta_n^{-1}}(\nu_n)
\eqno £7.2.4.$$
In particular, the group of inner automorphisms of $G$ acts on the set of 
{\it charactered Fitting $\omega\-$sequences\/} of~$\hat G$ and we denote by  ${\rm Chs}_\K^{_\omega}(\hat G)$ the set 
of ``$G\-$conjugacy classes'' of the {\it charactered Fitting $\omega\-$sequences\/} of~$\hat G\,.$

\medskip
£7.3. If $\B = \{(\K_n,\hat G_n,\hat O^n_{\,\eta_n},\nu_n)\}_{n\in \Bbb N}$ is a {\it charactered Fitting 
$\omega\-$sequence\/} of $\hat G\,,$ it is clear that $\hat O^n_{\,\eta_n}$ determines 
a block $b_n$ of  $\hat G_n$ for any $n\in \Bbb N$ and we call {\it defect\/} of $\B$ the sum
$${\rm d}(\B) = \sum_{n\in \Bbb N} {\rm d}(\nu_n)
\eqno £7.3.1\phantom{.}$$
which makes sense since ${\rm d}(\nu_m) = 0$ for $m$ big enough; as in [16,~6.4], let us call 
{\it irreducible character   $\omega\-$sequence associated to $\B$\/} any sequence  $\{\chi_n\}_{n\in \Bbb N}$ where $\chi_n$ 
belongs to ${\rm Irr}_{\K_n} (\hat G_n,b_n)$ in such a way that, up to identifications, we have (cf.~£6.13.1)
$$\Gamma_{\hat G_n}^\omega (\chi_n) = \chi_{n+1}
\eqno £7.3.2\phantom{.}$$
 for any $n\in \Bbb N\,;$ in this case, it easily follows from equality~£6.13.2 that
 $${\rm d}(\chi_0) = {\rm d}(\B)
 \eqno £7.3.3.$$

 \bigskip
 \noindent
 {\bf Theorem~£7.4.} {\it With the notation and the choice above, any charactered Fitting $\omega\-$sequence 
 $\B = \{(\K_n,\hat G_n,\hat O^n_{\,\eta_n},\nu_n)\}_{n\in \Bbb N}$ of~$\hat G$ admits a unique  
 irreducible character   $\omega\-$sequence   $\{\chi_n\}_{n\in \Bbb N}$ associated to $\B\,.$ Moreover, 
 the correspondence mapping $\B$ to $\chi_0$ induces a natural bijection
$${\rm Chs}_\K^\omega (\hat G)\cong {\rm Irr}_\K (\hat G)
\eqno £7.4.1\phantom{.}$$
which preserves the defect.\/}
\eject

\medskip
\noindent
{\bf Proof:} Since the sequence $\{\hat G_n\}_{n\in \Bbb N}$ ``stabilizes'', we can argue by induction 
on the minimal $n\in \Bbb N$ fulfilling $\hat G_{n+1} = \hat G_n\,.$ If $n = 0$ then $\hat O^0 =\O^*\,,$ 
$\nu_0 = {\rm id}_{\O^*}$ and the block $b_0$ of $\hat G_0$ has {\it defect zero\/}, so that ${\rm Irr}_k (\hat G_0,b_0)$ 
has a  unique element $\chi_0$ (cf.~Proposition~£3.9) and, setting $\chi_n = \chi_0$ for any $n\in \Bbb N\,,$ we get  an 
{\it irreducible character $\omega\-$sequence\/}
associated to $\B\,.$

\smallskip
If $n\not= 0$ then, according to the induction hypothesis, the {\it charactered  Fitting $\omega\-$sequence\/} 
 $\B_1 = \{(\K_{1+n},\hat G_{1 +n},\hat O^{1+n}_{\,\eta_{1+N}},\nu_{1+ n})\}_{n\in \Bbb N}$ of~$\hat G_1$ 
 already admits an {\it irreducible  character   $\omega\-$sequence\/}  
 $\{\chi_{1 +n}\}_{n\in \Bbb N}\,;$ thus, in order to get an {\it irrducible  character   $\omega\-$sequence\/} 
 associated to $\B\,,$ up to identifications it suffices to define (cf.~£7.3.2)
 $$\chi_0 = (\Gamma_{\hat G_0}^\omega)^{-1}(\chi_{1})
 \eqno £7.4.2.$$

\medskip
On the other hand, since the maps $\Gamma_{\hat G_n}^\omega$are bijective, equality~£7.3.2 shows that an 
{\it irreducible character $\omega\-$sequence\/} $\{\chi_n\}_{n\in \Bbb N}$ associated to $\B$ is 
uniquely determined by one of their terms; but, for $n$ big enough, we know that
$b_n$ is a block of {\it defect zero\/} of $\hat G_n$ and then $\chi_n$ is uniquely determined; 
consequently, $\{\chi_n\}_{n\in \Bbb N}$ is uniquely determined and it is quite clear that, up to identifications, 
it only depends on the $G\-$conjugacy class of $\B\,;$ thus,  we have obtained a {\it natural\/} map
$${\rm Chs}^\omega_\K(\hat G)\too {\rm Irr}_\K (\hat G)
\eqno £7.4.3\phantom{.}$$
which preserves the defect (cf.~equality~£7.3.3).

\smallskip
We claim that it is bijective; actually, we will define the inverse map. For~any $\chi\in {\rm Irr}_\K (\hat G)\,,$
we inductively define a sequence $\{\chi_n\}_{n\in \Bbb N}$ in the fol-lowing way; we set  
$\K_0 = \K\,,$ $\hat G_0 = \hat G$  and $\chi_0 = \chi\,,$ we denote by $b_0$ the block of $\chi\,,$ 
and we choose a {\it Fitting pointed $\O^*\-$group\/} $\hat O^0_{\,\eta_0}$ on $\O_*\hat G\,b$ 
and an irreducible character~$\nu_0$  of~$\hat O^0$ involved in the irreducible character of 
$N_{\hat G}(\hat O^0_{\,\eta_0})$ determined  by $\chi$ (cf.~Corollary~£5.3); moreover, for any  
$n\in \Bbb N\,,$  we denote by $\K_{n+1}$ the field 
extension of degree $\vert G_n\vert$ of $\K_n\,,$ we set
$$\hat G_{n+1} = \hat F_{k_*\hat G_{n}}(O^n_{\,\eta_n})\qq \chi_{n+1} 
= \Gamma_{\hat G_n}^\omega (\chi_n) 
\eqno £7.4.4,$$
we denote by $b_{n+1}$ the block of $\chi_{n+1}\,,$ and  we choose a {\it Fitting pointed 
 $\O^*\-$group\/} $\hat O^{n+1}_{\,\eta_{n+1}}$ on $({\O_n}_*\hat G_n)b_n$ and an irreducible 
 character $\nu_{n+1}$ of~$\hat O^{n+1}$ involved in the irreducible character of 
 $N_{\hat G_{n+1}}(\hat O^{n+1}_{\,\eta_{n+1}})$  determined by $\chi_{n+1}$ (cf.~Corollary~£5.3). 
 Then, it is clear that  $\B = \{(\K_n,\hat G_n,\hat O^n_{\,\eta_n},\nu_n)\}_{n\in \Bbb N}$ 
 is a {\it charactered Fitting $\omega\-$sequence\/} of~$\hat G$ and that  $\{\chi_n\}_{n\in \Bbb N}$ becomes the 
 {\it irreducible character $\omega\-$sequence\/} associated  to~$\B\,.$

 \smallskip
Our construction only depends  on the choices of a {\it Fitting pointed $\O_n^*\-$ group\/} $\hat O^n_{\,\eta_n}$ 
on  $({\O_n}_*\hat G_n)b_n\,,$ and of an irreducible  character $\nu_{n}$ of~$\hat O^{n}$ involved in the 
 irreducible character of  $N_{\hat G_n}(\hat O^{n}_{\,\eta_{n}})$  determined by $\chi_{n}\,,$ for a 
 {\it finite  set of values of $n\,.$\/} Moreover, since all the  {\it Fitting pointed $\O_n^*\-$groups\/} on 
 $({\O_n}_*\hat G_n)b_n$ are mutually $G_n\-$conjugate, $\chi$ determines a unique $G\-$conjugacy
 class of {\it charactered Fitting $\omega\-$sequence\/} of~$\hat G\,.$
 That is to say, we have obtained a map
 $${\rm Irr}_\K (\hat G)\too {\rm Chs}^\omega_\K (\hat G)
 \eqno £7.4.5\phantom{.}$$
 and it is easily checked that this map is the inverse of the map~£7.4.3. We are done.

\bigskip
\bigskip
\noindent
{\bf £8. The inductive weights}
\bigskip
£8.1.  Let $\hat G$ be an $\O^*$-{\it group\/} with a finite $p\-$solvable $\cal O^*\-${\it quotient\/} $G\,,$ 
 assume that $\K$ is a splitting field for all the $\O^*\-$subgroups of $\hat G\,,$ choose a {\it polarization\/} 
 $\omega\,,$ and consider a {\it charactered Fitting $\omega\-$sequence\/} of~$\hat G$
$$\B = \{(\K_n,\hat G_n,\hat O^n_{\,\eta_n},\nu_n)\}_{n\in \Bbb N}
\eqno £8.1.1;$$
we call {\it charactered weight $\omega\-$sequence\/} associated to $\B$ any sequence 
$$\W = \{(\hat Q^n,\mu_n,\zeta_n)\}_{n\in \Bbb N}
\eqno £8.1.2\phantom{.}$$
where $(\hat Q^n,\mu_n,\zeta_n)$ is a {\it charactered weight\/} of $\hat G_n$ such that, denoting by $b_n$ 
and by $(\hat Q^n_{\,\delta_n},\mu_n,\bar\zeta_n)$ the block and a
{\it charactered $b_n\-$weight\/} corresponding to~$(\hat Q^n,\mu_n,\zeta_n)$  (cf.~£6.14), we have
$$\hat O^n_{\,\eta_n}\i \hat Q^n_{\,\delta_n}\qq {\rm Res}_{\hat O^n}^{\hat Q^n}(\mu_n) 
= {\mu_n(1)\over \nu_n(1)}\.\nu_n
 \eqno £8.1.3\phantom{.}$$
 and the corresponding {\it charactered weight\/} $(\skew3\hat{\bar Q}^n,\bar\mu_n,\skew3\bar{\bar\zeta}_n)$
 (cf.~Proposition~£6.16) coincides with $(\hat Q^{n+1},\mu_{n+1},\zeta_{n+1})$ for any $n\in \Bbb N$

 \medskip
 £8.2. It is clear that, for any $h\in \Bbb N\,,$ the sequence
 $$\W_h = \{(\hat Q^{h+n},\mu_{h+n},\zeta_{h+n})\}_{n\in \Bbb N}
\eqno £8.2.1\phantom{.}$$
is a {\it charactered weight $\omega\-$sequence\/} associated to $\B_h\,.$ Moreover, calling 
{\it charactered weight $\omega\-$sequence\/} of $\hat G$ any sequence $\W = \{(\hat Q^n,\mu_n,\zeta_n)\}_{n\in \Bbb N}$
as-sociated to some {\it charactered Fitting $\omega\-$sequence\/} $\B$ of~$\hat G\,,$ as in £7.2 above
 the group of inner automorphisms of $G$ acts on this set; let us denote by ${\rm Cws}^\omega_\K (\hat G)$
 the corresponding set of $G\-$conjugacy classes.

 \bigskip
 \noindent
 {\bf Theorem~£8.3.} {\it With the notation and the choice above, any charactered Fitting $\omega\-$sequence 
 $\B = \{(\K_n,\hat G_n,\hat O^n_{\,\eta_n},\nu_n)\}_{n\in \Bbb N}$ of~$\hat G$ admits an associated 
 charactered weight  $\omega\-$sequence  $\W = \{(\hat Q^n,\mu_n,\zeta_n)\}_{n\in \Bbb N}\,,$ 
 unique up to identifications,
and then we have ${\rm d}(\B) = {\rm d}(\mu_0)\,.$ Moreover, the correspondence
 sending $\B$ to $\W$ induces a natural bijection
 $${\rm Chs}^\omega_\K (\hat G)\cong {\rm Cws}^\omega_\K (\hat G)
 \eqno £8.3.1.$$ \/}
 \eject

\par
\noindent
{\bf Proof:} Since the sequence $\{\hat G_n\}_{n\in \Bbb N}$ ``stabilizes'', we can argue by induction 
on the minimal $n\in \Bbb N$ fulfilling $\hat G_{n+1} = \hat G_n\,.$ If $n = 0$ then $\hat O^0 =\O^*\,,$ 
$\nu_0 = {\rm id}_{\O^*}$ and the block $b_0$ of $\hat G_0$ has {\it defect zero\/}, so that $b_0$ only admits
the trivial {\it charactered $b_0\-$weight\/} $(\O^*_{\{b_0\}},{\rm id}_{\O^*},{\rm id}_{\O^*})$ 
(cf.~Proposition~£3.9) and, setting 
$$(\hat Q^n,\mu_n,\zeta_n) = (\O^*,{\rm id}_{\O^*},{\rm id}_{\O^*})
\eqno £8.3.2\phantom{.}$$
 for any $n\in \Bbb N\,,$ we get  a {\it charactered weight $\omega\-$sequence\/}
associated to $\B\,;$ moreover, it is clear that ${\rm d}(\B) = 0 = {\rm d}({\rm id}_{\O^*})\,.$

\smallskip
If $n\not= 0$ then, according to our induction hypothesis, the {\it charactered  Fitting $\omega\-$sequence\/} 
 $\B_1 = \{(\K_{1+n},\hat G_{1 +n},\hat O^{1+n}_{\,\eta_{1+N}},\nu_{1+ n})\}_{n\in \Bbb N}$ of~$\hat G_1$ 
 already admits an associated {\it charactered weight   $\omega\-$sequence\/}  
  $$\W_1 = \{(\hat Q^{1+n},\mu_{1+n},\zeta_{1+n})\}_{n\in \Bbb N}
\eqno £8.3.3\phantom{.}$$
fulfilling ${\rm d}(\B_1) = {\rm d}(\mu_1)\,;$ but, it follows from Proposition~£16.6 that there is 
a {\it charactered $b_0\-$weight\/} $(\hat Q^0_{\,\delta_0},\mu_0,\bar\zeta_0)$  of~$\hat G$
fulfilling condition~£8.1.3 and determining a {\it charactered weight\/}
 $(\skew3\hat{\bar Q}^0,\bar\mu_0,\skew3\bar{\bar\zeta}_0)$ of 
 $\hat F_{\O_*\hat G}(\hat O^0_{\,\eta_0},\nu_0) = \hat G_1$ which coincides with 
 $(\hat Q^1,\mu_1,\zeta_1)\,;$ then, it suffices to choose the {\it charactered weight\/} 
 $(\hat Q^0,\mu_0,\zeta_0)$ of $\hat G$ determined by $(\hat Q^0_{\,\delta_0},\mu_0,\bar\zeta_0)$
 (cf.~£6.14) to get a charactered weight  $\omega\-$sequence $\W$ associated to~$\B\,;$ moreover, we have (cf.~equalities~£6.4.2 and~Proposition~£6.16)
 $${\rm d}(\B) = {\rm d}(\nu_0) + {\rm d}(\B_1) = {\rm d}(\nu_0) + {\rm d}(\mu_1) = {\rm d}(\mu_0)
 \eqno £8.3.4.$$

\medskip
On the other hand, since the maps $\Delta_{(\hat O^n,\nu_n)}^\omega$ are bijective, it is quite clear
 that a {\it  charactered weight $\omega\-$sequence $\W = \{(\hat Q^n,\mu_n,\zeta_n)\}_{n\in \Bbb N}$\/} 
associated  to~$\B$ is, up to identifications,  uniquely determined by one of their terms; but, for~$n$ big enough, 
we know that $b_n$ is a block of {\it defect zero\/} of $\hat G_n$ and then $(\hat Q^n,\mu_n,\zeta_n)$ is
 uniquely determined (cf.~£8.3.2);  consequently, $\W$ is, up to identifications, uniquely determined by~$\B\,.$ 
 Moreover, this correspondence is clearly compatible with $G\-$conjugation and therefore it defines a 
 {\it natural\/} map
$${\rm Chs}^\omega_\K(\hat G)\too {\rm Cws}^\omega_\K (\hat G)
\eqno £8.3.5.$$

\smallskip
Conversely, the elements of ${\rm Cws}_\K (\hat G)$ are the $G\-$conjugacy classes of 
{\it charactered weight $\omega\-$sequences\/} $\W = \{(\hat Q^n,\mu_n,\zeta_n)\}_{n\in \Bbb N}\,,$ 
namely associated to some {\it charactered Fitting 
$\omega\-$sequence\/} $\B  = \{(\K_n,\hat G_n,\hat O^n_{\,\eta_n},\nu_n)\}_{n\in \Bbb N}$ of~$\hat G\,,$ 
so that the map~£8.3.5 sends the $G\-$conjugacy class of $\B$ to the $G\-$conjugacy class of $\W\,;$
moreover, condition~£8.1.3 shows that $\W$ determines $\B\,.$ We are done.

\medskip
£8.4. Finally, we inductively define the announced  {\it inductive weights.\/} We call {\it inductive weight\/} 
of $\hat G$ any {\it charactered weight\/} $(\hat Q,\mu,\zeta)$  of $\hat G$ which is either {\it trivial\/} 
or, denoting by $b$ the block determined by $\hat Q_\delta$ and by $(\hat Q_{\delta},\mu,\bar\zeta)$ 
some {\it charactered $b\-$weight\/} associated to $(\hat Q,\mu,\zeta)$ (cf.~£6.14), $\hat Q_\delta$ contains 
a Fitting pointed $\O^*\-$group $\hat O_\eta$ of $\hat G\,,$ we have
$${\rm Res}_{\hat O}^{\hat Q}(\mu)= {\mu(1)\over\nu(1)}\.\nu
\eqno £8.4.1\phantom{.}$$
for a suitable irreducible character of $\hat O\,,$  and the corresponding {\it charactered weight\/} 
$(\skew3\hat{\bar Q}^n,\bar\mu_n,\skew3\bar{\bar\zeta}_n)$ is an {\it inductive
weight\/} of $\hat F_{\O_*\hat G}(\hat O_\eta,\nu)$ (cf.~Proposition~£6.16). Let us denote by 
${\rm Inw}^\omega_\K (\hat G)$ the set of $G\-$conjugacy classes of  {\it inductive weights\/} 
of $\hat G\,.$

\bigskip
\noindent
{\bf Proposition~£8.5.} {\it A {\it charactered weight\/} $(\hat Q,\mu,\zeta)$  of $\hat G$ is an inductive weight
if and only if there is a charactered weight $\omega\-$sequence 
$$\W = \{(\hat Q^n,\mu_n,\zeta_n)\}_{n\in \Bbb N}
\eqno £8.5.1\phantom{.}$$ 
associated to a charactered Fitting $\omega\-$sequence 
$\B  = \{(\K_n,\hat G_n,\hat O^n_{\,\eta_n},\nu_n)\}_{n\in \Bbb N}$ 
of~$\hat G$ such that $(\hat Q,\mu,\zeta) = (\hat Q^0,\mu_0,\zeta_0)\,.$\/}

\medskip
\noindent
{\bf Proof:} We argue by induction on $\vert G\vert$ and may assume that $(\hat Q,\mu,\zeta)$ is not
trivial. The existence of $\W$ and $\B$ yields the existence of the Fitting pointed $\O^*\-$group 
$\hat O^0_{\,\eta_0}$ contained in $\hat Q^0_{\,\delta_0} = \hat Q_{\delta_0}$ and of the 
irreducible character $\nu_0$ of~$\hat O^0$
fulfilling equality~£8.4.1 for $\mu_0 = \mu\,;$ moreover, the {\it charactered weight\/} 
$(\skew3\hat{\bar Q},\bar\mu,\skew3\bar{\bar\zeta})$ coincides with $(\hat Q^1,\mu_1,\zeta_1)\,;$
then, since the charactered weight $\omega\-$sequence 
$$\W_1 = \{(\hat Q^{1+n},\mu_{1+n},\zeta_{1+n})\}_{n\in \Bbb N}
\eqno £8.5.2\phantom{.}$$ 
is associated to $\B_1\,,$ our induction hypothesis guarantees that 
$(\skew3\hat{\bar Q},\bar\mu,\skew3\bar{\bar\zeta})$ is already an {\it inductive weight\/}
 of $\hat F_{\O_*\hat G}(\hat O^0_{\,\eta_0},\nu_0)$ and therefore $(\hat Q,\mu,\zeta)$ is an  
 {\it inductive weight\/} too.

 \smallskip
 Conversely, if $(\skew3\hat{\bar Q},\bar\mu,\skew3\bar{\bar\zeta})$ is an {\it inductive
weight\/} of $\hat F_{\O_*\hat G}(\hat O_\eta,\nu)$ for suitable $\hat O_\eta$ and $\nu$ fulfilling the 
conditions above, it follows from our induction hypothesis that there is a charactered weight $\omega\-$sequence 
$$\W_1 = \{(\hat Q^{1+n},\mu_{1+n},\zeta_{1+n})\}_{n\in \Bbb N}
\eqno £8.5.3\phantom{.}$$ 
associated to a charactered Fitting $\omega\-$sequence 
$$\B_1  = \{(\K_{1+n},\hat G_{1+n},\hat O^{1+n}_{\,\eta_{1+n}},\nu_{1+n})\}_{n\in \Bbb N}
\eqno £8.5.4\phantom{.}$$ 
of~$\hat F_{\O_*\hat G}(\hat O_\eta,\nu)$ such that 
$(\skew3\hat{\bar Q},\bar\mu,\skew3\bar{\bar\zeta}) = (\hat Q^1,\mu_1,\zeta_1)\,,$ and we may
assume that $\K_1$ is the field extension of $\K$ of degree $\vert G\vert\,.$ Then, setting $\K_0 = \K\,,$
$\hat G_0 = \hat G\,,$ $\hat O^0 = \hat O$ and $\nu_0 = \nu\,,$ it is clear that 
$$\B  = \{(\K_n,\hat G_n,\hat O^n_{\,\eta_n},\nu_n)\}_{n\in \Bbb N}
\eqno £8.5.5\phantom{.}$$
\eject
\noindent
 is a  charactered Fitting $\omega\-$sequence of $\hat G\,,$ and setting $\hat Q^0 = \hat Q\,,$
 $\mu_0 = \mu$ and $\zeta^0 = \zeta\,,$ the sequence $\W = \{(\hat Q^n,\mu_n,\zeta_n)\}_{n\in \Bbb N}$
 is a  charactered weight $\omega\-$sequence associated to $\B\,.$ We are done.

 \bigskip
 \noindent
{\bf Corollary~£8.6.}  {\it  With the notation and the choice above, the correspondence mapping
any charactered Fitting $\omega\-$sequence  $\B = \{(\K_n,\hat G_n,\hat O^n_{\,\eta_n},\nu_n)\}_{n\in \Bbb N}$ 
of~$\hat G$ on the inductive weight $(\hat Q^0,\mu_0,\zeta_0)\,,$ where 
$\W = \{(\hat Q^n,\mu_n,\zeta_n)\}_{n\in \Bbb N}$ is the charactered weight  $\omega\-$sequence
of $\hat G$  associated to $\B\,,$ induces a natural bijection
$${\rm Chs}^\omega_\K (\hat G)\cong {\rm Inw}^\omega_\K (\hat G)
\eqno£8.6.1.$$
In particular, we have ${\rm d}(\B) = {\rm d}(\mu_0)\,.$\/}

\medskip
\noindent
{\bf Proof:} It follows from Theorem~£8.3 that there is such a correspondence mapping 
$\B$ on $(\hat Q^0,\mu_0,\zeta_0)$ and that it fulfills ${\rm d}(\B) = {\rm d}(\mu_0)\,;$ moreover, it is
clear that this correspondence is compatible with the $G\-$conjugation, and therefore we get a {\it natural\/} map
$${\rm Chs}^\omega_\K (\hat G)\too {\rm Inw}^\omega_\K (\hat G)
\eqno £8.6.2;$$
furthermore, since $(\hat Q^0,\mu_0,\zeta_0)$ clearly determines $\W\,,$ this map is injective.

\smallskip
In order to prove the surjectivity, we argue by induction on $\vert G\vert\,;$ let $(\hat Q,\mu,\zeta)$ be
an {\it inductive weight\/} that we may assume not trivial;  denoting by $b$ and by 
$(\hat Q_{\delta},\mu,\bar\zeta)$ the block and a {\it charactered $b\-$weight\/} corresponding to $(\hat Q,\mu,\zeta)$ (cf.~£6.14), it follows from the very definition that $\hat Q_\delta$ contains 
a Fitting pointed $\O^*\-$group $\hat O_\eta$ of $\hat G\,,$ that we have
$${\rm Res}_{\hat O}^{\hat Q}(\mu)= {\mu(1)\over\nu(1)}\.\nu
\eqno £8.6.3\phantom{.}$$
for a suitable irreducible character of $\hat O\,,$  and that the corresponding {\it charactered weight\/} 
$(\skew3\hat{\bar Q},\bar\mu,\skew3\bar{\bar\zeta})$ is an {\it inductive
weight\/} of $\hat F_{\O_*\hat G}(\hat O_\eta,\nu)\,.$

\smallskip
Thus, by our induction hypothesis,
there exists a charactered Fitting $\omega\-$sequence  
$$\B_1 = \{(\K_{1+n},\hat G_{1+n},\hat O^{1+n}_{\,\eta_{1+n}},\nu_{1+n})\}_{n\in \Bbb N}
\eqno £8.6.4\phantom{.}$$ 
of~$\hat G_1 = \hat F_{\O_*\hat G}(\hat O_\eta,\nu)$ such that the unique charactered weight  $\omega\-$sequence
$\W_1 = \{(\hat Q^{1+n},\mu_{1+n},\zeta_{1+n})\}_{n\in \Bbb N}$ of $\hat G_1$  associated to $\B_1$
fulfills 
$$(\hat Q^1,\mu_1,\zeta_1) = (\skew3\hat{\bar Q},\bar\mu,\skew3\bar{\bar\zeta})
\eqno £8.6.5;$$
moreover,  we may assume that $\K_1$ is the field extension of $\K$ of degree $\vert G\vert\,.$
Once again, setting $\K_0 = \K\,,$
$\hat G_0 = \hat G\,,$ $\hat O^0 = \hat O$ and $\nu_0 = \nu\,,$ it is clear that 
$$\B  = \{(\K_n,\hat G_n,\hat O^n_{\,\eta_n},\nu_n)\}_{n\in \Bbb N}
\eqno £8.6.6\phantom{.}$$
 is a  charactered Fitting $\omega\-$sequence of $\hat G\,,$ and setting $\hat Q^0 = \hat Q\,,$
 $\mu_0 = \mu$ and $\zeta^0 = \zeta\,,$ the sequence $\W = \{(\hat Q^n,\mu_n,\zeta_n)\}_{n\in \Bbb N}$
 is a  charactered weight $\omega\-$sequence associated to $\B\,;$ thus, the map~£8.6.2 sends $\B$
 to $(\hat Q,\mu,\zeta)\,.$ We are done.

 \bigskip
 \noindent
{\bf Corollary~£8.7.} {\it With the notation and the choice above, we have a natural bijection
$${\rm Irr}_\K (\hat G)\cong {\rm Inw}_\K (\hat G)
\eqno £8.7.1\phantom{.}$$
preserving the defects, which maps $\chi\in {\rm Irr}_\K (\hat G)$ on the $G\-$conjugacy class of
an inductive weight $(\hat Q,\mu,\zeta)$ provided there is a   charactered Fitting $\omega\-$se-quence
$\B$ of $\hat G$ such that the irreducible character   $\omega\-$sequence   $\{\chi_n\}_{n\in \Bbb N}$
and the charactered weight $\omega\-$sequence $\W = \{(\hat Q^n,\mu_n,\zeta_n)\}_{n\in \Bbb N}$ 
associated to~$\B$ fulfill 
$$\chi_0 = \chi\quad and\quad (\hat Q^0,\mu_0,\zeta_0) = (\hat Q,\mu,\zeta)
\eqno £8.7.2.$$  \/}

\medskip
\noindent
{\bf Proof:} It suffices to compose the inverse of isomorphism~£7.4.1 with isomorphism~£8.6.1.

\bigskip\bigskip
\centerline{\large References}

\bigskip\noindent
[1]\phantom{.} Jon Alperin,  {\it Weights for finite groups\/}, in Proc. Symp.
Pure Math. 47 (1987) 369-379, Amer. Math. Soc., Providence.
\smallskip\noindent
[2]\phantom{.} Michel Brou\'e and Llu\'\i s Puig, {\it Characters and Local
Structure in $G\-$al-gebras,\/} Journal of Algebra, 63(1980), 306-317.
\smallskip\noindent
[3] Henri Cartan and Samuel Eilenberg, {\it ``Homological Algebra''\/},
Princeton Math. 19, 1956, Princeton University Press.
\smallskip\noindent
[4]\phantom{.} Everett Dade,  {\it Counting Characters in Blocks\/}, Invent. Math., 109(1992), 198-272. 
\smallskip\noindent
[5]\phantom{.} Daniel Gorenstein, {\it ``Finite groups''\/} Harper's Series,
1968, Harper and Row.
\smallskip\noindent
[6]\phantom{.} Burkhard K\"ulshammer and Llu\'\i s Puig, {\it Extensions of
nilpotent blocks}, Inventiones math., 102(1990), 17-71.
\smallskip\noindent
[7]\phantom{.} Llu\'\i s Puig, {\it Pointed groups and  construction of
characters}, Math. Z. 176 (1981), 265-292. 
\smallskip\noindent
[8]\phantom{.} Llu\'\i s Puig, {\it Local fusions in block source algebras\/},
Journal of Algebra, 104(1986), 358-369. 
\smallskip\noindent
[9]\phantom{.} Llu\'\i s Puig, {\it Pointed groups and  construction of
modules}, Journal of Algebra, 116(1988), 7-129.
\smallskip\noindent
[10]\phantom{.} Llu\'\i s Puig, {\it Nilpotent blocks and their source
algebras}, Inventiones math., 93(1988), 77-116.
\smallskip\noindent
[11]\phantom{.} Llu\'\i s Puig, {\it Affirmative answer to a question of Feit},
Journal of Algebra, 131(1990), 513-526.
\smallskip
\noindent
[12]\phantom{.} Llu\'\i s Puig, {\it The Center of a Block\/} in {\it ``Finite Reductive Groups''\/},
Pro-gress in Math. 141(1997), 361-372, Birkh\"auser.
\smallskip\noindent
[13]\phantom{.} Llu\'\i s Puig, {\it Source algebras of $p\-$central group 
extensions\/}, Journal of Algebra, 235(2001), 359-398.
\smallskip\noindent
[14]\phantom{.} Llu\'\i s Puig, ``{\it Frobenius categories versus Brauer blocks\/}'',
Progress in Mahtematics, 274(2009), Bikh\"auser, Basel
\smallskip\noindent
[15]\phantom{.} Llu\'\i s Puig, {\it Block Source Algebras  in p-Solvable
Groups}, Michigan Math. Journal, 58(2009), 323-338
\smallskip\noindent
[16]\phantom{.} Llu\'\i s Puig, {\it Weight parameterization of simple modules for
 p-solvable groups\/}, to appear
 \smallskip\noindent
[17]\phantom{.} Geoffrey Robinson, {\it Dade's projective conjecture for $p\-$solvable groups\/}, Journal of Algebra,
229(2000), 234-248.

\bigskip
\bigskip
\noindent
{\bf Abstract.} The {\it weights\/} for a finite group $G$ with respect to a prime number $p$ where introduced by Jon Alperin, in order to formulate his celebrated conjecture. In 1992, Everett Dade formulates a refinement of Alperin's conjecture involving {\it ordinary\/} irreducible characters --- with their {\it defect\/} ---  and, in 2000, Geoffrey Robinson proves that the new conjecture holds for $p\-$solvable groups. But this refinement is formulated in terms of a vanishing alternating sum, {\it without\/} giving any possible refinement for the  {\it weights\/}. In this note we show that, {\it in the case of the $p\-$solvable finite groups\/}, 
the method developed in a previous paper  can be suitably refined to provide,  up to the choice of a {\it polarization $\omega$\/},   a {\it natural bijection\/} --- namely compatible with the action of the group of {\it outer automorphisms\/}  of $G$ --- between the sets of {\it absolutely irreducible characters\/} of $G$ and of $G\-$conjugacy classes of suitable 
{\it inductive weights\/}, preserving {\it blocks\/} and {\it defects\/}.

\end